\documentclass[11pt]{amsart}

\usepackage{lineno}
\usepackage{hyperref}
\usepackage[top=4cm, bottom=4cm, left=3cm, right=3cm]{geometry}
\usepackage{amsmath} 
\usepackage{amsthm} 
\usepackage{mathrsfs}
\usepackage{amsfonts}
\usepackage{mathtools}
\usepackage{amssymb}
\usepackage{tikz}
\usepackage{xcolor}
\usepackage{cite}
\usepackage{url}

\makeatletter
\newcommand{\pcite}[2][]{%
	(\kern-0.35em
	\@ifempty{#1}{\cite{#2}}{\cite[#1]{#2}}%
	)%
}
\makeatother

\usepackage[final]{showlabels} 

 \usepackage[disable]{todonotes} 

\newtheorem{thm}{Theorem}[section]
\newtheorem{corollary}[thm]{Corollary}
\newtheorem{proposition}[thm]{Proposition}
\newtheorem{theorem}[thm]{Theorem}
\newtheorem{lemma}[thm]{Lemma}

\theoremstyle{definition}

\newtheorem{example}[thm]{Example}
\theoremstyle{remark}
\newtheorem{remark}[thm]{Remark}
\numberwithin{equation}{section}

\newcommand{\auxVector}{h}
\newcommand{\auxHilbert}{\mathcal{H}}
\newcommand{\auxMult}{\mathcal{M}}
\newcommand{\RR}{\mathbb{{R}}}
\newcommand{\ZZ}{\mathbb{{Z}}}
\newcommand{\NN}{\mathbb{{N}}}
\newcommand{\TT}{\mathbb{{T}}}
\newcommand{\CC}{\mathbb{{C}}}
\newcommand{\DD}{\mathbb{{D}}}

\newcommand{\dom}{\operatorname{dom}}

\newcommand{\ran}{\operatorname{ran}}

\newcommand{\linearOp}{\mathcal{L}}

\newcommand{\proj}{P}
\newcommand{\generator}{E}
\newcommand{\embedding}{\mathcal{J}}

\renewcommand{\restriction}{\mathord{\upharpoonright}}

\newcommand\Real{{\mathfrak R}{\mathfrak e}\,} %
\newcommand\Imag{{\mathfrak I}{\mathfrak m}\,} %

\usepackage{hyperref}
\hypersetup{pdfpagemode={UseOutlines},	bookmarksopen, pdfstartview={FitH},	colorlinks,	linkcolor={red},	citecolor={blue},	urlcolor={blue}}

\makeatletter
\newcommand{\setword}[2]{%
	\phantomsection
	#1\def\@currentlabel{\unexpanded{#1}}\label{#2}%
}
\makeatother

\begin{document}
	\title[Similarity to contraction semigroups]{Similarity to contraction semigroups: \\ structural properties, criteria, \\ and applications to control theory}
	
	\author{J. Oliva-Maza}
	\address[J. Oliva-Maza]{Departamento de Matemáticas, Instituto Universitario de Matemáticas y Aplicaciones, Pedro Cerbuna 12, 50009 Zaragoza, Spain}
	\email{joliva@unizar.es}
	\thanks{The first author has been mainly supported by XXXV Scholarships for Postdoctoral Studies, by Ramón Areces Foundation. He has also been partially supported by Project PID2022-137294NB-I00, DGI-FEDER, of the MCYTS and Project E48-23R, D.G. Aragón, Universidad de Zaragoza, Spain.
	}
	
	\author{Y. Tomilov}
	\address[Y. Tomilov]{Institute of Mathematics, Polish Academy of Sciences, {\'S}niadeckich 8, 00-956 Warsaw, Poland, and 
		Faculty of Mathematics and Computer Science, Nicolas Copernicus University, Chopin Street 12/18, 87-100 Toru{\'n}, Poland}
	\email{ytomilov@impan.pl}
	\thanks{The second author was supported by the NCN Opus grant UMO-2023/49/B/ST1/01961. He was also partially supported by the NAWA/NSF grant BPN/NSF/2023/1/00001 and the NCN Weave-Unisono grant
		2024/06/Y/ST1/00044.}

	
	

\begin{abstract}
We reveal new aspects of the structure of Hilbert space $C_0$-semigroups $\mathcal T = (T(t))_{t\ge 0}$
similar to semigroups of contractions. In particular, we prove that $\mathcal T$ is similar to a semigroup
of contractions if and only if $\mathcal T$ is similar to a quasi-contraction $C_0$-semigroup and $T(t)$
is similar to a contraction for a single $t>0.$ Moreover, our methods allow us to estimate
the corresponding similarity constants and clarify their role in the study of similarity to
contractions. Along the way, we obtain similarity conditions
involving unbounded operators
and imposing minimal assumptions on regularity of $\mathcal T$. Such a general setting allows us to find significant applications to control theory,
including characterizations of exactly observable and exactly controllable systems. Finally we
establish several criteria of the same flavor for similarity to isometric semigroups, and illustrate
the developed theory by a number of pertinent examples.
\end{abstract}

\maketitle

\tableofcontents

\section{Introduction}\label{introduction}

\subsection{State of the art, motivations and illustrations}\label{state}
The paper presents a novel and systematic approach to the study of strongly continuous ($C_0$-)  operator semigroups on a Hilbert space $H$
similar to contraction semigroups, denoted for ease of reference by
$\mathcal{SC}(H).$

To motivate our results,  put them into the right context,
and underline their specific issues,
we first briefly recall some milestones of the discrete theory.
Contractions constitute one of the most well-understood
classes of bounded linear operators $\linearOp(H)$ on $H$,
with a number of useful properties
and rich theory. So, describing  operators similar to contractions
is one of the main tasks of operator theory.
The study of such operators
goes
back to Sz.-Nagy and Rota, who noted correspondingly
that  if $T\in \linearOp(H)$ is a power bounded compact operator on a Hilbert space $H$ or the spectral radius of $T$ is less
than $1$, then $T$ is similar to a contraction.
However, as was shown by Foguel in \cite{foguel1964counterexample}, there exist a power bounded operator on $H$ not similar to a contraction.
Moreover, Foguel's example turned out not to be not polynomially bounded, \cite{lebow1968power-bounded}.
So, after elaborating this example to a more elegant shape \pcite{halmos1964foguel},
Halmos inquired in \cite{halmos1970ten} whether polynomial boundedness instead of power boundedness
might be the right analytic description of similarity to contractions.
His question
spurred  intensive activity around similarity problems,
mainly concentrated on the study of Foguel-Hankel type operators,
functional calculus, and related norm-estimates, see e.g. \cite{aleksandrov1996hankel}, \cite{badea2001schur},
\cite{bourgain1986similarity},  \cite{bozejko1987littlewood} and  \cite{peller1982estimates}.

A high point of these developments was Pisier's landmark example
 of a polynomially bounded operator not
similar to a contraction,
see \cite{pisier2001similarity} for a nice account, and \cite{davidson1997polynomially} for related examples.
A number of  remarkable similarity criteria
were obtained on the way.
Among the criteria crucial for this paper, one can mention a theorem by Holbrook \pcite{holbrook1973operators} saying that
for $A, B \in \linearOp(H)$ the ``quadratic nearness''
\begin{equation}\label{nearness}
	\sum_{n=0}^{\infty}\|T^n - AC^nB\|^2< \infty,
\end{equation}
of $T$ to a contraction $C$ implies
the similarity of $T$ to a (possibly different) contraction.
The important partial case of that nearness,  $T^n =A C^n B$ for $A, B \in \linearOp(H),$ a contraction $C$ and $n \in \mathbb N\cup \{0\},$
was proved in \cite{holbrook1971spectral}.
Holbrook's results were further elaborated in e.g. \cite{badea2003operators}, \cite{hadwin1981dilations},
\cite {lebow1975spatial} and \cite{pisier2001similarity}.

A basic result of the theory due to
Paulsen \cite{paulsen1984every} asserts that a bounded operator $T$ on $H$ is similar to a contraction on $H$ if and only if $T$ is completely
polynomially bounded, i.e.,
for all $N \in \mathbb  N$ and matrices $(p_{ij})_{1\le i,j \le N}$ of polynomials, one has
\begin{equation}\label{pauls}
	\|(p_{i,j}(T))_{1 \le i,j \le N} \|\le C \sup_{z\in \DD} \| (p_{i,j}(z))_{1 \le i,j \le N} \|,
\end{equation}
for an absolute constant $C>0$ where the norms in \eqref{pauls} stand for the operator norms
on $\linearOp(H^N)$ and $\linearOp(\mathbb C^N).$
This criterion appeared to be very useful, especially in abstract theory,
although it is often non-trivial to apply it in concrete situations, and an appropriate continuous analogue ceases to exist.
(See Section \ref{cbdSubsection} for a relevant discussion.)

Despite Pisier's negative result, the study of operators similar to contractions
and their sub-classes became a well-developed area,
and it continued to bear fruit, see e.g. \cite{badea2003operators},
\cite{cassier2005generalized},
\cite{delaubenfels1998similarity},
\cite{gamal2019examples}, \cite{muller2007quasisimilarity} and \cite{popescu2014similarity}
as samples.
It is instructive to observe  that
if $T \in \linearOp(H)$
satisfies the so-called Ritt's
resolvent
condition $\|(z-T)^{-1}\|\le C|z-1|^{-1}$ for $z$ outside the unit disc, then the polynomial boundedness of $T$
does imply its similarity to a contraction, \cite[Theorem 5.1]{lemerdy1998similarity}. However, there are Ritt operators which are not polynomially bounded, and thus not similar to a contraction, as for instance \cite[Proposition 5.2]{lemerdy1998similarity} shows.
Among other results linking similarity to contraction to resolvent estimates, polynomial and power boundedness, we note
\cite{bouabdillah2024polygonal},
\cite{cassier2005generalized},
\cite{delaubenfels1998similarity} and \cite{gamal2018sufficient}.
However, it is worth noting that polynomial boundedness is a rather strong assumption, and it is often difficult to verify in concrete applications.

There are also numerous conditions that ensure similarity to unitary and isometric operators, and go back to B. Sz.-Nagy and J. Dixmier. These conditions, however, tend to be quite specific to their restricted setting, and typical results can be found in e.g. \cite{cassier2005generalized},
\cite{kubrusly2020asymptotic}, \cite{petitcunot2008intrinsic}, and
\cite{popescu1992similarity}.

Naturally, there were also attempts to develop
a similar theory for
continuous one parameter ($C_0$-) semigroups.
Recall that $C_0$-semigroups govern well-posed abstract Cauchy problems,
and thus they are backbone of an abstract approach to PDE.
Semigroups of contractions form a distinguished class, fundamental in applications to PDEs.
The main reason is that these semigroups allow for a simple and handy characterization
via the Lumer-Phillips theorem, and thus
appear  frequently in applications.
Moreover, they possess a rich functional calculus, admit unitary dilations and
functional models,
and thus represent one of the most well-studied  classes of $C_0$-semigroups on a Hilbert space.
So, the problem of describing $C_0$-semigroups $\mathcal T=(T(t))_{t \ge 0}$ similar to semigroups of contractions,
i.e., $T(t)=R S(t)R^{-1}, t \ge 0,$ for a contraction semigroup $\mathcal S=(S(t))_{t \ge 0}$ and invertible $R \in \linearOp(H),$
is of primary importance, and it will be addressed in this paper thoroughly.
To simplify our presentation, the class of such semigroups will be denoted by $\mathcal{SC}(H)$ in the sequel.

Unfortunately, despite its importance for applications,
$\mathcal{SC}(H)$ is not well-understood
and  in contrast to the discrete case
there are very few conditions ensuring
similarity to a contraction semigroup.
At the same time, several pertinent examples were found illustrating the difficulty
of  the problem.
Relying on ideas of Foguel \cite{foguel1964counterexample}, Packel constructed in \cite{packel1969semigroup} an example
of a bounded $C_0$-semigroup $\mathcal T \not \in \mathcal{SC}(H),$
with the generator being rather implicit. In \cite{benchimol1977feedback},
Benchimol provided a similar construction of $\mathcal T \not \in \mathcal{SC}(H).$
However, he first constructed the generator of $\mathcal T$ as a bounded perturbation
of a generator of a contraction semigroup, thus making both the generator and the semigroup transparent and explicit.
Moreover, among other instructive examples, Chernoff showed in \cite[p. 254]{chernoff1976two} how to produce an exponentially stable $C_0$-semigroup
that does not belong to $\mathcal{SC}(H).$
This somewhat  discouraged research on $\mathcal{SC}(H),$ and much later even stronger counterexamples emerged.
In particular,  Le Merdy observed in  \cite[Theorem 1.1]{lemerdy1998similarity} that a similar example can be realised for exponentially stable holomorphic semigroup.
In addition, in contrast to the case of holomorphic semigroups,
one can find $\mathcal T \not \in \mathcal{SC}(H)$ whose negative generator is invertible and has bounded imaginary powers.
The strongest known example of $\mathcal T \not \in \mathcal{SC}(H)$  was again  found by Le Merdy in \cite{lemerdy2000bounded},
see \cite[Section 9.1]{haase2006functional} for further discussion.
In his example, $\mathcal T=(T(t))_{t \ge 0}$ is exponentially stable, consists of compact 
operators for $t >0,$
and is holomorphic of angle $\pi/2.$ (The latter was not stated in \cite{lemerdy2000bounded}, but can be easily deduced
using e.g. \cite[Lemma 9.1.2]{haase2006functional}.)
Very recently, among other improvements to Le Merdy's example, it was shown in \cite[Theorem 1.3]{oliva2025tensor}
that  one can construct an even  stronger example
by replacing exponential decay with nilpotency while preserving compactness, though necessarily omitting holomorphicity. However, holomorphicity is present in another example 
from \cite[Theorem 1.3]{oliva2025tensor}, where nilpotency is relaxed to quasi-nilpotency.

Among very few positive results, Le Merdy proved in \cite[Theorem 1.1]{lemerdy1998similarity} that a sectorially bounded
holomorphic semigroup belongs to $\mathcal{SC}(H)$ if and only if its negative generator
has bounded imaginary powers, or equivalently
admits a bounded $H^\infty$-calculus on an appropriate sector. This allowed one to formulate similarity criteria in terms of square
function estimates, \cite[Theorem 7.3.1]{haase2006functional}, \cite[Theorem 4.3]{lemerdy1998similarity}, and, a posteriori, invoking resolvent of the generator \cite{boyadzhiev1994logarithms} (where boundedness of  imaginary powers was described via integral resolvent bounds).
Moreover, aiming at an analogue of Nagy's theorem, it was shown in \cite{vu1998similarity} that every bounded, uniformly continuous and quasi-compact
semigroup in a Hilbert space is similar to a contraction semigroup. Furthermore, \cite[Theorem 4.3]{arendt2001functional} showed that the assumption of uniform continuity in \cite{vu1998similarity} can be relaxed to the weaker conditions of analyticity and quasi-contractivity.

The examples by Chernoff and Le Merdy, along with their improvements in \cite[Theorem 1.3]{oliva2025tensor}, demonstrate two key points. First, they show that continuous analogs of the classical criteria by Nagy and Rota fail in a dramatic way.
Second,
they clarify that our study concerns in fact the \emph{joint similarity}
of  the  semigroup $(T(t))_{t \ge 0}$ to a family of contractions on $H$, i.e. $\|RT(t)R^{-1}\|\le 1$ for an invertible $R \in \linearOp(H)$ and $t \ge 0$. This joint similarity does not, in general, follow
from the similarity of each individual $T(t),$ $t \ge 0,$ to contraction - unlike the discrete case of the semigroup $(T^n)_{n \ge 0}.$
As it will be made clear in Section \ref{jointSect} below,
the property $\mathcal T \in \mathcal {SC}(H)$ depends on the behavior of the similarity constants
$C(T(t))$
near zero. For an operator
$T \in \linearOp(H)$ its similarity constant is defined as
\[
C(T)=\inf\left\{\|R\|\|R^{-1}\|\, : \, R \in \linearOp(H) \mbox{ invertible with } \|RTR^{-1}\|\le 1 \right\}
\]
with $C(T)=\infty$ if $T$ is not similar to a contraction.
If $\mathcal T \in \mathcal {SC}(H)$ then the (joint) similarity constant $\mathcal C (\mathcal T)$ is defined similarly
as
\[
\mathcal C(\mathcal T):=\inf \left\{\|R\|\|R^{-1}\|\,: \,R \in \linearOp(H) \mbox{ invertible with } \|RT(t)R^{-1}\|\le 1 \text{ for all } t \ge 0\right\},
\]
setting $\mathcal C(\mathcal T) = \infty$ if $\mathcal T \notin \mathcal{SC}(H)$.
Thus, the continuous setting requires  similarity criteria of a new type that
take into account the behavior of $\mathcal T$ near zero.
While both the behavior of $\mathcal T$ near zero and near infinity are crucial for understanding the similarity properties of $\mathcal T$,
their distinct roles have not been emphasized in the literature.
In particular, the significance of small
$t$ remained unclear, and  the potential relevance of quadratic nearness in this context is still an open question.

To date, the structure of semigroups in $\mathcal{SC}(H)$ has remained largely unexplored
and revealing its underlying issues has been a long-standing challenge.
To better understand this structure, we employ the notion of quasi-contractivity and explore the class of semigroups on
$H$ that are similar to quasi-contraction semigroups, denoted by
$\mathcal{SQC}(H)$.
Recall that a $C_0$-semigroup $\mathcal T = (T(t))_{t \ge 0}$ on $H$ is said to be \textit{quasi-contractive} if there is $\lambda \in \mathbb R$
such that
\begin{align*}
	\|T(t)\|\le e^{\lambda t}, \qquad t \ge 0.
\end{align*}
Quasi-contraction semigroups originate from the classical work \cite{lumer1961dissipative} and frequently appear in applications. For instance, they play a crucial role in proving the well-posedness of concrete linear PDEs via the (shifted) Lumer-Phillips theorem. They also arise in various stability problems, such as those related to Chernoff's formula and its applications in numerical analysis,
as well as in the study of nonlinear PDEs and their linearizations.
To illustrate the concept of
quasi-contractivity,
observe that every $C_0$-semigroup $\mathcal T$ with a bounded generator $\generator$ is quasi-contractive
as it satisfies the bound $\|T(t)\|\le e^{\|\generator\|t}$ for all $t \ge 0.$
Furthermore,  since any bounded perturbation of the generator of a quasi-contraction semigroup results in another quasi-contraction semigroup, bounded perturbations of normal operators whose spectrum lies in some left half-plane generate quasi-contraction semigroups.
By a similar argument, one can show
that every $C_0$-group on $H$ is similar to a quasi-contraction semigroup (cf.  \cite[Theorem 7.2.8 and Corollary 7.2.9]{haase2006functional}).

Thus, quasi-contractivity is a significantly weaker requirement than contractivity and naturally arises in semigroup theory. On the other hand,
Chernoff's example mentioned above \pcite[p. 254]{chernoff1976two}
actually yields a semigroup not similar to a quasi-contraction semigroup. Furthermore,
by \cite[Theorem 7.3.13]{haase2006functional}
one can construct an analogous example where, in addition, negative generator is invertible and has
bounded imaginary powers.
Moreover, as we show in Section \ref{ExamplesSect},
the Packel and Le Merdy semigroups  do not
belong to $\mathcal{SQC}(H)$.
Even stronger examples of semigroups outside  $\mathcal{SQC}(H)$
can be found in \cite{oliva2025tensor}.
At the same time, Benchimol's semigroup is similar to a quasi-contraction semigroup
since its generator arises as  a bounded perturbation of the generator
of a contraction semigroup.
Thus the relations between $\mathcal{SC}(H)$ and $\mathcal{SQC}(H)$ appear to be rather intricate,
and the aim to clarify them underlines the specifics and novelty of our approach.
The paper links quasi-contractivity to the study of semigroups in  $\mathcal{SC}(H)$ and sheds new light on their structure,
revealing the different roles of small and large times and their interplay. 
It aims to emphasize these roles and to make the subject accessible to different audiences.

\subsection{Ideas and results}
To develop our approach, we start with a Holbrook-type characterization of $\mathcal{SC}(H)$, stating that a $C_0$-semigroup $\mathcal T$ belongs to $\mathcal{SC}(H)$ if and only if there exist a contraction $C_0$-semigroup $\mathcal S=(S(t))_{t \ge 0}$, a Hilbert space $K$, and bounded operators $A \in \linearOp(K,H)$ and $B \in \linearOp(H,K)$ such that
\begin{equation}\label{h}
	T(t)=A S(t) B, \qquad t \ge 0.
\end{equation}
The result can be proved by following the arguments from the proof of its discrete analogue, established in \cite{holbrook1971spectral}. See also \cite[Proposition 4.2]{pisier2001similarity} for a more general setting.
This characterization, together with its generalization in the spirit of \eqref{nearness}, is fundamental to our approach. Accordingly, we provide both results with new, more explicit proofs, placing them within the context of representations of (algebraic) semigroups. In particular, we make quotient norms explicit and thereby avoid the indirect argument requiring verification of the parallelogram law, as in \cite{holbrook1973operators}. In this way, given $\tau >0$, we infer that $(T(t))_{t \ge \tau}$ is jointly similar to contractions if \eqref{h} holds only for $t \ge \tau$. This observation becomes crucial in the next step of our program, where we describe $\mathcal T=(T(t))_{t \ge 0}$ with at least one $T(t)$ similar to a contraction.

Since a direct characterization of $\mathcal T \in \mathcal{SC}(H)$ is not available, we turn to representations of $\mathcal T$ of the form \eqref{h}. The preceding discussion makes it clear that, besides the similarity properties of $\mathcal T$ for large times, it is also necessary to address these properties for small times. Accordingly, we fix $\tau >0$ and consider the restrictions of \eqref{h} to $[0,\tau]$ and $[\tau,\infty)$. Somewhat surprisingly, \eqref{h} for $t$ near zero and for $t$ separated from zero can be characterized via similarities of $T(t)$ in an elegant way, see Theorems \ref{awayzero} and \ref{nearzero} below, which are of independent interest. Unfortunately, \eqref{h}, when formulated for the partition of $[0,\infty)$ into $[0,\tau]$ and $[\tau,\infty)$, in general requires two distinct pairs of intertwining operators $A$ and $B$, whereas to recover \eqref{h} (under suitable assumptions) one needs only a single pair. We therefore combine the two characterizations of \eqref{h} on $[0,\tau]$ and on $[\tau,\infty)$ to obtain \eqref{h} for all $t \ge 0$ with a single appropriate pair $A$ and $B$. A matrix trick shows that this is indeed possible, leading to a new characterization of $\mathcal T \in \mathcal{SC}(H).$

When considering $t$ away from zero, the joint similarity of $(T(t))_{t \ge \tau}$ with contractions can be analyzed along the lines of the discrete case of $T \in \linearOp(H)$, invoking intuition developed in the classical discrete setup. The following statement, proved in Section \ref{AwayzeroSect}, makes this precise.
\begin{theorem}\label{awayzero}
	Let $\mathcal T = (T(t))_{t\geq0}$ be a $C_0$-semigroup on a Hilbert space $H.$
	The following are equivalent.
	\begin{enumerate}
		\item [(i)] There exists $\tau >0$ such that $T(\tau)$ is similar to a contraction.
		\item [(ii)] For each $\tau>0$, the family $(T(t))_{t\geq \tau}$ is jointly similar to contractions.
		\item [(iii)] For each $\tau>0$, there exist a Hilbert space $K$, $B \in \linearOp(H, K),$  $A \in \linearOp(K, H)$ and a unitary $C_0$-group $\mathcal U = (U(t))_{t\in \RR}$ on $K$ such that
		\begin{equation}\label{unitEq_int}
			T(t) = A U(t)B, \qquad t\geq \tau.
		\end{equation}
		\item [(iv)] There exist $\tau>0$, a Hilbert space $K$, $B \in \linearOp(H, K)$, $A \in \linearOp(K, H)$, and a $C_0$-semigroup $\mathcal S = (S(t))_{t\geq0}$ in $\mathcal{SC}(K)$ such that
		$$T(t) = A S(t)B, \qquad t \geq \tau.
		$$
	\end{enumerate}
\end{theorem}
Recall that if $T\in \linearOp(H)$ is such that $T^n$ is similar to a contraction for some $n \in \mathbb N,$ then all of
$(T^n)_{n\in \NN}$ are jointly similar to contractions. Theorem \ref{awayzero} reveals an analogous rigidity for the similarity property for
continuous parameter semigroups. 
More precisely,  if $(T(t))_{t\geq0}$ is a $C_0$-semigroup on  $H,$
then $(T(t))_{t \ge 0}$ are similar to contractions either for all  $t>0$ or for no $t>0,$ see Corollary \ref{SimSetCor}. However, joint similarity can be lost as e.g. Chernoff's example shows.

Remark that if $\mathcal T=(T(t))_{t \ge 0}$ is such that $T(\tau)$ is bounded from below for
some $\tau>0$ (and thus for all $\tau>0$), then
the problem of similarity of $\mathcal T$ to a contraction semigroup  reduces essentially to the same problem
for a single operator. In this case,  $\mathcal T \in \mathcal {SC}(H)$ if and only if
there exists $\tau>0$ such that $T(\tau)$ is similar to a contraction (Theorem \ref{boundedBelowContractProp}).
As a consequence, we obtain a variant of Liapunov's theorem (Corollary \ref{LiapunovThCor})
for left invertible $C_0$-semigroups (so also for $C_0$-groups) $\mathcal T$ saying that
for any $a$ strictly larger than the exponential type $\omega_0(\mathcal T)$ of $\mathcal T$, there exists
an equivalent Hilbertian norm $\|\cdot \|_{\mathscr H}$ on $H$
such that $\|T(t)\|_{\mathscr H} \leq e^{at}$ for all $t \geq 0.$ This result improves, in particular, \cite[Corollary 7.2.5]{haase2006functional} and compare it with \cite[Theorem 4.1]{arendt2001functional}.

After clarifying joint similarity to contractions  for the family $(T(t))_{t \ge \tau},$
we turn to the study of the same issue for $(T(t))_{t \in [0,\tau)},$ thus concentrating on small $t$.
This situation has not been emphasized in the literature so far.
To deal with small $t$ we rely on the notion of quasi-contractivity,
introduced above. Quasi-contractivity appeared to be the right tool to
characterize local joint similarity to contractions via global
joint similarity to quasi-contractions, which is easier to prove.
The following statement proved in Section \ref{NearzeroSect} recast quasi-contractivity in the spirit of Holbrook-type relation \eqref{h}.
\begin{theorem}\label{nearzero}
	Let $\mathcal T = (T(t))_{t\geq0}$ be a $C_0$-semigroup on a Hilbert space $H$.
	The following are equivalent.
	\begin{enumerate}
		\item [(i)] $\mathcal T \in \mathcal{SQC}(H).$
		\item [(ii)] For each $\nu>0$, there exist a Hilbert space $K$,  $B \in \linearOp(H, K)$, $A \in \linearOp(K, H),$ and a contractive and nilpotent $C_0$-semigroup $\mathcal N = (N(t))_{t\geq0}$ on $K$ such that
		$$T(t) = A N(t)B, \qquad t\in [0,\nu].
		$$
		\item [(iii)] For each $\nu>0$, there exist a Hilbert space $K$, $B \in \linearOp(H,K)$, $A \in \linearOp(K, H)$, and a unitary $C_0$-group $\mathcal U = (U(t))_{t\in \RR}$ on $K$ such that
		$$T(t) = A U(t)B, \qquad t\in [0,\nu].
		$$		
		\item [(iv)] There exist $\nu>0$, a Hilbert space $K$,  $B \in \linearOp(H, K)$, $A \in \linearOp(K, H),$ and a $C_0$-semigroup $\mathcal S=(S(t))_{t\geq0}$ in $\mathcal{SQC}(K)$ such that
		$$T(t) = A S(t)B, \qquad t\in [0,\nu].
		$$
	\end{enumerate}
\end{theorem}
The proof of Theorem \ref{nearzero} is given in Section \ref{NearzeroSect}. In particular, it depends on Proposition \ref{OT24quasi}, a partial case of a new  characterization of tensor products of semigroups similar to semigroups of contractions obtained in \cite[Theorem 1.1]{oliva2025tensor}.
Note a Banach space version of Theorem \ref{nearzero} given in \cite[Theorem 5.4]{batty2017lower},
where Banach space semigroups bounded from below
were characterized in terms of existence of equivalent norm making their negative generators dissipative.

Theorem \ref{nearzero} allows us, in particular, to describe semigroups in $\mathcal{SQC}(H)$ 
in geometric terms: they are precisely those admitting a group dilation, as the following corollary shows.
\begin{corollary}\label{C0groupdilation_int}
	Let $\mathcal T = (T(t))_{t\geq0}$ be a $C_0$-semigroup on a Hilbert space $H$. Then the following are equivalent.
	\begin{itemize}
		\item [(i)] $\mathcal T \in \mathcal {SQC}(H)$.
		\item [(ii)] $\mathcal T$ has a $C_0$-group dilation, i.e., there exist a Hilbert space $K$, an isometry $V \in \linearOp(H, K)$ and a $C_0$-group $\mathcal G = (G(t))_{t\in \RR}$ on $K$ with
		$$
		T(t) = V^\ast G(t) V, \qquad t \geq0.
		$$
				\item [(iii)]
	There exist $\nu, \nu'>0$, a Hilbert space $K$,  $B \in \linearOp(H, K)$, $A \in \linearOp(K, H),$ and a $C_0$-semigroup $\mathcal S=(S(t))_{t\geq0}$ on $K$
	such that $S(\nu')$ is bounded from below and
	$$
		T(t) = A S(t)B, \qquad t\in [0,\nu].
		$$
			\end{itemize}
\end{corollary}
It is instructive to recall that the assumption of  boundedness from below of some $S(t)$ is equivalent to the left-invertibility of $S(t)$ for all $t\ge 0.$ Corollary \ref{C0groupdilation_int} should be compared to
\cite[Theorem 7.3]{batty2017lower}
characterizing extensions to $C_0$-groups. 

Finally, based on the preceding considerations, we combine the statements (and the corresponding arguments) for small and large $t$
and represent $\mathcal T$ as in \eqref{h},  thus showing that $\mathcal T \in \mathcal{SC}(H).$
This step requires
an additional, new argument based on  an auxiliary operator-theoretical construction.
The next statement, proved in Section \ref{mainSection}, is one of the main results of the paper.
\begin{theorem}\label{quasiContrTh_int}
	Let $H$ be a Hilbert space and let $ \mathcal T = (T(t))_{t\geq0}$ be a $C_0$-semigroup on $H.$
	Then $\mathcal T \in \mathcal{SC}(H)$ if and only if the following conditions hold.
	\begin{enumerate}
		\item [(i)] $\mathcal T \in \mathcal{SQC}$(H).
		\item [(ii)] There exists $\tau>0$ such that $T(\tau)$ is similar to a contraction on $H$.
	\end{enumerate}
\end{theorem}
{Thus, condition (i) addresses the similarity behavior of $\mathcal T$ near zero, 
whereas condition (ii) governs its similarity properties at infinity.} Moreover, in the setting of Theorem \ref{quasiContrTh_int}, our methods allow us to provide explicit estimates for the similarity constant of $\mathcal T$ in terms of the corresponding constants for $T(\tau)$ and for the quasi-contraction semigroup appearing in Theorem \ref{quasiContrTh_int}(i), see Theorem \ref{SimConstTh} below. With the exception of \cite{holbrook1977distortion}, the existing literature has not focused on such estimates, though they are of considerable importance. The importance of this is illustrated in Section \ref{jointSect}, where we prove that $\mathcal T$ is similar to a contraction semigroup if and only if $\liminf_{t\to 0} C(T(t)) < \infty$, showing that the small-time behavior of similarity constants determines the similarity properties of the entire semigroup $\mathcal T$. We also provide a parallel characterization in terms of the asymptotic behavior of the similarity constants of the resolvents $\lambda (\lambda - \generator)^{-1}$ as $\lambda \to \infty$, where $\generator$ denotes the generator of $\mathcal T$. As a corollary, we further establish a curious trichotomy (Corollary \ref{trichotomyCor}), classifying the possible similarity properties of $\mathcal T$ into three distinct cases, and where the asymptotic behavior of $C(T(t))$ as $t \to 0$ suffices to determine which of the three cases applies.

Recalling the discussion in Section \ref{state}, note that the semigroup $(T(t))_{t \ge 0}$ constructed by Benchimol is quasi-contractive, while $T(t)$ is not similar to a contraction for any $t > 0$. On the other hand, Le Merdy's semigroup does not belong to $\mathcal{SC}(H)$, although each operator $T(t)$ is similar to a contraction, see Section \ref{standardexamples} for further details. Hence, conditions (i) and (ii) in Theorem \ref{quasiContrTh_int} are, in general, independent of each other.

\subsection{Examples and applications}

The theory developed in the preceding sections can be illustrated by pertinent and explicit examples.
Some of them originate from \cite{oliva2025tensor}, while others worked out in this paper.
In particular, via a detailed analysis of Packel's semigroup and its variants,
we show in Section \ref{PackelSect} that generically they do not belong to $\mathcal{SC}(H),$ and that these variants may be constructed so that they fail to satisfy either condition (i) or (ii) in Theorem \ref{quasiContrTh_int}, or both simultaneously.
As a result we, among other things, simplify  constructions in \cite{chernoff1976two}, making some of them redundant.
At the same time, we exhibit  a new class of bounded semigroups in $\mathcal{SQC}(H) \setminus \mathcal{SC}(H)$
using the interpolation technique of Bhat-Skeide \cite{bhat2015pure}, and its extension in \cite{dahya2024interpolation}.

After clarifying the structure of $\mathcal{SC}(H)$ and $\mathcal{SQC}(H),$
we turn to characterizations of semigroups from these classes
in terms of bounds for the associated averages
arising in applications.
Similar bounds  appeared frequently in the studies of similarity of  $T \in \linearOp(H)$
to contractions and, more specifically, to isometries and unitaries. The proofs often relied on the fact
that if appropriate averages of $T$ are bounded away from  zero and infinity,
then they can be employed to define an equivalent inner product on $H$,
making $T$ a contraction. The framework of $C_0$-semigroups $\mathcal T$ suggests a new, more general form
of these conditions involving averages ``weighted'' with unbounded operators and thus
with important consequences for control theory and its applications to PDEs.

In particular,  invoking unbounded operator weights and using a Holbrook-type similarity condition,
we obtain  new characterizations of $\mathcal{SQC}(H)$ and $\mathcal{SC}(H).$
In turn, these characterizations allow us to clarify the notions of observability
and controllability of control systems by relating them to semigroups  in $\mathcal{SC}(H)$ and $\mathcal{SQC}(H)$.
The notions go back to Kalman for finite-dimensional spaces, and have become
a basic part of the theory due to the  influential book by  J.-L. Lions \cite{lions1988controlabilite}.
A nice account of the semigroup approach to the study of these notions can be found
in \cite{tucsnak2009observation}.

To give a flavor of our results, consider the control system
\begin{eqnarray}\label{observab}
	\dot x(t)&=&\generator x(t), \qquad t \ge 0, \qquad
	x(0)=x_0,\\
	y(t)&=&Cx(t),\notag
\end{eqnarray}
where $\generator$ is the generator of a $C_0$-semigroup $\mathcal T=(T(t))_{t \ge 0}$ on $H$,
and 
$C$ is a bounded operator from ${\rm dom}(\generator)$ (with the graph norm) into a Hilbert space $(K, \|\cdot\|_K)$ 
such that
\[
\|y\|_{L^2([0,\tau])}^2=\int_{0}^{\tau} \|CT(s)x_0\|_K^2 \, ds \le C_\tau \|x_0\|^2,
\]
for all $x_0 \in {\rm dom}(\generator)$ and some $C_\tau >0,$ i.e., $C$ is an admissible observation operator for $\generator$.
A basic problem in control theory
is to describe situations when the initial state $x_0$ of \eqref{observab} can be recovered in a continuous way
from the output signal $y$ on the time interval $[0, \tau],$
thus making \eqref{observab} observable.
One of the strong forms of observability for \eqref{observab}
was formalized as
\begin{equation}\label{observ_con}
	\int_{0}^{\tau} \|CT(s)x_0\|_K^2 \, ds \ge c_\tau \|x_0\|^2,
\end{equation}
for some $c_\tau>0,$
and called exact observability of \eqref{observab} in time $\tau>0$,
and infinite time exact observability if $\tau=\infty.$

To illustrate the notion of controllability, consider the dynamical system
\begin{eqnarray}\label{controlab}
	\dot x(t)&=&\generator x(t)+By(t), \qquad t \ge 0\\
	x(0)&=&x_0, \notag
\end{eqnarray}
where $B$ is a bounded operator from a Hilbert space $K$ to $H_{-1}.$ (Here, $H_{-1}$ denotes the completion of $H$ in the norm induced by the resolvent of $\generator.$) If the convolution of $\mathcal T$ with $By$ belongs to $H$ for all $y \in L^2([0,\tau], K)$, i.e.,
\[x= \int_{0}^\tau T(t-s) By(s)\, ds \in H,\]
then $B$ is called an admissible control operator for \eqref{controlab}.  Exact controllability of \eqref{controlab} in time $\tau > 0$ means that any $x \in H$ can be realized as above, or equivalently, setting $x_0=0,$ any final state $x \in H$ of \eqref{controlab} can be reached by an appropriate input $y \in L^2([0,\tau], K)$.
For an advanced theory of this and other common ways to look  at observability and controllability
see again  \cite{tucsnak2009observation}.

Developing our approach to the study of $\mathcal{SC}(H)$ and $\mathcal{SQC}(H)$, we prove that a control system $\eqref{observab}$ is \emph{finite time exactly observable}
(for some admissible $C$) if and only if the semigroup $\mathcal T$ generated by $\generator$
belongs to $\mathcal{SQC}(H).$ So $\mathcal{SQC}(H)$ separates a class of semigroups and their generators suitable  for the study of observability of \eqref{observab}.
Moreover, we show \eqref{observab} is \emph{infinite time exactly observable} for some admissible $C$ if and only if $\mathcal T$ is strongly stable and belongs to $\mathcal{SC}(H),$
thus covering \cite[Theorem 3.1 and Corollary 3.1]{grabowski1996admissible} and \cite[Proposition 5.1]{haak2012exact}. (See Section \ref{ObservSubsect} for more details.)  This remains to be true even if the observability property of \eqref{observab} is slightly relaxed,
and both, infinite time controllability and observability, imply similarity of the corresponding $C_0$-semigroup to a semigroup of contractions.
The case of bounded $C$ was addressed, in particular, in \cite{dolecki1977general} and \cite{zwart2013left},
where for such $C$ the exact observability of \eqref{observab} was characterized by the left-invertibility of $\mathcal T=(T(t))_{t \ge 0},$
the property stronger than being in $\mathcal{SQC}(H).$ However, in many of applications $C$ is required to be unbounded,
and may not even be closable, which emphasizes the naturality of our setting. By duality, similar results hold when observability is replaced by controllability.
For this and other relations to  control-theoretical results see Section \ref{DefectOperSect}.

Finally, we provide criteria for similarity to an isometric semigroup involving two-sided bounds for averages of $C\mathcal T$
and being in the same spirit as our results on observability mentioned above.
By replacing $\mathcal T$ with $C\mathcal T$ for unbounded $C$, these results substantially generalize
well-known similarity conditions due to van Casteren \cite{van1983operators}, Naboko \cite{naboko1984conditions} and Malamud \cite{malamud2003similarity},
and could also be crucial in control theory. The arguments employed in the proofs resemble those used in these papers. However, the technical details differ and had to be adjusted to a more demanding setting.

\section{Preliminaries: basic objects, results, and tools}\label{PreliminariesSect} 

In this section, we recall several basic definitions and results concerning operators similar to contractions, which are crucial for what follows. While most of this material is standard and some of it was already discussed in the introduction, we elaborate on it here for clarity and ease of reference. As the paper is directed at communities with different backgrounds, we recall several facts that might be well-known to some, but obscure to others.

Recall that a bounded operator $T$ on a Hilbert space $H$ is said to be \textit{similar to a contraction} on $H$ if there exists a bounded and invertible $R\in \linearOp(H)$ such that $\|RTR^{-1}\|\leq 1$, or equivalently, if there exists an equivalent Hilbertian norm $\|\cdot\|_{\mathscr H}$ on $H$ satisfying $\|T\|_{\mathscr H} \leq 1$.  Given an operator $T\in \linearOp(H)$  that is similar to a contraction, we define its {similarity constant} $C(T)$  by
\begin{align*}
	C(T) := \inf \left\{ \|R\| \|R^{-1}\| \, : \, R \in \linearOp(H) \mbox{ invertible and } \|RTR^{-1}\|\leq 1\right\}.
	\end{align*}
If $T$ is not similar to a contraction, we let $C(T) := \infty$. 
It is direct and well-known that $T$ is similar to a contraction if and only if there
exists an equivalent Hilbertian norm $\|\cdot\|_{\rm eq}$ on $H$ such that
 $\|Th \|_{\rm eq} \le \|h\|_{\rm eq}$  for all $h \in H.$
Moreover,
\begin{align*}
	C(T)=&\inf\{M \geq 1\, : \, \mbox{there exists a Hilbertian norm } \|\cdot\|_{\rm{eq}} \mbox{ on } H \mbox{ such that} \\
	& \qquad \|T h\|_{\rm{eq}} \leq \|h\|_{\rm{eq}} \mbox{ and } \|h\| \leq \|h\|_{\rm{eq}} \leq M \|h\| \mbox{ for all } h \in H\}.
\end{align*}
See \cite[Proposition 04]{pisier2001similarity}, \cite[p. 230]{holbrook1977distortion}, and 
\cite[pp. 235-236]{benchimol1977feedback} for these basic properties.
Note that, since $\|\cdot\|_{\rm eq}$ is an equivalent norm, we have for all $h\in H$,
\begin{equation}\label{eqq}
	\|Th\| \le \|Th\|_{\rm eq} \le \|h\|_{\rm eq} \le M \|h\|,
\end{equation}
and thus $\|T\|\le M$. Taking the infimum over such $M$ yields $C(T)\ge \|T\|$.

The similarity constants were studied in particular in \cite{holbrook1977distortion}, where they are called \textit{distortion coefficients}, and in \cite{badea2003operators}. Note that if $C(T)$ is finite, then it is attained, that is, there exists a Hilbertian norm $\|\cdot\|_{\rm{eq}}$ on $H$ such that $\|h\| \leq \|h\|_{\rm{eq}} \leq C(T) \|h\|$ and $\|T h \|_{\rm{eq}} \leq \|h\|_{\rm{eq}}$ for all $h \in H$, see \cite[Proposition 2.4]{holbrook1977distortion}.

Let now $\mathcal T = (T_\gamma)_{\gamma \in \Gamma}$ be a family of bounded operators on a Hilbert space $H$, where $\Gamma$ is an index set. Similarly to the above, we say that
$\mathcal T$ is \textit{jointly similar to contractions} on $H$
if there exists a bounded and invertible $R \in \linearOp(H)$ such that
\begin{equation}\label{SimEq}
	\|R T_\gamma R^{-1}\|\leq 1, \qquad \gamma \in \Gamma.
\end{equation}
The latter condition is clearly equivalent to the existence of an equivalent Hilbertian norm $\|\cdot\|_\mathscr H$ on $H$
such that $\|T_\gamma\|_\mathscr H \leq 1$ for all $\gamma \in \Gamma$.
Note that similarity to a contraction for each $T_\gamma$ does not, in general, imply the joint similarity of $(T_\gamma)_{\gamma\in \Gamma}$ to contractions; see, e.g., Section~\ref{ExamplesSect}.
This fact is one of the origins of the main difficulties in similarity theory for $\Gamma=\mathbb R_+$,
and is in contrast to the case $\Gamma=\mathbb Z_+$. 

If $\mathcal T$ is jointly similar to contractions, then, analogously to the discrete case, the \textit{similarity constant} $\mathcal C(\mathcal T)$ of $\mathcal T$ is defined by
\begin{align*}
	\mathcal C(\mathcal T) := \inf \left\{ \|R \| \|R^{-1}\| \, : \, R \in \linearOp(H) \mbox{ invertible and } \|RT_\gamma R^{-1}\|\leq 1 \mbox{ for all } \gamma \in \Gamma\right\}, 	
\end{align*}
and satisfies
\begin{equation}\label{eqNormSim}
	\begin{aligned}
	\mathcal C(\mathcal T) =& \inf\{M \geq 1 \, :  \, \mbox{there exists a Hilbertian norm } \|\cdot\|_{\rm{eq}} \mbox{ on } H \mbox{ such that }
	\\  & \qquad  \|T_\gamma h\|_{\rm{eq}} \leq \|h\|_{\rm{eq}},  \quad  \|h\| \leq \|h\|_{\rm{eq}} \leq M \|h\| \mbox{ for all } h \in H, \, \gamma\in \Gamma\}.
	\end{aligned}
\end{equation}
Moreover, using elementary bounds analogous to \eqref{eqq}, we have
\begin{equation}\label{const_simm}
	\mathcal C(\mathcal T)\ge \sup_{\gamma\in \Gamma}\|T_\gamma\|.
\end{equation}
If $\mathcal T$ is not jointly similar to contractions, then we set $\mathcal C(\mathcal T) = \infty$. Similarly to the case of a single operator, $\mathcal C(\mathcal T)$ is always attained. This can be inferred from the proof of \cite[Proposition 2.4]{holbrook1977distortion}.

We will be mainly interested in one-parameter families $\mathcal T = (T(t))_{t\geq0}\subset \mathcal L(H)$ that are $C_0$-semigroups. If the family $(T(t))_{t\geq0}$ is jointly similar to contractions, we say that $\mathcal T$ is \textit{similar to a semigroup of contractions},
where the latter is necessarily a $C_0$-semigroup. The class of
such $C_0$-semigroups $\mathcal T$  will be denoted by $\mathcal{SC}(H).$

One of the fundamental properties of contraction $C_0$-semigroups on Hilbert spaces is the existence of their unitary dilations, provided by the classical Sz.-Nagy dilation theorem; see, for example, \cite{nagy2010harmonic} or \cite[Corollary~6.14]{davies1980one}. To formulate this result, we introduce some notation that will be used from now on.  Given two Hilbert spaces $H$ and $K$ such that $H \subseteq K$, and given $T \in \linearOp(K)$, we denote by $\proj_H \in \linearOp(K)$  the orthogonal projection from $K$ onto $H$, and by $T\restriction_H \in \linearOp(H,K)$ the restriction of $T$ to $H$.

In this setting, if $\mathcal{T} = (T(t))_{t \geq 0}$ is a contraction $C_0$-semigroup on a Hilbert space $H$, then there exists a unitary dilation of $\mathcal T$, that is, there exist a Hilbert space $K$ containing $H$ as a subspace, and a unitary $C_0$-group $\mathcal{U} = (U(t))_{t \in \mathbb{R}}$ on $K$ such that
\begin{equation}\label{dil}
	T(t) = \proj_H U(t) \restriction_H, \qquad t \geq 0.
\end{equation}

There's a geometric description of dilations due to Sarason, which in our setting takes the next form. Suppose $\mathcal{T} = (T(t))_{t \geq 0}$ is a semigroup on a Hilbert space $K$, and $H \subseteq K$ is a subspace. Then the compression $\mathcal{T}_H = \proj_H \mathcal{T} \restriction_H$ is a semigroup on $H$ if and only if there exist $\mathcal T$-invariant subspaces $H_2$ and $H_1$ of $K$ such that the orthogonal difference $H_2 \ominus H_1$ is precisely $H$. See e.g. \cite[Theorem 1.7, Remark 1.8]{pisier2001similarity} for the proof of a more general result
and comments. Such a subspace $H$ is often called semi-invariant in the literature. Clearly, $\mathcal T_H$ is strongly continuous if $\mathcal T$ is so.

The property \eqref{dil} can be recast in a slightly more general form.
We will say that  $\mathcal S$ is a \emph{dilation} to $K$ of another semigroup $\mathcal T$ on a possibly different space $H$ if there exists an isometry $V \in \linearOp(H, K)$ such that
\begin{equation}\label{dil_isom}
	T(t) = V^*S(t)V, \qquad t \ge 0.
\end{equation}
A detailed formal argument explaining how \eqref{dil} fits into this definition can be found e.g. in
\cite[Section D.8]{eisner2015operator}.

There are very few criteria for similarity to semigroups of contractions on Hilbert spaces.
Among them is the next well-known result due to Sz.-Nagy
characterizing semigroups similar to semigroups of isometric (and unitary) operators,
and relevant for the sequel. Specifically, a $C_0$-semigroup $\mathcal{T} = (T(t))_{t \geq 0}$ on a Hilbert space $H$ is similar to a semigroup of isometries if and only if there exist constants $\alpha, \beta > 0$ such that
\begin{equation}\label{nagy_isom_in}
	\alpha \|h\| \le \|T(t)h\| \le \beta \|h\|, \qquad h \in H, \; t \ge 0.
\end{equation}
In addition, $\mathcal{T}$ is similar to a unitary group if and only if $\mathcal T$ is invertible and \eqref{nagy_isom_in} holds. This characterization highlights the importance of uniform boundedness from below in the study
of similarity to contractions, which will play an important role in this paper as well.

Our considerations will depend on a generalization of the notion of contractivity, called quasi-contractivity.
Recall that a $C_0$-semigroup $\mathcal T = (T(t))_{t\geq0}$ is said to be a \textit{quasi-contraction semigroup} if there exists $\lambda \in \RR$ with
\begin{equation}\label{quasi}
\|T(t)\| \leq e^{\lambda t}, \qquad t \geq 0.
\end{equation}
Note that  quasi-contractivity is, in fact, a local condition and addresses the behavior of $\mathcal T$ near zero. It is easy to show that \eqref{quasi} is equivalent to
\begin{equation}\label{quasi_loc}
\|T(t)\|=1 + O (t), \qquad t \to 0^+.
\end{equation}
In particular, if $\|T(t)\| \leq f(t)$ for $t\geq0$, where $f:[0,\infty) \to [0,\infty)$ is differentiable at $0$ and satisfies $f(0) =1$, then
\begin{equation}\label{quasi_loc_2}
	\|T(t)\| \leq e^{f'(0)t}, \qquad t \geq 0.
\end{equation} 

The class of $C_0$-semigroups on $H$ similar to a quasi-contraction semigroup will be denoted by $\mathcal{SQC}(H).$
Clearly, a semigroup $\mathcal T$
 is in $\mathcal{SQC}(H)$ if and only if there exists an equivalent Hilbertian norm
on $H$ making $\mathcal T$ a quasi-contraction semigroup.
The class $\mathcal {SQC}(H)$ is quite large, and it includes e.g. $C_0$-groups on $H$ as  noted in Section \ref{introduction}. Apart from the perturbation argument mentioned there, this follows also from the fact that the generator of a $C_0$-group on $H$ admits a bounded $H^\infty$-calculus in an appropriate strip. See \cite[Section 7.2]{haase2006functional} for more on that, or \cite[Section 8]{yakubovich2004linearly} for an alternative approach relevant for this paper.

Note that both contractivity and quasi-contractivity are highly unstable with respect to similarities. It is instructive to note that, as proved in \cite[Theorem 2]{matolcsi2003relation}, for every $C_0$-semigroup $\mathcal T=(T(t))_{t \ge 0}$ with unbounded generator on $H$ there is a semigroup on $H$ similar to $\mathcal T,$ which is not quasi-contractive.

We will be using frequently shift semigroups defined on vector-valued
$L^2$-spaces. Let $H$ be a Hilbert space, $\Delta\subset\mathbb{R}$ an interval (possibly unbounded), and $\chi_\Delta$ the characteristic function of $\Delta.$  
 For $\lambda\in \RR$ set
 \[
 e_\lambda(t):=e^{\lambda t},\qquad t\ge0,
 \]
 and define the weighted $L^2$-space of $H$-valued functions by
 \[
 L^2(\Delta,e_{-\lambda};H):=\Big\{ f:\Delta\to H\ \text{(Bochner) measurable} : 
 \int_\Delta \|f(x)\|^2\, e^{-2\lambda x}\,dx <\infty\Big\}.
 \]
 On this space we define the right- and left-shift semigroups 
 $\mathcal R=(R(t))_{t\ge0}$ and $\mathcal L=(L(t))_{t\ge0}$ by
 \[
 (R(t)f)(x) = f(x-t)\,\chi_\Delta(x-t),\qquad
 (L(t)f)(x) = f(x+t)\,\chi_\Delta(x+t),
 \]
 for $f\in L^2(\Delta,e_{-\lambda};H)$, $t\ge0$, and $x\in\Delta.$
 It is straightforward to check that both $\mathcal R$ and $\mathcal L$ are 
 $C_0$-semigroups on $L^2(\Delta,e_{-\lambda};H).$ Moreover, one has
 \[
 \|R(t)\|\le e^{-\lambda t},\qquad \|L(t)\|\le e^{\lambda t},\qquad t\ge 0,
 \]
 so in particular both semigroups are contractive when $\lambda=0.$
 If it is crucial to indicate the interval explicitly we use subscripts, 
 for instance $R_{\lambda,\nu}(t)$ or $L_{\lambda,\nu}(t)$ when $\Delta=[0,\nu].$ 
 If the underlying space is clear we simply write $\mathcal R$ and $\mathcal L.$

Shift semigroups are closely related to the so-called \emph{evolution semigroups}, which have proved to be a useful tool in semigroup theory and will play a crucial role in our arguments. Let $\Delta\subseteq\mathbb{R}$ be an interval and $H$ a Hilbert space. Define the multiplication operator $\auxMult_T\in\linearOp\big(L^2(\Delta,e_{-\lambda}; H)\big)$ by
\[
(\auxMult_T f)(x)=T f(x),\qquad x\in\Delta,\; f\in L^2(\Delta, e_{-\lambda}, H),
\]
and note that $\|\auxMult_T\|=\|T\|$. If $\mathcal{T}=(T(t))_{t\ge 0}$ is a $C_0$-semigroup on $H$, then $(\auxMult_{T(t)})_{t\ge 0}$ is a $C_0$-semigroup on $L^2(\Delta,e_{-\lambda}; H)$, which we denote by $\auxMult_{\mathcal{T}}$.
Since $\auxMult_{\mathcal{T}}$ commutes with the right shift semigroup $\mathcal R_\Delta=(R_\Delta(t))_{t\ge 0}$ on $L^2(\Delta,e_{-\lambda}; H)$, the product $\mathcal R_\Delta\auxMult_{\mathcal{T}}=(R_\Delta(t)\auxMult_{T(t)})_{t\ge 0}$ is also a $C_0$-semigroup on $L^2(\Delta, e_{-\lambda}; H)$, called the \emph{evolution semigroup}. For further discussion of evolution semigroups and related topics we refer to \cite[Section~2]{chicone1999evolution}. Since in both $\mathcal R_\Delta$ and $\auxMult_{\mathcal{T}}$ the parameters $\Delta$ and $\mathcal{T}$ vary frequently, we refrain from introducing a separate notation for the evolution semigroup. This keeps the exposition clear and explicit.

We will also need several auxiliary notions and properties from semigroup theory.
Recall that the \textit{exponential growth bound} $\omega_0(\mathcal T)$ of a $C_0$-semigroup $\mathcal T = (T(t))_{t\geq0}$ is defined by
\begin{equation*}
	\omega_0(\mathcal T) = \lim_{t\to \infty} \frac{\log \|T(t)\|}{t}  \in [-\infty, \infty).
\end{equation*}
Note that $\omega_0(\mathcal T) = \frac{\log r(T(t))}{t}$, $t>0$, where $r(\cdot)$ denotes the spectral radius,
and
$$\omega_0(\mathcal T) = \inf\{ \omega \in \RR \, : \, \exists K\geq 1 \mbox{ such that } \|T(t)\| \leq K e^{\omega t}, \, t \geq 0\}.
$$
Every $C_0$-semigroup $\mathcal T=(T(t)_{t \ge 0})$ on $H$, with generator $\generator,$
can be recovered from $\generator$ via various approximation formulas. Among most well-known ones is
the Post-Widder inversion formula asserting that
	\begin{equation}\label{post}
		T(t)h = \lim_{n\to \infty} 
		 \left(I -\frac{t}{n}\generator\right)^{-n} h, 
			\end{equation}
	where the convergence is uniform in $t$ for compact subsets of $(0,\infty)$; see, for instance, \cite[Corollary III.5.5]{engel2000one}. This formula will allow us to transfer similarity properties from resolvents to semigroups.

As usual in the study of similarity, we will be using a Banach limit on $\ell^\infty(\mathbb N).$
We will not distinguish between different Banach limits on $\ell^\infty(\mathbb N)$,  will denote a fixed
Banach limit by $\rm{LIM},$ and will write ${\rm LIM} [(a_k)]$ for $(a_k)_{k\in \NN} \in \ell^\infty(\mathbb N).$
While Banach limits on $L^\infty(\mathbb R_+)$ could also be used for the similar purposes,
we find their discrete analogues more transparent and are easier to deal with.
	
\medskip

\section{Notation} \,\, We collect  notations which will be used throughout  the paper.
By $\mathcal L (H, K)$ we denote the Banach space of bounded linear operators
between Hilbert spaces $H$ and $K,$ and we will write just
$\mathcal L(H)$
when $H=K.$ The symbol $I_H$ denotes the identity operator on $H$, or simply $I$ when the choice of space is apparent.
Writing $H \oplus K$ we will always  mean the orthogonal direct sum of Hilbert spaces $H$ and $K$,
 and if $K \subseteq H$, we will use $H \ominus K$ to denote the orthogonal difference of $H$ and $K$. The tensor product  of Hilbert spaces $H$ and $K$ will be denoted by $H \otimes K,$
 and $T\otimes S \in \linearOp(H\otimes K)$ will stand 
 for the tensor product of  $T \in \linearOp(H)$ and $S \in \linearOp(K).$ 
  For a closed subspace $H$ of a Hilbert space $\mathcal H$, $\proj_H \in \linearOp(\mathcal H)$ will  denote the orthogonal projection from $\mathcal H$ onto $H$, and given $T \in \linearOp(\mathcal H)$, $T\restriction_H \in \linearOp(H,\mathcal H)$ will denote the restriction of $T$ to $H$.

As is customary in the literature, with a slight abuse of notation,
we write $\|T\|$ for the norm of a bounded operator $T$ on a Hilbert space with fixed norm $\|\cdot\|.$ This helps avoid cumbersome notation when, in particular, the space is considered with multiple norms.

For an arbitrary linear operator $E$ on a Hilbert space $H$, we denote its domain by $\dom(E)$, its range by $\ran(E)$, and its resolvent set by $\rho(E)$. Even if $E$ is not densely defined, $\rho(E)$ will be understood as the set of $\lambda \in \mathbb C$ such that $(\lambda -E)^{-1}$ exists and belongs to $\mathcal L(H)$. Then its spectrum $\sigma(E)$ is defined as $\sigma(E) := \CC \setminus \rho(E)$. 
If $E$ is densely defined, then $E^\ast$ denotes its Hilbert space adjoint. If in addition $E$ is bounded, then we denote by $r(E)$ its spectral radius.

For $\lambda \in \CC$, $e_\lambda$ denotes the exponential function $t \mapsto e^{\lambda t}$, and its domain will typically be clear from context. Consequently, given a $C_0$-semigroup $\mathcal T = (T(t))_{t\geq0}$ and $\lambda \in \CC$, $e_{\lambda} \mathcal T$ denotes the rescaled semigroup $(e^{\lambda t} T(t))_{t\geq0}$. Also, $\omega_0(\mathcal T)$ stands for the exponential growth bound of $\mathcal T$.

The Lebesgue measure will be denoted by ${\rm meas}$ and its domain
will be clear from context. 

For any complex number $z$ its real part will be denoted by $\Real z,$  $\mathbb T$ will stand for the unit circle and $\CC^+$ for the right-half plane. We write $\ZZ_+$ for the set of non-negative integers, and $\RR_+$ for $[0,\infty)$. Finally, given a subset $\Delta \subseteq \RR$, $\chi_\Delta$ denotes the characteristic function of $\Delta$.

\section{Holbrook-type similarity conditions}

\subsection{Holbrook-type condition for semigroups}\label{Holbrook_cond}

We begin with general results on similarity to contractive semigroup representations, which form the core of many arguments in this paper

It is well-known that the equality \eqref{SimEq} can be relaxed by replacing $R$ and $R^{-1}$ by any pair of bounded operators $A$ and $B$ whenever the families of operators $(T_\gamma)_{\gamma \in \Gamma}$ and $(AT_\gamma B)_{\gamma \in \Gamma}$ are semigroup representations; see, for instance, \cite[Theorem 5 \& Corollary 1]{lebow1975spatial}, \cite[Proposition 4.2]{pisier2001similarity} or \cite[Proposition 2.4]{arhancet2017isometric}
(where the latter result was proved in \cite[Proposition 5.5.6]{fackler2015thesis}).

Below, we include  a proof of this fact for a general case of representations of arbitrary (in general, non-abelian) semigroup. It is a variation upon known arguments though looks somewhat simpler. The statement in such a form  will be needed to prove Proposition \ref{HolSemigroup} crucial for the sequel
and covering, in particular,  the case of semigroups $[\tau, \infty)$ for some fixed $\tau\geq0.$ 
As usual, for any representation $(T(g))_{g \in G}$ of a semigroup $G$
with unit $e$ on $\linearOp(H),$ we assume that $T(e)=I.$

\begin{proposition}\label{compressionProp}
	Let $G$ be a unital semigroup, let $H$ and $K$ be two Hilbert spaces, and let $\mathcal T = (T(g))_{g\in G} \subseteq \linearOp(H)$ and $\mathcal S = (S(g))_{g\in G} \subseteq \linearOp(K)$ be
	representations of $G$ in $H$ and $K$, respectively. Assume there exist  $A \in \linearOp(K, H)$ and $B \in \linearOp(H, K)$ such that
	\begin{equation}\label{holbeq}
		T(g) = A S(g) B, \qquad g\in G.
	\end{equation}
		Then there exist $\mathcal S$-invariant subspaces $K_1,K_2$ of $K$ with $K_1\subseteq K_2,$
	and an isomorphism $\mathscr A  \in \linearOp(K_2\ominus K_1, H)$ satisfying $\|\mathscr A\| \|\mathscr A^{-1}\| \leq \|A\| \|B\|$ and
	$$T(g) = \mathscr A \proj_{K_2 \ominus K_1} S(g) \restriction_{K_2 \ominus K_1} \mathscr A^{-1}, \qquad g \in G.
	$$ 
	Thus $\mathcal T$ is similar to the compression of $\mathcal S$ to $K_2\ominus K_1.$
\end{proposition}
\begin{proof}
	Note that
	\begin{equation}\label{spacef}
	K_2: = \overline{\operatorname{span}} \{S(g)Bh \, : \, g \in G, \, h \in H\},
	\end{equation}
	is a closed, $\mathcal S$-invariant subspace of $K$ that contains the range of $B$. Set also $K_1:= K_2 \cap \ker A$.

	Since $AB = T(e) = I$, where $e$ is the unit of $G$, we infer that $A(K_2) = H$. Thus it follows that
	the induced map $\mathscr A \in \linearOp( K_2 \ominus K_1, H),$ given by
	$\mathscr A=A\restriction_{K_2 \ominus K_1}$ is an
	isomorphism. In addition, as $\mathscr A \proj_{K_2\ominus K_1} B = AB = I$, one obtains $\mathscr A^{-1} = \proj_{K_2\ominus K_1} B$.
	
	Furthermore, for all $g\in G$, $\widetilde g \in G$, and $h \in H,$
	$$T(g) A S(\widetilde g)Bh =
	T(g \widetilde g)h = AS(g  \widetilde g) Bh = A S(g) S(\widetilde g) Bh.
	$$
	One deduces from the above equality that $T(g)A f = A S(g)f$ for all $g\in G$ and all $f\in K_2$ and, consequently,  $K_1$ is $\mathcal S$-invariant.
	Since $K_1 \subseteq \ker A$, one has
	$$A S(g) \restriction_{K_1} = 0, \qquad g \in G,
	$$
	and since $\ran(B) \subseteq K_2$, one obtains
	$$
	A S(g) B = A S(g) \restriction_{K_2\ominus K_1} \, \proj_{K_2\ominus K_1}  B + A S(g) \restriction_{K_1} \, \proj_{K_1} B = A S(g) \restriction_{K_2\ominus K_1} \, \proj_{K_2\ominus K_1}  B, \qquad g \in G,
	$$
	implying that
	\begin{align*}
		\mathscr A \proj_{K_2 \ominus K_1} S(g) \restriction_{K_2\ominus K_1} \,  \mathscr A^{-1} &= 
		A S(g) \restriction_{K_2\ominus K_1} \, \proj_{K_2\ominus K_1} B = A S(g) B = T(g),\qquad g \in G.
	\end{align*}
	Thus, we conclude that $\mathcal T$ is similar to the compression
	$\left(\proj_{K_2\ominus K_1} S(g)\restriction_{K_2\ominus K_1} \right)_{g\in G}$
	of $\mathcal S$ to $K_2 \ominus K_1$.
	
	Finally, by construction, we have $\|\mathscr A\| \leq \|A\|$ and $\|\mathscr A^{-1}\| \leq \|B\|$.
\end{proof}
\begin{remark}\label{similar_dil}
Proposition \ref{compressionProp} implies that 
if a representation $\mathcal S$ dilates $\mathcal T$ and $\mathcal S$ is similar to a contractive representation, then so is $\mathcal T.$
This property will be important in what follows.
\end{remark}

If $G=\mathbb Z_+,$ then Proposition \ref{compressionProp} can be improved substantially.
Holbrook
showed in \cite{holbrook1973operators}
that a bounded operator $T$ on a Hilbert space $H$ is similar to a contraction if and only if there exist a Hilbert space $K$, operators $B\in \linearOp(H, K)$ and $A \in \linearOp(K,H)$, and a contraction $S$ on $K$ such that
\begin{align}\label{HolbrooksIdentity}
	\sum_{n=0}^\infty \|T^n - A S^n B\|^2 < \infty.
\end{align}
Holbrook proves this by explicitly constructing an equivalent norm on $H$ under which $T$ becomes a contraction. This norm appears to be induced by an inner product,
since it was shown to satisfy the parallelogram law.
See also the proof of \cite[Theorem 9.1]{paulsen2002completely} for a similar approach.

We will need the next extension of Holbrook's result to the setting of representations of $G$ on $H$. Its proof relies on Proposition \ref{compressionProp} and employs an appropriate  quotient space to produce an equivalent norm making the representation contractive.
Thus the Hilbertian nature of the norm becomes apparent.
This approach simplifies the technical details, makes the argument more transparent, and reveals hidden aspects of Holbrook's original construction\footnote{It was claimed in \cite[p. 164]{holbrook1973operators} (without proof) that the proof of \eqref{HolbrooksIdentity} can be reduced to the case $T^n = A S^n B$ for $n \in \mathbb{Z}_+$.}.
Note that implicitly the quotient space appeared in the arguments leading to \cite[Theorem 2.1.3 \& Theorem 4.1]{badea2003operators} for $G = \ZZ_+$.

{To realize our approach we have to impose additional assumptions  on the semigroup $G$. Recall that 
a semigroup $(G,\cdot)$ is left-cancellative if $f\cdot a=f\cdot b$ implies  $a=b$ for all $f, a, b \in G$. The class of left-cancellative semigroups is quite wide,
and includes, in particular, subsemigroups of groups.}
%

\begin{proposition}\label{HolSemigroup}
Let $G$ be an infinite left-cancellative semigroup and let 
$\mathcal T = (T(g))_{g\in G}$ be a representation of $G$ on a Hilbert space $H$. 
Assume there exist a Hilbert space $K$ and a contractive representation 
$\mathcal S = (S(g))_{g\in G}$ on $K$ such that
\begin{equation}\label{quadnear}
   \sum_{g \in G} \|T(g) - A S(g) B\|^2 < \infty
\end{equation}
for some $A \in \mathcal L(K,H)$ and $B \in \mathcal L(H,K)$. 
Then the operators $(T(g))_{g\in G}$ are jointly similar to contractions.
\end{proposition}
\begin{remark}
		Note that we require neither regularity of $\mathcal T$ nor commutativity of $G$. 	
		The condition \eqref{quadnear}, however, forces 
		$T(g)=AS(g)B$ for all but countably many $g \in G$. 
		Yet it remains of independent value, as Corollary~\ref{HolbrooksCont} shows.
\end{remark}

\begin{proof}
We may assume that $G$ has an identity $e$, by adjoining one if necessary 
and setting $T(e)=S(e)=I$. The summability condition \eqref{quadnear} remains unaffected.
For clarity, we divide the proof in several steps.

Let $\ell^2(G, H)$ be the Hilbert space of functions $h: G \to H$ satisfying $\sum_{u \in G}\|h(u)\|^2<\infty,$
and consider a subset $\widehat K_{00}\subset \ell^2(G, H)$ consisting of  functions
	with finite support.
	For $h^1,h^2\in \widehat K_{00}$ put
\begin{equation}\label{inner}
\langle h^1,h^2\rangle_{\widehat K}
   := \Big\langle \sum_{g\in G} S(g)B\,h^1(g),\; \sum_{g\in G} S(g)B\,h^2(g)\Big\rangle_K
      \;+\; \sum_{g\in G} \langle h^1(g),h^2(g)\rangle_H.
\end{equation}
This defines an inner product, with associated norm denoted $\|\cdot\|_{\widehat K},$
so that for $h \in \widehat K_{00}$ its squared norm $\|h\|_{\widehat K}^2$ is a sum of the ``$K$-part'' 
	and ``$H$-part''.
Write $\widehat K$ for the completion of $(\widehat K_{00},\langle\cdot,\cdot\rangle_{\widehat K})$.

Next, for $f \in G$ define the left shift operator $T_L(f):\widehat K_{00}\to \widehat K_{00}$ 
		by
\begin{equation}\label{shift}
   (T_L(f)h)(u) := \sum_{g:\,fg=u} h(g), \qquad h\in \widehat K_{00}, u\in G.
\end{equation}
Because $G$ is left-cancellative, the selfmap $L_f(g)=fg$  of $G$ is injective for every $f \in G$, hence for each $u$ 
the sum in \eqref{shift} has at most one nonzero term.
Note that equivalently 
\[
   T_L(f)h \;=\; \sum_{g\in G} \delta_{fg}\, h(g), \qquad f\in G, \quad h\in \widehat K_{00},
\]
with $\delta_{u}$ denoting the standard basis vector of
$\ell^2(G)$ supported at $u.$

We claim that $T_L(f)$ is contractive on $(\widehat K_{00},\langle\cdot,\cdot\rangle_{\widehat K})$.
To prove this we consider the $K$-part and the $H$-part of the norm $\|\cdot\|_{\widehat K}$ separately.
To estimate the $K$-part of $\|T_L(f)\|_{\widehat K},$  
we use \eqref{shift} and note that by reindexing
\begin{align*}
   \sum_{u\in G} S(u)B (T_L(f)h)(u)
   &= \sum_{g\in G} S(fg)B\,h(g)
    = S(f)\sum_{g\in G} S(g)B\,h(g).
\end{align*}
Hence, since $S(f)$ is a contraction,
\begin{equation}\label{kpart}
   \Big\|\sum_{u\in G} S(u)B\,(T_L(f)h)(u)\Big\|_K
   \le \Big\|\sum_{g\in G} S(g)B\,h(g)\Big\|_K .
\end{equation}
For the $H$-part of $\|T_L(f)\|_{\widehat K},$ 
 we observe that 
for every $f \in G,$ using the injectivity of $L_f,$
\begin{align}\label{hpnorm}
  \sum_{u\in G}\|(T_L(f)h)(u)\|^2
    = \sum_{u\in \operatorname{ran} L_f}\|h(L_f^{-1}(u))\|^2 
    = \sum_{g\in G}\|h(g)\|^2. 
\end{align}

Combining \eqref{kpart} and \eqref{hpnorm} with the definition of
$\|\cdot\|_{\widehat K}$, we have
\[
   \|T_L(f)h\|_{\widehat K}^2
      \le \Big\|\sum_{g\in G} S(g)B\,h(g)\Big\|_K^2 + \sum_{g\in G}\|h(g)\|^2
   = \|h\|_{\widehat K}^2 .
\]
Thus $T_L(f)$ is contractive on $\widehat K_{00}$.
Since $\widehat K_{00}$ is dense in $\widehat K$, $T_L(f)$ extends uniquely by continuity to a contraction on $\widehat K$, and, denoting the extension by the same symbol, $\mathcal T_L=(T_L(f))_{f\in G}$ is a contractive representation on $\widehat K$.

Next we construct auxiliary intertwining and embedding operators.
First, define the linear operator $Q:\widehat K_{00}\to H$ by
\[
   Qh := \sum_{g\in G} T(g)h(g), \qquad h \in \widehat K_{00}.
\]
For any $h\in\widehat K_{00}$, we estimate
\begin{align*}
   \|Qh\|
   &\le \Big\|\sum_{g\in G}AS(g)B\,h(g)\Big\| + \sum_{g\in G} \|(T(g)-AS(g)B)h(g)\|  \\
   &\le \|A\|\Big\|\sum_{g\in G} S(g)B\,h(g)\Big\|_K
       + \Big(\sum_{g\in G} \|T(g)-AS(g)B\|^2\Big)^{1/2}
         \Big(\sum_{g\in G} \|h(g)\|^2\Big)^{1/2}   \\
   &\le \Big(\|A\|^2 + \sum_{g\in G} \|T(g)-AS(g)B\|^2\Big)^{1/2}\,\|h\|_{\widehat K}.
\end{align*}
Hence, using again the density of  $\widehat K_{00}$ in $\widehat K$, we conclude that $Q$ extends to a bounded operator from $\widehat K$ to $H,$ which we also denote by $Q$.
Moreover, for all $f\in G$ and $h\in \widehat K_{00}$, we have
\begin{align*}
   Q(T_L(f)h)
	   = \sum_{g\in G} T(fg)\,h(g)
          = T(f)(Qh).
\end{align*}
Thus $Q T_L(f) = T(f)Q$ on $\widehat K$, and in particular 
$\ker Q$ is $\mathcal T_L$-invariant.

In addition, define the bounded linear operator $\pi:H\to \widehat K$
 by
\[
   [\pi (h)](g)=
   \begin{cases}
      h, & g=e,\\
      0, & g\neq e,
   \end{cases}
   \qquad h\in H.
\]
(Equivalently, identifying $\ell^2(G,H)$ with $\ell^2(G)\otimes H,$ one can write 
$\pi h=\delta_e\otimes h$ for all $h \in H$.) 
Then $Q \pi$ is the identity map on $H$, so
	\begin{equation*}
		T(f) = T(f) Q \pi = Q T_L(f) \pi , \qquad f \in G.
	\end{equation*}
	
	It then follows from Proposition \ref{compressionProp} that $\mathcal T$ is similar to a compression of $\mathcal T_L$. More precisely,
	if $F\subseteq \widehat K$ is
	given by
	\begin{equation}\label{spaceff}
	F: = \overline{\operatorname{span}} \{T_L(f) \pi  h \, : \, f \in G, \, h \in H\},
	\end{equation}
		(cf. \eqref{spacef} for a similar construction), then $F$ is a closed, $\mathcal T_L$-invariant subspace of $\widehat K$, and  $\mathcal T$ is similar to the compression of $\mathcal T_L$ 
		to $F\ominus (F \cap  \ker Q)$.
		\end{proof}

\begin{remark}\label{BoundSimConstantRemark}
	Under the assumptions of Proposition \ref{HolSemigroup} 
	(in particular, with $G$ left-cancellative),
		 it follows from Proposition \ref{compressionProp} that the similarity constant $\mathcal C(\mathcal T)$ of  $\mathcal T$ is bounded by $\|Q\| \|\pi\|$, where the operators $Q$ and $\pi$ are defined in the proof of Proposition \ref{HolSemigroup}. Since
	for any $\gamma>0$ one may replace $A$ and $B$ in \eqref{quadnear} by $\gamma A$ and $\frac{1}{\gamma} B$, we infer that
		\begin{equation}\label{constbounds}
	\mathcal C(\mathcal T) \leq \|Q\| \| \pi \| \leq \left( \gamma^2 \|A\|^2 + \sum_{g\in G \cup \{e\}} \left\| T(g) - A S(g) B \right\|^2 \right)^{1/2} \left(1 + \gamma^{-2} \|B\|^2\right)^{1/2},
	\end{equation}
	where $e$ denotes the identity element of $G$ if $G$ has identity, or the adjoined identity otherwise.
	If $A\neq 0$, $B\neq0$ and $T(g) \neq A S(g) B$ for some $g \in G\cup\{e\},$ then setting 
	$$\gamma^2 = \frac{\|B\|}{\|A\|} \left(\sum_{g\in G \cup \{e\}} \left\| T(g) - A S(g) B \right\|^2\right)^{1/2},
	$$
	in the right hand-side of \eqref{constbounds},
	we obtain
	\begin{equation}\label{Cbound}
		\mathcal C(\mathcal T) \leq \|A\| \| B \| + \left( \sum_{g\in G \cup \{e\}} \left\| T(g) - A S(g) B \right\|^2 \right)^{1/2}.
	\end{equation}
	Moreover, if either $A = 0$ or $B=0$,  then  passing to the limit in \eqref{constbounds}  as $\gamma \to \infty$ or $\gamma\to 0$ respectively, we show \eqref{Cbound}. If $T(g) = A S(g) B$ for all $g \in G\cup \{e\}$, then \eqref{Cbound} follows directly from Proposition~\ref{compressionProp}.
\end{remark}
Returning to the proof of Proposition~\ref{HolSemigroup}, note that the left–cancellativity of $G$ is exactly what makes the left shift $T_L(f)$
contractive on the $\ell^2$–part of the inner product on $\widehat K$.
Indeed, for $h \in \widehat K_{00}$ one has
\[
(T_L(f)h)(u)=\sum_{g:\,fg=u} h(g), \qquad u \in G,
\]
so if $G$ is not left–cancellative, different ``indices'' $g$ may merge under
$g\mapsto fg$, and the $\ell^2$–sum of squares can increase.
The next proposition illusrates and clarifies this point.
\begin{proposition}\label{lef_can}
For every Hilbert space $H$ with $\dim H \geq 2$, there exist an infinite, non-left-cancellative semigroup $G$ with unit,  a contractive representation
$\mathcal S = (S(g))_{g\in G}$ on $H$, and  a representation
$\mathcal T = (T(g))_{g\in G}$ such that
\[
\sum_{g\in G} \|T(g)-AS(g)B\|^2 < \infty,
\]
for some $A, B \in \linearOp(H),$ but $\mathcal T$ is not similar to a contractive representation on $G$.
\end{proposition}
\begin{proof}
Let $G$ be the left-zero semigroup on $\mathbb N$ with the adjoined unit $e$,
i.e. we define 
$
m \cdot n := m, $ for all $m,n \in \mathbb N,$
and adjoin $e$ setting $e\cdot n=n\cdot e$ for all $n \in \mathbb N.$ 
Let $M$ be a non-zero subspace of $H$  such that
$M^\perp$ is also non-zero.
For each $n \in \mathbb N$, pick a bounded operator
$
D_n : M \to M^\perp
$
such that $D_n, n \in \mathbb N, $ are pairwise different and
\[
\sum_{n=1}^\infty \|D_n\|^2 < \infty.
\]
(For instance, one may choose rank-one operators $D_n$ with $\|D_n\| = 2^{-n}$.)

Define idempotents $E_n \in {\mathcal L}(H), n \in \mathbb N,$ by the $2 \times 2$ block matrices 
\[
E_n =\begin{pmatrix}
I_M & 0 \\
D_n & 0
\end{pmatrix}
\]
relative to the decomposition $H=M \oplus M^\perp.$
Then it is easy to see that
\begin{equation}\label{e}
E_m E_n = E_m, \qquad m,n \in \mathbb N.
\end{equation}
Hence $\mathcal T = (T(n))_{n\in \NN}$, given by $T(e)=I$ and $T(n) := E_n$ for each $n \in \mathbb N,$
is a representation of $(\mathbb N, \cdot)$.

Now define the contractive representation $\mathcal S  =(S(n))_{n\in \NN}$ of $G$ on $H$ by $S(e):=I$ and
$S(n) := P_M$
for each $n \in \mathbb N,$
where $P_M$ is the orthogonal projection onto $M$. 
If $A=B=I,$ 
then
\[
\sum_{n=1}^\infty \|T(n)-AS(n)B\|^2 = \sum_{n=1}\|E_n-P\|^2=\sum_{n=1}^\infty \|D_n\|^2 < \infty.
\]
Assume for contradiction that there exists an invertible $R$ such that $\mathcal T_0 = (T_0(n))_{n\in \NN}$,
defined by $
T_0(n):= R^{-1}T(n)R$ for every $n \in \mathbb N,$
is
a contractive representation.
For each $n \in \mathbb N,$ since $T_0(n)$ is an idempotent, and $\|T_0(n)\|\le 1,$ we conclude that $T_0(n)$ is an orthogonal projection. In view of \eqref{e}, we have
\[
T_0(m)T_0(n) = T_0(m),  \qquad m,n \in \mathbb N.
\]
Hence
$\ran T_0(m) \subset \ran T_0(n)$ for all $m,n \in \mathbb N$, and then all of $T_0(n)$ coincide.
Letting
$P:= T_0(m)$
for any $m \in \NN$, it follows that $T(n) = RPR^{-1}$ is independent of $n$, which contradicts the choice of pairwise different $E_n.$ This contradiction shows $\mathcal T$ is not similar to a contractive representation, and proves the statement.
 \end{proof}
\begin{remark}
If $G$ is finite, then any sum of the form $\sum_{g\in G}\|T(g)-AS(g)B\|^2$ is finite a priori. Therefore the ``quadratic'' closeness hypothesis \eqref{quadnear} is substantive only for infinite semigroups/groups.
On the other hand, it is easy to show that if $G$ is a finite, unital, and left-cancellative semigroup, then $G$ is a group, and any representation of $G$ is similar to a unitary one. 
The left-cancellation property cannot be dropped even in the finite case, since the construction of representation  in Proposition \ref{lef_can} can be easily adapted to cover the case of left-zero semigroup defined on a finite set $\{1, \dots N\}$ for any $N \in \mathbb N, N \ge 2,$ 
with adjoined unit. 
\end{remark}

The assumptions of Proposition \ref{HolSemigroup} are natural in the context of a single bounded operator \cite{holbrook1973operators} (i.e., when the semigroup $G$ is $\mathbb{Z}_+$), but they are less suitable for $C_0$-semigroups. More precisely, let $G=\mathbb R_+$ so that $\mathcal{T} = (T(t))_{t \ge 0}$ is a $C_0$-semigroup on a Hilbert space. Then the hypotheses of Proposition \ref{HolSemigroup} hold if and only if $T(t) = AS(t)B$ for all $t \ge 0$. Nevertheless, Proposition \ref{HolSemigroup} remains useful for semigroups such as $(\mathbb{Z}_+^n,+)$ or $(\prod_{i=1}^n[\tau_i,\infty),+)$ with $\tau_i>0$ for all $i=1,\dots, n,$
and adjoined $0$ as identity. 
 This observation is crucial for the remainder of the paper,
 and the next statement makes it precise.

\begin{corollary}\label{HolbrooksCont}
	Let $\mathcal T = (T(t))_{t\geq0}$ be a semigroup of bounded operators (not necessarily strongly continuous) on a Hilbert space $H$. Assume there exist $\tau \geq 0$, a Hilbert space $K$,  and a contraction semigroup $(S(t))_{t\geq0}$ on $K$ satisfying
	$$T(t) = A S(t)B, \qquad t \geq \tau,
	$$
	for some $A \in \linearOp(K,H)$,  $B \in \linearOp(H,K).$
	Then the operators $\mathcal T_\tau = (T(t))_{t\geq \tau}$ are jointly similar to contractions.
		Moreover,
		$$
	\mathcal C (\mathcal T_\tau)\leq \|A\| \|B\| + \|I - AB\|.
	$$
\end{corollary}
Since the semigroup  $([0,\infty),+)$ is left-cancellative, the statement is an immediate corollary of Proposition \ref{HolSemigroup}.

\subsection{Holbrook-type condition away from the origin}\label{AwayzeroSect}

In this section we prove Theorem~\ref{awayzero}, a continuous-parameter analogue of Holbrook-type similarity criteria. The theorem concerns the behavior of 
$\mathcal{T}=(T(t))_{t\ge0}$ for $t$ bounded away from zero. 
Unlike in the discrete case, however, the continuous setting introduces new phenomena: in general it is not possible to extend joint similarity of 
$(T(t))_{t\ge\tau}$ to contractions to the entire family $(T(t))_{t\ge0}$. 
This obstruction will be explored in detail in Section~\ref{ExamplesSect}.

In fact, we prove a stronger quantitative version of the implication (i) $\implies$ (iii) in Theorem~\ref{awayzero}, formulated as follows.
\begin{theorem}\label{awayzeroSCboundTh}
	Let $\mathcal T = (T(t))_{t\geq0}$ be a $C_0$-semigroup on a Hilbert space $H$,
	and let $\tau_1 > 0$ be such that $T(\tau_1)$ is similar to a contraction on $H.$
	Then,  for every $\tau_2 \in (0,\tau_1]$, there exist a Hilbert space $K$, $A \in \linearOp(K,H)$, $B \in \linearOp(H,K)$ and a unitary $C_0$-group $\mathcal U = (U(t))_{t\in \RR}$ on $K$ such that
	\begin{equation}\label{InterpLarget}
		T(t) = A U(t)B, \qquad t \geq \tau_2,
	\end{equation}
	and
		\begin{equation}\label {boundab}
	\|A\| \|B\| \leq C(T(\tau_1)) M^2 \sqrt{\frac{\tau_1}{\tau_2}},
	\end{equation}
	where $M = \sup_{t \in [0,\tau_1]} \|T(t)\|$.
\end{theorem}
\begin{proof}[Proof of Theorems \ref{awayzero} and \ref{awayzeroSCboundTh}]
	
	Let us first prove Theorem \ref{awayzeroSCboundTh}, and thus
	(i) $\implies$ (iii).
	Since $T(\tau_1)$ is similar to a contraction, there exists
			a Hilbertian norm  $\|\cdot\|_{{\rm eq}}$
			on $H$ such that $\|T(\tau_1)\|_{\rm {eq}}\leq 1$ and $\|h\| \leq \|h\|_{{\rm eq}} \leq C(T(\tau_1)) \|h\|$ for all $h \in H.$
	Define then the norm $\|\cdot\|_\mathscr H$ on $H$ by
	\begin{align*}
		\|h\|_\mathscr H^2 := \frac{1}{M^2  \tau_1} \int_0^{ \tau_1} \|T(s)h\|^2_{{\rm eq}} \,ds, \qquad h \in H,
	\end{align*}
	where $M =\sup_{t\in [0,\tau_1]} \|T(s)\| < \infty$, and note that
	$\|h\|_\mathscr H \leq C(T(\tau_1)) \|h\| $ for all $h\in H$.
	
	Let $\mathscr H$ be the completion of $H$ equipped with $\|\cdot\|_\mathscr H$. Then $\mathscr H$ is a Hilbert space, the embedding  $\embedding : H \to  \mathscr H$ is bounded, with $\|\embedding\|\leq C(T(\tau_1))$, and has dense range. Moreover, for every $t \in (0, \tau_1)$,
	\begin{equation}\label{boundst}
		\begin{aligned}
		\|\embedding T(t)h\|_\mathscr H^2 & =
		\frac{1}{M^2  \tau_1} \left( \int_{t}^{ \tau_1} \|T(s)h\|^2_{\rm eq}\,ds +
		\int_{0}^{t} \|T(s+\tau_1)h\|^2_{{\rm eq}}\,ds \right)
		\\ & \leq \frac{1}{M^2  \tau_1} \left( \int_{t}^{ \tau_1} \|T(s)h\|^2_{{\rm eq}}\,ds
		+ \int_{0}^{t} \|T(s)h\|^2_{{\rm eq}}\,ds \right) = \|\embedding h\|_\mathscr H^2, \qquad h \in H.
		\end{aligned}
	\end{equation}
	For fixed $t \ge 0$ let $S(t)\embedding h: = \embedding T(t)h$, $h\in H$.
	As $\embedding$ has dense range, \eqref{boundst} implies that $S(t)$ extends to a contraction on $\mathscr H,$
	and thus $\mathcal S=(S(t))_{t \ge 0}$ is a semigroup of contractions on $\mathscr H.$
	Since by the dominated convergence theorem,
		\begin{align*}
		\lim_{t \to 0} \|S(t)\embedding h - \embedding h\|_\mathscr H^2 = \lim_{t \to 0} \frac{1}{M^2  \tau_1} \int_0^{ \tau_1} \|(T(t)-I)T(s) h\|_{{\rm eq}}^2 \,ds=0
	\end{align*}
for all $h \in H,$	we conclude that  $\mathcal S$ is a $C_0$-semigroup on $\mathscr H$.
		
Furthermore, fix  $\tau_2 \in (0,\tau_1]$	
and	define the operator $\mathscr A$ from $\embedding(H)$ to $H$ by $\mathscr A (\embedding h): = T( \tau_2)h$ for all $h \in H$. In view of
	\begin{align*}
		\|\mathscr A (\embedding h) \|^2 &= \|T( \tau_2)h\|^2  = \frac{1}{ \tau_2} \int_0^{ \tau_2} \|T( \tau_2-s) T(s) h\|^2 \,ds
		 \\ &\leq  \frac{M^2}{ \tau_2} \int_0^{ \tau_2} \|T(s) h\|^2\,ds
		\leq  \frac{M^2}{ \tau_2} \int_0^{ \tau_1} \|T(s) h\|^2_{{\rm eq}}\,ds
		\notag \\ & = M^4 \frac{ \tau_1}{ \tau_2} \|\embedding h\|_\mathscr H^2, \qquad h \in H,
			\notag
	\end{align*}
we infer that  $\mathscr A$ extends to a bounded operator from $\mathscr H$ to $H$ and $\|\mathscr A\| \leq  M^2 \sqrt{\frac{ \tau_1}{\tau_2}}$.

	Now Nagy's dilation theorem yields a unitary dilation $\mathcal U = (U(t))_{t\in \RR}$ of $\mathcal S$ to a Hilbert space $K$ containing $\mathscr H$ as a subspace, so that $S(t) = \proj_\mathscr H U(t)\restriction_\mathscr H$ for $t\geq0$. 
	Thus, letting $A:= \mathscr A \proj_\mathscr H U(- \tau_2)$ and $B:= \embedding$, we find $A \in \linearOp(K, H)$ and $B \in \linearOp(H, K)$ such that $\ran(B) \subseteq \mathscr H$ and
	\begin{equation}\label{unitEq}
		\begin{aligned}
		AU(t)B &= \mathscr A \proj_\mathscr H U(- \tau_2) U(t)  \embedding
		= \mathscr A \proj_\mathscr H U(t- \tau_2)\restriction_\mathscr H \,  \embedding
		\\ &= \mathscr A S(t- \tau_2)  \embedding
		=\mathscr A \embedding T(t- \tau_2) = T(t), \qquad t \geq  \tau_2.
		\end{aligned}
	\end{equation}
	In addition, we have $\|A\| \leq \|\mathscr A\| \leq M^2 \sqrt{\frac{\tau_1}{\tau_2}}$ and $\|B\| =  \|\embedding\| \leq C(T(\tau_1))$, hence \eqref{boundab} holds, and Theorem \ref{awayzeroSCboundTh}
	follows.
		\medskip
	
The rest of implications are easy consequences of (i) $\implies$ (iii) and
Corollary \ref{HolbrooksCont}:
\medskip
	
	(iii) $\implies$ (ii): If \eqref{InterpLarget} holds, then the statement is a direct consequence
	of Corollary \ref{HolbrooksCont}.
	\medskip
	
	(ii) $\implies$ (iv): If for each $\tau>0$,
 the operators $(T(t))_{t\geq \tau}$ are jointly similar to contractions, then the statement is immediate corollary of the  implication (i) $\implies$ (iii).
	\medskip

	(iv) $\implies$ (i): If there exist $\tau>0$, a Hilbert space $K$, $B \in \linearOp(H, K)$, $A \in \linearOp(K, H)$, and a $C_0$-semigroup $\mathcal S = (S(t))_{t\geq0} \in \mathcal{SC}(K)$
	satisfying $T(t) = A S(t)B$ for all $t \geq \tau,$
	then employing Corollary \ref{HolbrooksCont} again we infer that $T(\tau)$ is similar to a contraction.
		\end{proof}

Recall that if $A$ is a bounded operator on a Hilbert space such that $A^N$ is similar to a contraction for some $N \in \NN$, then $A$ is similar to a contraction; see, for instance, \cite[p. 912]{halmos1970ten}. One may wonder whether a similar property holds for $C_0$-semigroups. More generally, given a semigroup $\mathcal T = (T(t))_{t\geq0}$ one may ask whether the subset of $[0,\infty)$ consisting of those $t\in [0,\infty)$ such that $T(t)$ is similar to a contraction can be different from $[0,\infty)$ and $\emptyset$. Somewhat surprisingly,  Theorem \ref{awayzero} directly implies that the answer is negative, and we formalize it in the next corollary.
\begin{corollary}\label{SimSetCor}
	Let $\mathcal T = (T(t))_{t\geq0}$ be a $C_0$-semigroup on a Hilbert space $H.$
	Then only one of the following holds.
	\begin{enumerate}
		\item [(i)] For each $\tau > 0$, the operators $(T(t))_{t\geq \tau}$ are jointly similar to contractions.
		\item [(ii)] There is no $\tau > 0$ such that $T(\tau)$ is similar to a contraction.
	\end{enumerate}
\end{corollary}
\begin{proof}
	This is an immediate consequence of (i) $\implies$ (ii) of Theorem \ref{awayzero}.
\end{proof}

Many results concerning similarity to contractions for single operators can be naturally extended to strongly continuous semigroups, provided that the semigroup consists of operators that are bounded from below. This key observation, which we establish in the proposition below, appears to have been overlooked in the literature. In particular, our result generalizes and strengthens Liapunov’s theorem for $C_0$-groups \cite[Subsect. 7.2]{haase2006functional}
and extends the main result of \cite{vu1998similarity}, where the authors proved that a uniformly continuous, quasi-compact, and bounded semigroup belongs to $\mathcal{SC}(H)$.

Note that boundedness from below has a number of special features.
It is well-known and easy to prove that the boundedness from below of $T(t)$ for a single $t >0$ implies the same property for all $t >0$ and in fact implies the existence of $\alpha>0$ and $\beta \in \mathbb R$ such that
\begin{equation}\label{leftInvEq}
	\|T(t)h\|\ge \alpha e^{\beta t}\|h\|, \qquad h \in H, \, t \geq 0,
\end{equation}
see e.g. \cite[Lemma 1.2.2]{nikolski1998controllabilite} or \cite[p. 487]{batty2017lower}
 (and cf. \cite[Section 3]{haak2012exact}). Moreover, \eqref{leftInvEq} is equivalent to the left invertibility of $T(t)$ for all $t \geq 0$, meaning that there exists $(S(t))_{t \geq 0}\subset \mathcal L(H)$ satisfying $S(t)T(t) = I$ for all $t \geq 0$. Moreover,  $(S(t))_{t \geq 0}$ can be chosen to be a $C_0$-semigroup. See e.g. \cite[Theorem 7.3]{batty2017lower} or \cite[Theorem 1]{zwart2013left}
for a relevant discussion.

\begin{proposition}\label{boundedBelowContractProp}
	Let $\mathcal T = (T(t))_{t\geq0}$ be a $C_0$-semigroup on a Hilbert space $H$, and assume that there exists $t'>0$ such that $T(t')$ is bounded from below (or, equivalently, $\mathcal T$ is left invertible). Then
	\begin{itemize}
		\item [(i)] $\mathcal T \in \mathcal{SC}(H)$ if and only if there exists $\tau>0$ such that $T(\tau)$ is similar to a contraction.
		\item [(ii)] $\mathcal T \in \mathcal{SQC}(H)$.
	\end{itemize}
\end{proposition}
\begin{proof}
	(i) The ``only if'' part of the claim is trivial. So assume there exists $\tau>0$ such that $T(\tau)$ is similar to a contraction. Since our statement is invariant with respect to equivalent renormings,  without loss of generality, we assume that $T(\tau)$ is a contraction.
	Let $\|\cdot\|_{\mathscr H}$ be the norm on $H$ given by
	$$\|h\|_{\mathscr H} := \left(\frac{1}{\tau} \int_0^{\tau} \|T(t)h\|^2 \, dt\right)^{1/2}, \qquad h \in H.
	$$
	Clearly, $\|h\|_{\mathscr H} \leq M \|h\|$ for all $h \in H$, where $M:= \sup_{t\in[0,\tau]} \|T(t)\| < \infty$.
	Also, it follows from \eqref{leftInvEq} that there exists $c>0$ such that $\|T(t)h\| \geq c \|h\|$ for all $h \in H$ and all $t \in [0,\tau]$.
	Thus,  $\|h\|_{\mathscr H} \geq c \|h\|$ for $h \in H$, implying that $\|\cdot\|_{\mathscr H}$ is an equivalent Hilbertian norm on $H$.
	
	Moreover, since
	\begin{align*}
		\|T(s)h\|_{\mathscr H}^2 &
				= \frac{1}{\tau} \int_s^{\tau} \|T(t)h\|^2 \, dt + \frac{1}{\tau} \int_0^s \|T(\tau) T(t)h\|^2\, dt
		\\ &\leq \frac{1}{\tau} \int_0^{\tau} \|T(t)h\|^2 \, dt = \|h\|_{\mathscr H}^2, \qquad h \in H, \, s \in [0,\tau],
	\end{align*}
	we infer that $\mathcal T$ is a contraction $C_0$-semigroup on
	$H$ endowed with the norm $\|\cdot\|_{\mathscr H},$ i.e. $\mathcal T$ is similar to a semigroup of contractions on $H$. \medskip
	
	(ii) Choose $\lambda\in \RR$ with $e^{\lambda}\ge \|T(1)\|$. Then $e_{-\lambda} \mathcal T = (e^{-\lambda t}T(t))_{t\geq0} \in \mathcal{SC}(H)$ by (i), and our claim follows.
\end{proof}

Proposition \ref{boundedBelowContractProp} implies a version of the Liapunov's theorem
for bounded from below (i.e. left-invertible) $C_0$-semigroups. It generalizes the result in \cite[Corollary 7.2.5]{haase2006functional},
where only exponentially stable $C_0$-groups were addressed.
\begin{corollary}
\label{LiapunovThCor}
	Let $\mathcal T = (T(t))_{t\geq 0}$ be a $C_0$-semigroup such that $T(t')$ is bounded from below for some $t'>0.$ 
		For every $a > \omega_0(\mathcal T)$, there exists an equivalent Hilbertian norm $\|\cdot \|_\mathscr H$ on $H$ such that
	$$\|T(t)\|_\mathscr H \leq e^{at}, \qquad t \geq 0.
	$$
\end{corollary}
\begin{proof}
	Fix $a > \omega_0(\mathcal T)$ and note that $r(e^{-a\tau}T(\tau)) = e^{(-a + \omega_0(\mathcal T))\tau} < 1$ for all $\tau>0$. Then $e^{-a\tau}T(\tau)$ is similar to a contraction for every $\tau>0$ by Rota's theorem, and our claim follows from Proposition \ref{boundedBelowContractProp}(i) applied
	to the	semigroup $e_{-a} \mathcal T$.
\end{proof}

If $\mathcal T = (T(t))_{t\geq0}$ is a $C_0$-semigroup with generator $\generator$, with $T(t')$ bounded from below for some $t'>0,$ then there exists $B \in \linearOp(H)$ such that the semigroup generated by $\generator + B$ is similar to an isometric semigroup \cite[Theorem 3]{zwart2013left} (see also \cite[p. 497-498]{batty2017lower}, so that $(T(t))_{t \ge 0}$ is similar to quasi-contraction semigroup. Thus, Proposition \ref{boundedBelowContractProp}(i) and Corollary \ref{LiapunovThCor} are, in fact, corollaries of our main Theorem \ref{quasiContrTh_int} below, whose proof is independent of both.
However, it is instructive to formulate them here and provide  with direct proofs
as an illustration of condition (i) in Theorem \ref{awayzero}.

\subsection{Holbrook-type condition near the origin}\label{NearzeroSect}

In contrast to the discrete setting, the similarity properties of  semigroups $\mathcal T=(T(t))_{t \ge 0}$
depend crucially on the behavior of $\mathcal T$  near zero. This fact invalidates several continuous counterparts of well-known similarity criteria for single operators and motivates
our studies of ``local'' similarity properties of $\mathcal T$.
 Building on this motivation, we now turn to Theorem \ref{nearzero}, complementing Theorem \ref{awayzero}. Specifically, we focus on the property $T(t) = AS(t)B$
for $t$ in a neighborhood of $0$ and a contraction semigroup $\mathcal S=(S(t))_{t \ge 0}.$
Next, we consider a family of spaces and corresponding operators that will play a crucial role in the sequel. 
Given a Hilbert space $\auxHilbert$, $\nu > 0$, and $\lambda \in \mathbb R$, let
$L^2_{\lambda, \nu}(\auxHilbert)$ denote $L^2([0,\nu],\auxHilbert)$ endowed with the equivalent Hilbertian norm
\begin{equation}\label{QuasiAuxRemarkEq1}
	\|f\|_{L^2_{\lambda, \nu}(\auxHilbert)}^2 := \int_0^{\nu} \|f(x)\|^2 e^{-2\lambda x} \, dx,
	\qquad f \in L^2([0,\nu], \auxHilbert).
\end{equation}
Clearly, $L^2_{\lambda, \nu}(\auxHilbert)$ coincides with $L^2([0,\nu], e_{-\lambda}; H)$ defined in 
Section \ref{PreliminariesSect}.
Let the operator $Q_{\lambda, \nu} \in \linearOp(L^2_{\lambda, \nu}(\auxHilbert), \auxHilbert)$ be given by
\begin{equation}\label{QuasiAuxRemarkEq1bis}
	Q_{\lambda,\nu}f = \int_0^{\nu} f(x) \, dx, \qquad f \in L^2_{\lambda, \nu}(\auxHilbert).
\end{equation}
Recall from Section \ref{PreliminariesSect} that, for $T \in \linearOp(\auxHilbert)$, 
the multiplication operator $\auxMult_T$ acts on any Hilbert space of $\auxHilbert$-valued functions, 
and for a $C_0$-semigroup $\mathcal T = (T(t))_{t\ge 0} \subset \linearOp(\auxHilbert)$, the induced multiplication $C_0$-semigroup is denoted by $\auxMult_{\mathcal T} = (\auxMult_{T(t)})_{t\ge0}$.
If $\mathcal R_{\lambda,\nu} = (R_{\lambda,\nu}(t))_{t\geq0}$ stands for the right shift semigroup on $L^2_{\lambda,\nu}(\auxHilbert)$, then the evolution semigroup $\mathcal R_{\lambda,\nu} \auxMult_\mathcal T$, acting on $L^2_{\lambda,\nu}(\auxHilbert),$ can be written as
\begin{equation}\label{Evsemeq}
	(R_{\lambda,\nu}(t) \auxMult_{T(t)} f)(x) = \begin{cases} T(t) f(x-t), & \text{if } x \geq t,\\[2mm]
		0, & \text{if } x < t, 
	\end{cases}
	\qquad x \in [0,\nu], \, f \in L^2_{\lambda,\nu}(\auxHilbert).
\end{equation}
To simplify notation, in the case $\lambda = 0$, we write $L_\nu^2(\auxHilbert)$, $Q_\nu$, $\mathcal R_\nu$, and $R_\nu(t)$ instead of $L^2_{0,\nu}(\auxHilbert)$, $Q_{0,\nu}$, $\mathcal R_{0,\nu}$, and $R_{0,\nu}(t)$, respectively.

Our arguments for the implication (iv) $\implies$ (i) in Theorem \ref{nearzero} will be based
on the next result describing a splitting phenomenon for an evolution semigroup in $\mathcal{SQC}(L^2_\nu(\auxHilbert))$, see Remark \ref{OTremark}. Note \cite[Proposition 8.9]{kubrusly1997introduction},
where a similar statement was proved for a less demanding
situation of the left shift on $\ell^2(\mathbb N)$ with
operator weight.

\begin{proposition}\label{OT24quasi}
	Let $\nu>0$, let $\auxHilbert$ be a Hilbert space and let $\mathcal T = (T(t))_{t\geq0}$ be $C_0$-semigroup on $\auxHilbert$. Then $\mathcal R_\nu \auxMult_{\mathcal T} \in \mathcal{SQC}(L^2_\nu(\auxHilbert))$ if and only if $\mathcal T \in \mathcal{SQC}(\auxHilbert)$.
\end{proposition}
\begin{proof}
	Assume first that $\mathcal T \in \mathcal{SQC}(\auxHilbert)$, and let $\lambda \geq 0$ and $\|\cdot\|_{\mathscr H}$ 
	be a Hilbertian norm on $\mathscr H$ equivalent to the original norm $\|\cdot\|$ and such that
	$\| T(t)\|_{\mathscr H} \leq e^{\lambda t}$ for all $t\geq0$. Setting $\mathscr H := (\auxHilbert, \|\cdot\|_{\mathscr H})$, observe that the induced norm on $L^2_\nu(\mathscr H)$ is equivalent to the norm on $L_\nu^2(\auxHilbert)$. Thus, in view of
	$$\|R_\nu(t) \auxMult_{T(t)}\|_{L^2_\nu(\mathscr H)} \leq \|T(t)\|_{\mathscr H} \leq e^{\lambda t}, \qquad t \geq0,
	$$
	we conclude that $\mathcal R_\nu \auxMult_\mathcal T \in \mathcal{SQC}(L^2_\nu(\auxHilbert))$. \medskip
	
	Let conversely $\mathcal R_\nu \auxMult_{\mathcal T} \in \mathcal{SQC}(\auxHilbert)$, and let $\|\cdot\|_{\rm{eq}}$ be a Hilbertian norm on $L^2_\nu(\auxHilbert)$ such that 
	\begin{equation}\label{evolSemEq}
		\|\cdot\|_{L_\nu^2(\auxHilbert)} \leq \|\cdot\|_{\rm eq} \leq \mathcal C (\mathcal R_\nu \auxMult_{\mathcal T}) \|\cdot\|_{L_\nu^2(\auxHilbert)}
	\end{equation}
	and for some $\lambda \in \RR$ the semigroup $e_{-\lambda} \mathcal R_\nu \mathcal M_{\mathcal T}$ is contractive on $L^2_\nu(\auxHilbert)$ with respect to $\|\cdot\|_{\rm{eq}}$. 
	Observing that, given $h \in \auxHilbert$,  $t \to h\chi_{[t,\nu]}$ is a continuous $(L_\nu^2(\auxHilbert), \|\cdot\|_{\rm{eq}})$-valued function, define
	$$\|h\|_{\mathscr H}^2 := 
	\int_0^\nu \|h\chi_{[t,\nu]}\|_{\rm{eq}}^2 \, dt, \qquad h \in \auxHilbert. 
	$$
	Then
	\begin{align*}
		\|h\|^2_{\mathscr H} &\leq (\mathcal C(\mathcal R_\nu \auxMult_{\mathcal T}))^2 \int_0^\nu \|h\chi_{[t,\nu]}\|^2_{L_\nu^2(\auxHilbert)} \, dt
		= \frac{(\mathcal C(\mathcal R_\nu \auxMult_{\mathcal T}))^2 \nu^2}{2} \|h\|^2, \qquad h \in \auxHilbert,
	\end{align*}
	where we used that $\|h\chi_{[t,\nu]}\|_{L_\nu^2(\auxHilbert)} = \|h\| (\nu-t)^{1/2}$ for $t\in [0,\nu]$. 
	
	Similarly, from the first inequality in \eqref{evolSemEq} one infers that
	\begin{align*}
		\|h\|^2_{\mathscr H} &\geq 
		\frac{\nu^2}{2} \|h\|^2, \qquad h \in \auxHilbert.
	\end{align*}
	Hence, $\|\cdot\|_{\mathscr H}$ is an equivalent Hilbertian norm on $\auxHilbert$ and we let $\mathscr H := (\mathcal H, \|\cdot\|_\mathscr H)$. Next, for $s\in (0,\nu)$ and $h\in \auxHilbert$, write 
	\begin{align}\label{intSum}
		\|T(s)h\|_{\mathscr H}^2 &= \int_0^s \|T(s)h \chi_{[t,\nu]}\|_{\rm{eq}}^2 \, dt
		+  \int_s^\nu \|T(s)h \chi_{[t,\nu]}\|_{\rm{eq}}^2 \, dt.
	\end{align}
	and estimate each of the two summands on the right hand side of \eqref{intSum} separately.
	For the first summand, one has
	\begin{align*}
		\int_0^s \|T(s)h \chi_{[t,\nu]}\|_{\rm{eq}}^2 \, dt &\leq (\mathcal C(\mathcal R_\nu \auxMult_{\mathcal T}))^2 \|T(s) h\|^2 \int_0^s (\nu-t)\, dt
		\\ &\leq (\mathcal C(\mathcal R_\nu \auxMult_{\mathcal T}))^2 \left(\sup_{t\in [0,\nu]} \|T(t)\| \right)^2  \|h\|^2 \nu s
		\leq D s \|h\|_{\mathscr H}^2,
	\end{align*}
	where we set 
	$$D:= \frac{\nu^3}{2} (\mathcal C(\mathcal R_\nu \auxMult_{\mathcal T}))^2 \left(\sup_{t\in [0,\nu]} \|T(t)\| \right)^2.
	$$
	For the second summand, using the contractivity of $e_{-\lambda} \mathcal R_\nu \mathcal M_{\mathcal T}$  on $(L^2_\nu(\auxHilbert), \|\cdot\|_{\rm{eq}})$, one obtains
	\begin{equation*}
		\begin{aligned}
		\\ &\int_s^\nu \|T(s)h \chi_{[t,\nu]}\|_{\rm{eq}}^2 \, dt = \int_s^\nu \|R_\nu(s) \auxMult_{T(s)}  h \chi_{[t-s,\nu]}\|_{\rm{eq}}^2 \, dt
		\\ \leq& \int_s^\nu e^{2\lambda s} \|h \chi_{[t-s,\nu]}\|_{\rm{eq}}^2 \,dt
		= e^{2\lambda s}\int_0^{\nu-s} \|h \chi_{[t,\nu]}\|_{\rm{eq}}^2 \,dt \leq e^{2\lambda s}\|h\|_{\mathscr H}^2.
		\end{aligned}
	\end{equation*}
	From these inequalities, we infer that $\|T(s)\|_{\mathscr H} \leq \left(e^{2\lambda s}+Ds\right)^{1/2}$ for all $s\in (0,1)$. Since the mapping $s \mapsto (e^{2\lambda s} +Ds)^{1/2}$ is differentiable at $0$, the characterization of quasi-contractivity given by \eqref{quasi_loc} and \eqref{quasi_loc_2} implies that $\mathcal T \in \mathcal{SQC}(\auxHilbert)$ and that $\|T(t)\|_{\mathscr H} \leq e^{(\lambda + D/2)t}$ for $t\geq0$, completing the proof.
	
\end{proof}

\begin{remark}\label{OTremark}
	Fix $\nu>0$ and let $\auxHilbert$ be a Hilbert space and $\mathcal T=(T(t))_{t\geq0}$ be a $C_0$-semigroup on $\auxHilbert$. Then the evolution semigroup $\mathcal R_\nu \auxMult_\mathcal T$ is unitarily equivalent to the tensor product semigroup $\mathcal R_\nu \otimes \mathcal T$, which acts on $L^2[0,\nu] \otimes \auxHilbert$. In fact, Proposition \ref{OT24quasi} is a particular case of the results obtained in \cite{oliva2025tensor}, where we prove that the tensor product of two $C_0$-semigroups $\mathcal T$ and $\mathcal S$ is similar to a contraction one if and only if there exist $\lambda \in \RR$ such that the rescaled semigroups $e_{-\lambda} \mathcal T$ and $e_\lambda \mathcal S$ are similar to contraction ones.
	
	 While our proofs of Theorems \ref{nearzero} and \ref{quasiContrTh_int}, as well as the ones of the auxiliary results of Section \ref{mainSection}, could be rephrased in terms of semigroup tensor products, 
		 we chose to work with 
	 evolution semigroups, which we find more transparent and instructive. 
	 \end{remark}

Now we are ready to give the proof of Theorem \ref{nearzero}.
\begin{proof}[Proof of Theorem \ref{nearzero}] $ $\newline
	Let $\mathcal T = (T(t))_{t\geq0}$ be a $C_0$-semigroup on a Hilbert space $H$.
	
	(i) $\implies$ (ii): Assuming $\mathcal T \in \mathcal{SQC}(H)$, we prove that for each $\nu>0$, there exist a Hilbert space $K$, operators $A \in\linearOp(K,H)$ and $B\in \linearOp(H,K)$, and a nilpotent  contraction $C_0$-semigroup $\mathcal N = (N(t))_{t\geq0}$ on $K$ satisfying
	\begin{equation}\label{InterpSmallT1}
		T(t) = AN(t)B, \qquad t \in [0,\nu].
	\end{equation}
	Fix $\lambda \in \RR$ and an equivalent Hilbertian norm $\|\cdot\|_{\mathscr H}$ on $H$ such that $\|T(t)\|_{\mathscr H} \leq e^{\lambda t}$
	for $t \geq 0$, and set $\mathscr H:=(H, \|\cdot\|_{\mathscr H})$. Let $\mathcal R_{\lambda, \nu+1} = (R_{\lambda, \nu+1}(t))_{t\geq0}$ be the  right shift semigroup on $L_{\lambda,\nu+1}^2(\mathscr H)$.
	Then 
\[
\|R_{\lambda, \nu+1}(t)\|_{L_\lambda^2(\mathscr H)} \leq e^{-\lambda t} \quad \text{for} \,\, t\geq0, \qquad \text{and}\qquad R_{\lambda, \nu+1}(t) = 0\quad \text{if}\,\, t\geq \nu+1.
\]
 Hence, we have
	$$\| R_{\lambda, \nu+1}(t) \auxMult_{T(t)}\|_{L_{\lambda,\nu+1}^2(\mathscr H)} = \|R_{\lambda, \nu+1}(t)\|_{L_{\lambda,\nu+1}^2(\mathscr H)} \|T(t)\|_{\mathscr H} \leq 1, \qquad t \geq 0,
	$$
	so that $\mathcal N := \mathcal R_{\lambda, \nu+1} \auxMult_\mathcal T = (R_{\lambda,\nu+1}(t) \auxMult_{T(t)})_{t\geq0}$,  given by \eqref{Evsemeq}, is a nilpotent contraction $C_0$-semigroup on $L_{\lambda,\nu+1}^2 (\mathscr H)$ vanishing on  $[\nu+1,\infty).$
In addition, note that, for all $h\in H$,
	\begin{equation}\label{QuasiAuxRemarkEq2}
		Q_{\lambda,\nu+1} R_{\lambda, \nu+1}(t) h\chi_{[0,1]} = \int_0^{\nu+1} h \chi_{[t, t+1]}(x), \,dx =  \begin{cases}
			h, \quad & \mbox{if } 0 \leq t\leq \nu,
			\\ (\nu+1 - t)h, \quad & \mbox{if } \nu \leq t \leq \nu+1,
			\\0, \quad & \mbox{if } t \geq \nu+1.
		\end{cases}
	\end{equation}
	Let now  $F_{\lambda, \nu+1} \in \linearOp(H, L_{\lambda, \nu+1}^2 (\mathscr H))$ be defined by
	\begin{equation}\label{gammah}
	F_{\lambda,\nu+1} h = h \chi_{[0,1]}, \qquad h \in H,
	\end{equation}
	where $h\chi_{[0,1]}$ is naturally regarded as an element of $L^2_{\lambda,\nu+1}(\mathscr H)$.
	Then taking into account \eqref{QuasiAuxRemarkEq2}, we conclude that
	\begin{equation}\label{Interp1}
		\begin{aligned}
			Q_{\lambda,\nu+1} N(t)
			F_{\lambda, \nu+1} h 
			&= Q_{\lambda,\nu+1} R_{\lambda, \nu+1} (t) \, T(t)h \chi_{[0,1]} 
			\\ &= \begin{cases}
				T(t)h, \quad & \mbox{if } t\leq \nu,
				\\ (\nu+1 - t) T(t) h, \quad & \mbox{if } \nu \leq t \leq \nu+1,
				\\0, \quad & \mbox{if } t \geq \nu+1,
			\end{cases}
					\end{aligned}
	\end{equation}
	for all $h \in H$ and $t \geq 0.$ Letting
	\begin{equation}
	A:= Q_{\lambda,\nu+1} \qquad \text{and}  \qquad B:= F_{\lambda, \nu+1},
		\end{equation}
	we infer that \eqref{InterpSmallT1} holds.
		
	Note that
	\begin{equation}\label{ABboundsEq}
		\begin{aligned}
			\|A\| &= \|Q_{\lambda,\nu+1}\| =  \left(\frac{e^{2\lambda(\nu+1)}-1}{2\lambda}\right)^{1/2},
			\\ \|B\| &= \|F_{\lambda, \nu+1}\| = \|\chi_{[0,1]}\|_{L_{\lambda, \nu+1}^2 (\CC)} =
			\left(\frac{1-e^{-2\lambda}}{2\lambda}\right)^{1/2}.
		\end{aligned}
	\end{equation}
	These bounds will be of importance in the proof of Theorem \ref{SimConstTh}.\medskip
	
	(ii) $\implies$ (iii): Let $\nu>0$ be fixed and let $K,\, A,\, B$ and $\mathcal N$ be as in \eqref{InterpSmallT1}. Invoking Nagy's dilation theorem we construct
	a Hilbert space $\mathcal H$ containing $K$ as a subspace,
	and a unitary dilation $\mathcal U  = (U(t))_{t\in \RR}$ of $\mathcal N$ on $\mathcal H.$
	 Set $\mathscr A := A\proj_K \in \linearOp(\mathcal H,H)$ and $\mathscr B := B \in \linearOp(H,\mathcal H).$
	 Since $\ran(B) \subseteq K$, one has
	\begin{equation}\label{InterpolationSmallT2}
		T(t) = AN(t)B = A\proj_K U(t)\restriction_K B = \mathscr A U(t) \mathscr B, \qquad t \in [0,\nu],
	\end{equation}
	as required.
		  \medskip

	(iii) $\implies$ (iv): This implication is direct.
		\medskip
	
	(iv) $\implies$ (i): 
	Suppose there exist  $\nu>0$, a Hilbert space $K$, operators $A \in\linearOp(K,H)$
	and $B\in \linearOp(H,K)$, and a semigroup $\mathcal S \in \mathcal{SQC}(K)$
	satisfying
		\begin{equation}\label{InterpSmallT3}
		T(t) = A S(t)  B, \qquad t \in [0,\nu].
	\end{equation}
		Up to passing to an equivalent Hilbertian norm on $K$,
		we may assume that $\mathcal S$ is a quasi-contraction $C_0$-semigroup. To simplify our presentation, with an abuse of notation, we let $\mathcal R_\nu = (R_\nu(t))_{t\geq0}$  stand for both the right shift semigroup on $L_\nu^2(K)$ and on $L_\nu^2(H)$. Consider the evolution semigroup $\mathcal R_\nu \auxMult_{\mathcal T} = (R_\nu(t) \auxMult_{T(t)})_{t\geq0}$ on $L^2_\nu(H)$. As $R_\nu(t) = 0$ for $t \geq \nu$, we have
	\begin{equation}\label{LebowTensorEq}
		\begin{aligned}
		R_\nu(t) \auxMult_{T(t)}  =& \begin{cases}
			R_\nu(t) \auxMult_{ A S(t) B}, &\mbox{ if } 0\leq  t \leq \nu, 
			\\0, & \mbox{ if } t \geq \nu,
		\end{cases} 
		\\=& R_\nu(t) \auxMult_{ A S(t) B} =  \auxMult_{A}  R_\nu(t) \auxMult_{S(t)}  \auxMult_{ B}, \qquad t \geq 0.
		\end{aligned}
	\end{equation}
	Since $\mathcal R_\nu \auxMult_{\mathcal S} = (R_\nu(t) \auxMult_{S(t)})_{t\geq0}$ is a quasi-contraction $C_0$-semigroup,
	Proposition \ref{compressionProp} implies that $\mathcal R_\nu \auxMult_{\mathcal T} \in \mathcal{SQC}
	\left(L_\nu^2(H)\right)$. Then, by Proposition \ref{OT24quasi}, we infer that $\mathcal T \in \mathcal{SQC}(H).$
				\end{proof}
	\begin{remark}
			The parameters $\lambda$ and $\nu$ introduced in the proof above will play a role in the proof of Theorem \ref{SimConstTh}. In particular, $\lambda$ will be used in $\lambda$-dependent estimates for similarity constants, while $\nu$ will allow to consider simultaneously  $Q_{\lambda, \nu}$ and $F_{\lambda, \nu}$ for different $\nu.$
					\end{remark}

\begin{remark}
		Theorem \ref{nearzero} has no discrete counterpart. 
	In fact, for an arbitrary operator $T$ on a Hilbert space $H$ and $N \in \NN$, there exist a Hilbert space $K$, bounded operators $B \in \linearOp(H, K)$ and $A \in \linearOp(K, H)$, and a contraction operator $S$ on $K$ such that
	\begin{align*}
		T^n = A S^n B, \qquad n = 1, \ldots, N,
	\end{align*}
	see, for example, \cite[p.910]{halmos1970ten} and \cite[Sect. 4]{holbrook1973operators}.
\end{remark}

Next we prove Corollary \ref{C0groupdilation_int}, which has the same flavor as the one of \cite[Theorem 7.3]{batty2017lower} for left-invertible semigroups.

\begin{proof}[Proof of Corollary \ref{C0groupdilation_int}] \,\,
	(i) $\implies$ (ii): Letting $\mathcal T = (T(t))_{t\geq0}$ be a $C_0$-semigroup on a Hilbert space $H$, we assume that $\mathcal T \in \mathcal{SQC}(H)$, so there exist a bounded and invertible operator $R$ and $\lambda \in \RR$ such that $\|R^{-1} T(t) R\| \leq e^{\lambda t}$ for all $t\geq 0$. By applying Nagy's dilation theorem in the form \eqref{dil_isom}
to $(e^{-\lambda t} R^{-1} T(t) R)_{t\ge0}$, we infer that there exist a Hilbert space $K$, an isometry $V\in \linearOp(H, K)$ and a unitary $C_0$-group $\mathcal U =(U(t))_{t\in \RR}$ on $K$ such that
	$$T(t) = e^{\lambda t} R V^\ast U(t) V R^{-1}, \qquad t \geq 0.
	$$
	Set $Q := VRV^\ast + I - VV^\ast \in \linearOp(K)$. Since $V^\ast V = I$, it is readily seen that $Q$ is invertible with
	$$Q^{-1} = VR^{-1}V^\ast + I - VV^\ast.
	$$
	Define now $\mathcal G = (G(t))_{t\in \RR}$ as $G(t) = Q e^{\lambda t} U(t) Q^{-1}$ for $t\in \RR$. Then $\mathcal G$ is a $C_0$-group on $K$ that is similar to $(e^{\lambda t} U(t))_{t\in \RR}.$ So observing that $V^\ast Q = RV^\ast$ and $Q^{-1} V = VR^{-1}$, we obtain
	\begin{align}
		 T(t) = e^{\lambda t} R V^\ast U(t) V R^{-1} = V^\ast Q e^{\lambda t} U(t) Q^{-1} V =  V^\ast G(t) V, \qquad t \geq 0,
	\end{align}
	as required. \medskip
		
	(ii) $\implies$ (iii): This implication is trivial consequence of the identity
	$T(t)=V^*G(t)V, t \ge 0,$ for a $C_0$-semigroup $\mathcal T=(T(t))_{t \ge 0}$ and its dilation $G=(G(t))_{t \ge 0}.$
	\medskip
	
	(iii) $\implies$ (i): Let  $K$ and $H$ be Hilbert spaces,  $\mathcal T$ be a $C_0$-semigroup on $H$,
	$\mathcal S$ be a $C_0$-semigroup on $K$ bounded from below, and let $T(t)=V^*S(t)V$ for
	$t \in [0,\nu]$ and an isometry
	$V: H \to K.$ Then Proposition \ref{boundedBelowContractProp}(ii)	yields $\mathcal S \in \mathcal{SQC}(K),$
and the claim follows from the implication (iv) $\implies$ (i) of Theorem \ref{nearzero}.
\end{proof}

\section{Similarity through quasi-similarity: putting the two sides together}\label{mainSection}

After analyzing the similarity properties of $\mathcal T =(T(t))_{t \ge 0}$ for small and large $t$ we are able to prove one of the main results of the paper, Theorem \ref{quasiContrTh_int}, characterizing $\mathcal T \in \mathcal{SC} (H),$ or, in other words, joint similarity of $(T(t))_{t \ge 0}$ to contractions for all $t \ge 0.$ The two cases are linked to each other through the notion of quasi-contractivity which, being a much weaker property than contractivity, reveal new effects in the studies of similarity to contractions in the continuous setting.

In fact,
 we prove a more general statement by providing an upper bound on the similarity constant $C(\mathcal T)$ of $\mathcal T,$
which is of independent interest. While there are very few  works
addressing quantitative aspects of similarity constants,
the importance
of such bounds was emphasized and clarified, in particular, in \cite{badea2003operators},
\cite{holbrook1977distortion} and \cite[Chapter 9]{pisier2001similarity}.

\begin{theorem}\label{SimConstTh}
	Let $\mathcal T = (T(t))_{\geq0}$ be a $C_0$-semigroup on a Hilbert space $H$. Assume that there exist $\lambda >0$ and $\tau>0$ such that
	\begin{enumerate}
		\item [(i)] $e_{-\lambda} \mathcal T \in \mathcal{SC}(H)$;
		\item [(ii)] $T(\tau)$ is similar to a contraction.
	\end{enumerate}
	Then $\mathcal T \in \mathcal{SC}(H).$
	Moreover, the similarity constant $\mathcal C(\mathcal T)$ of $\mathcal T$ satisfies
	\begin{equation}\label{boundSimEq}
	\mathcal C(\mathcal T) \leq  \sqrt{2} \, \mathcal C(e_{-\lambda} \mathcal T) \frac{e^{2\lambda}-1}{2\lambda} + 2 \sqrt{2} \, C(T(\tau)) M^2 \max \{1,\sqrt{\tau}\},
	\end{equation}
	where $M = \sup_{t\in [0,\tau]} \|T(t)\|$.
\end{theorem}

Note that if $e_{-\lambda} \mathcal T \in \mathcal{SC}(H)$ for some $\lambda \leq 0$, then clearly $\mathcal T \in \mathcal{SC}(H)$ and $\mathcal C(\mathcal T) \leq \mathcal C (e_{-\lambda} \mathcal T)$.
{ We also emphasize that condition (i) concerns similarity properties of $\mathcal T$ for small times, 
while condition (ii) does so for large times.}

To prove Theorem \ref{SimConstTh} (and thus Theorem \ref{quasiContrTh_int}, since (i) and (ii) are clearly necessary), we need to construct an auxiliary semigroup $\mathcal W$ on $L^2([0,1], \auxHilbert) \oplus L^2([0,1], \auxHilbert)$ that is similar to a contraction semigroup, where $\auxHilbert$ is a given Hilbert space from now on in this section.
The semigroup $\mathcal W$ will serve to link our constructions of intertwining relations for $\mathcal T$ near zero and away from zero.

We begin with the several definitions and a technical lemma.
Define  the unitary (periodic) $C_0$-group $\mathcal R_p = (R_p(t))_{t\in \RR}$ on $L^2([0,1], \auxHilbert)$ as
\begin{equation}\label{defRp}
(R_p(t)f)(x) = f([x-t] \operatorname{mod} 1), \qquad x \in [0,1), \, t \in \RR, \, f \in L^2([0,1], \auxHilbert).
\end{equation}
In addition, for each $t\geq0$, consider the family $\mathcal V=(V(t))_{t \ge 0} \subset \mathcal L(L^2([0,1],\auxHilbert))$ given by
\begin{equation}\label{defp}
(V(t)g)(x) = \begin{cases}
	0, \quad &x \geq t,
	\\ g\left(\left[x - t\right]  \operatorname{mod} 1\right) , \quad &x < t,
\end{cases}
\qquad \qquad x \in [0,1), \, t \geq 0, \, g \in L^2([0,1],\auxHilbert).
\end{equation}
Clearly,
$V(t)g$ tends strongly to $0$
for every $g \in L^2([0,1], \auxHilbert)$ as $t\to 0$.

\begin{lemma}\label{PIdentityLemma}
	If $\mathcal R = (R(t))_{t\geq 0}$ is the right shift semigroup on $L^2([0,1], \auxHilbert),$ and the families of linear bounded operators $\mathcal R_p$ and $\mathcal V$ on $L^2([0,1],\auxHilbert)$ are given by \eqref{defRp} and \eqref{defp} respectively, then for all $s,t\geq 0$, 
	\begin{equation}\label{almsem}
	V(s+t) = R_p(s) V(t) + V(s) R(t).
	\end{equation}
\end{lemma}
\begin{proof}
		For $x \in [0,1)$ and $s \geq 0$, define $a_{x,s} := \left[x-s\right]  \operatorname{mod} 1\in [0,1)$. Then, for all $s,t\geq0$, $g \in L^2([0,1], \auxHilbert)$ and $x \in [0,1)$,
	\begin{align*}
		(R_p(s)V(t)g)(x) &= (V(t)g) (a_{x,s}) = \begin{cases}
			0, \quad &a_{x,s} \geq t,
			\\g\left(\left[ a_{x,s} - t\right] \operatorname{mod} 1\right),  \quad & a_{x,s} < t,
		\end{cases}
	\end{align*}
	and
	\begin{align*}
		(V(s)R(t)g)(x) &= \begin{cases}
			0, \quad & x \geq s,
			\\ (R(t)g)(a_{x,s}), \quad & x < s,
		\end{cases}
		\\ &= \begin{cases}
			0, \quad & x \geq s \,\mbox{ or }\, a_{x,s} < t,
			\\ g(a_{x,s}-t), \quad & x < s\, \mbox{ and } \, a_{x,s} \geq t.
		\end{cases}
		\\ &= \begin{cases}
			0, \quad & x \geq s \,\mbox{ or }\, a_{x,s} < t,
			\\ g([a_{x,s}-t]  \operatorname{mod} 1), \quad & x < s\, \mbox{ and } \, a_{x,s} \geq t,
		\end{cases}
	\end{align*}
	where, in the last step, we used that $a_{x,s}-t \in [0,1)$ if $a_{x,s}\geq t$, so $[a_{x,s}-t]  \operatorname{mod} 1 = a_{x,s}-t$. Furthermore, one has
	$x \geq s +t$ with $x \in [0,1]$ and $s,t\geq0$ if and only if $x \geq s$ and $a_{x,s} \geq t$ (note that $x \geq s$ implies $a_{x,s} = x - s$). Since,
	 $[x+(y  \operatorname{mod} 1)]  \operatorname{mod} 1= [x+y]  \operatorname{mod} 1$ for all $x,y\in \RR$,
		we conclude that
	\begin{align*}
		(R_p(s)V(t)g + V(s)R(t)g)(x) &= \begin{cases}
			0,\quad & x \geq s \, \mbox{ and } \, a_{x,s} \geq t,
			\\g([a_{x,s}-t] \operatorname{mod} 1), \quad & x < s \, \mbox{ or } \, a_{x,s}<t,
		\end{cases}
		\\ &= \begin{cases}
			0,\quad & x \geq s +t,
			\\g\left(\left[x-s-t\right] \operatorname{mod} 1\right), \quad & x < s+t,
		\end{cases}
		\\ &= (V(s+t)g)(x),
	\end{align*}
	i.e. \eqref{almsem} holds.
	\end{proof}

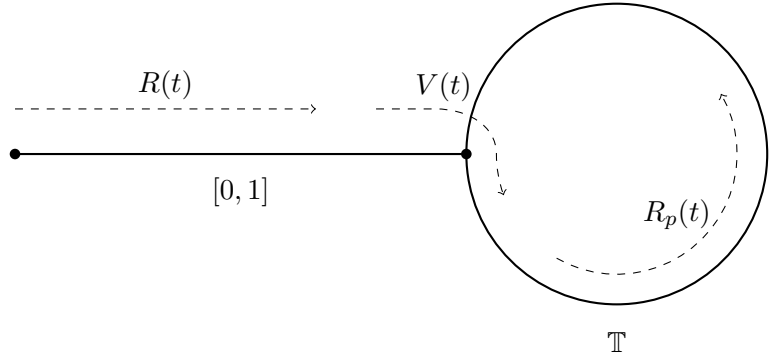
\begin{figure}
	\begin{center}
		
		\begin{tikzpicture}
			\pgfmathsetmacro\scale{2}

			\draw[thick] (-2*\scale,0) -- (1*\scale,0);
			\node at (-0.5*\scale, -0.25*\scale) {$[0,1]$};
			
			\fill (-2*\scale,0) circle (2pt);
			\fill (1*\scale,0) circle (2pt);
			
			\draw[thick] (2*\scale,0) circle (1*\scale);
			\node at (2*\scale,-1.25*\scale) {$\TT$};
			
			\draw[dashed,->] (-2*\scale,0.3*\scale) -- (0*\scale,0.3*\scale) node[midway, above] {$R(t)$};
			\draw [dashed,->,domain=240:390] plot ({2*\scale+0.8*\scale*cos(\x)}, {0.8*\scale*sin(\x)});
			\node at ({2*\scale+0.4*\scale},-0.4*\scale) {$R_p(t)$};

			\draw[dashed,->] (0.4*\scale, 0.3*\scale) -- (0.8*\scale,0.3*\scale) to[out=0, in=90] (1.2*\scale,0) arc(180:200:0.8*\scale);
			\node at (0.85*\scale, 0.45*\scale) {$V(t)$};
		\end{tikzpicture}
	\caption{Visualization of the action of the semigroup $\mathcal W$ on $L^2_\oplus(\auxHilbert) \simeq L^2(\TT, \auxHilbert) \oplus L^2([0,1],\auxHilbert)$, illustrating how it ``shifts'' the support of functions.}
	\label{WFigure}
	\end{center}
\end{figure}

Given a Hilbert space $\auxHilbert$, we set
$$L^2_\oplus (\auxHilbert) := L^2([0,1], \auxHilbert) \oplus L^2([0,1], \auxHilbert),
$$
equipped with a natural, direct sum Hilbertian norm $\|\cdot\|_{\oplus}$ and, for each $t\geq0$, define the operator $W(t) \in \linearOp(L^2_\oplus (\auxHilbert))$ by
\begin{align}\label{WDef}
	W(t): = \begin{pmatrix}
		R_p(t) & V(t)
		\\ 0 & R(t)
	\end{pmatrix},
\end{align}
see Figure \ref{WFigure} for an illustration of $\mathcal W = (W(t))_{t\geq0}.$

Now we are ready to show that $\mathcal W$ is indeed similar to a contraction semigroup.
\begin{lemma}\label{WcontractLemma}
	Let $\auxHilbert$ be a Hilbert space. If $\mathcal W = (W(t))_{t\geq0}$ is given by \eqref{WDef}, then $\mathcal W$ is a
	$C_0$-semigroup  and $\mathcal W \in \mathcal{SC}\left(L^2_\oplus (\auxHilbert)\right)$.
\end{lemma}
\begin{proof}
	It follows from Lemma \ref{PIdentityLemma} that $\mathcal W$ satisfies the semigroup law, i.e., $W(s+t) = W(s)W(t)$ for $s,t\geq0$. The strong continuity of $\mathcal W$ follows immediately from the strong continuity of $\mathcal R_p$, $\mathcal V$, and $\mathcal R$. It remains to prove that $\mathcal W \in \mathcal{SC}\left(L^2_\oplus (\auxHilbert)\right).$
		
	Note first that $R_p(1) = V(1) = I,$
	and
	$$R_p(1)V(t) + V(1)R(t) =
		V(t+1) = R_p(t)V(1) + V(t) R(1) = R_p(t), \qquad t \geq 0,
	$$
	by Lemma \ref{PIdentityLemma}. Hence
	\begin{align}\label{ident}
		R_p(t)f + V(t)g + R(t)g = R_p(t)f + (R_p(1)V(t) + V(1)R(t))g = R_p(t)(f + g),
	\end{align}
	for all $f, \, g \in L^2([0,1], \auxHilbert)$ and $t \geq 0.$

	Now let $\|\cdot\|_{\oplus, {\rm eq}}$ be the equivalent Hilbertian norm on $L^2_\oplus(\auxHilbert)$ induced by the invertible operator $\Lambda_I:= \begin{pmatrix}I & I \\0 & I \end{pmatrix}$, i.e.,
	\begin{align}\label{SimDefEq}
		\|(f,g)\|_{\oplus, {\rm eq}}^2 := \|\Lambda_I (f,g)\|^2_{\oplus} = \|f+g\|^2 + \|g\|^2, \qquad (f,g) \in L^2_\oplus(\auxHilbert).
	\end{align}
	
	 Since $\mathcal R_p$ is  unitary and $\mathcal R$ is contractive,
		taking into account \eqref{ident}, we obtain
	\begin{align*}
		\|W(t)(f,g)\|_{\oplus, {\rm eq}}^2 &= \|(R_p(t)f + V(t)g, R(t)g)\|_{\oplus, \rm{eq}}^2
		\\ &= \|R_p(t) f + V(t)g + R(t)g\|^2 + \|R(t)g\|^2
		\\ &= \|R_p(t) (f + g)\|^2 + \|R(t)g\|^2
		\\ &\leq \|f + g\|^2 +  \|g\|^2 = \|(f,g)\|_{\oplus, {\rm eq}}^2, \qquad (f, g) \in L^2_\oplus(\auxHilbert), \, t \geq 0.
	\end{align*}
	So setting
	\begin{equation}\label{space}
	L^2_{\oplus,\rm{eq}}([0,1],\auxHilbert):=\left(L^2_\oplus(\auxHilbert), \|\cdot\|_{\oplus, {\rm eq}}\right),
	\end{equation}
	it follows that
	$\mathcal W$ is a contraction $C_0$-semigroup on $L^2_{\oplus, \rm{eq}}(\auxHilbert)$,
		which implies the claim.
\end{proof}

We will need one more identity involving the operators $(V(t))_{t \ge 0}.$

Recalling the definition of the family of operators $Q_{\lambda, \nu}$ in \eqref{QuasiAuxRemarkEq1bis},
consider the operator $Q=Q_{0,1}$ from $L^2([0,1],\auxHilbert)$ to $\auxHilbert$ given by
\begin{align}\label{QOpDef}
	Qf = \int_0^1 f(x)\, dx, \qquad f \in L^2([0,1], \auxHilbert).
\end{align}

\begin{lemma}\label{IntQLemma}
	Let $\auxHilbert$ be a Hilbert space. If $\mathcal V=(V(t))_{t \ge 0}$ and $Q$ are given by \eqref{defp} and \eqref{QOpDef} respectively, then for all $t\geq0$ and all $\auxVector\in \auxHilbert$, one has
	$$Q V(t) \auxVector \chi_{[0,1/2]} = \begin{cases}
		0, \quad & 0\leq t \leq 1/2,
		\\ (t-\frac{1}{2}) \auxVector, \quad & 1/2\leq t \leq 1,
		\\ \frac{1}{2}\auxVector, \quad & t \geq 1.
	\end{cases}
	$$
\end{lemma}
\begin{proof}
	Assume first $t\geq 1.$ Then by Lemma \ref{PIdentityLemma},
	$$
	V(t) = V(t-1+1) = R_p(t-1) V(1) + V(t-1)R(1) =  R_p(t-1),
	$$
	as $R(1)=0$ and $V(1) = I$. Also, a simple change of variable shows that $QR_p(s) =Q$ for every $s\geq0$. Hence, for each $\auxVector \in \auxHilbert$,
	\begin{align*}
		Q V(t) \auxVector \chi_{[0,1/2]} = Q R_p(t-1) \auxVector\chi_{[0,1/2]} = Q \auxVector\chi_{[0,1/2]} =\int_0^{1} h \chi_{[0,1/2]}(x) \, dx =  \frac{1}{2}\auxVector.
			\end{align*}
		
	Let now $t \in [0,1]$. In this case, $x \in [0,t)$ implies that $x - t \in [-1,0)$, thus $\left[x - t\right]  \operatorname{mod} 1= 1+ x - t$. Then, for each $\auxVector \in \auxHilbert$,
	\begin{align*}
		V(t) \auxVector \chi_{[0,1/2]}(x)  &=  \begin{cases}
			0, \quad &x \geq t,
			\\ \auxVector \chi_{[0,1/2]}\left(\left[x - t\right] \operatorname{mod} 1  \right), \quad & x < t,
		\end{cases}
		\\ &=  \begin{cases}
			0, \quad &x \geq t,
			\\ \auxVector \chi_{[0,1/2]}\left(1+ x - t\right), \quad & x < t,
		\end{cases}
		\\ &=  \auxVector \chi_{[t-1, t-1/2]} (x),
		\\ &= \auxVector \chi_{[0, \max\{0,t-1/2\}]}(x), \qquad x \in [0,1).
	\end{align*}
	Hence, for each $t \in [0,1]$ and $\auxVector\in \auxHilbert$,
	\begin{align*}
		Q V(t)  \auxVector \chi_{[0,1/2]} = Q  \auxVector\chi_{[0, \max\{0,t-1/2\}]} = \int_0^1 \auxVector \chi_{[0, \max\{0,t-1/2\}]} (x) \, dx = \max\{0,t-1/2\} \auxVector,
			\end{align*}
	and the statement follows.
\end{proof}

We are now ready to prove Theorem \ref{SimConstTh} (and thus Theorem \ref{quasiContrTh_int}).

\begin{proof}[Proof of Theorem \ref{SimConstTh}]
	Let $\mathcal T = (T(t))_{t\geq0}$ be a $C_0$-semigroup on a Hilbert space $H$ and assume that there exist $\lambda > 0$ and $\tau>0$ such that $e_{-\lambda} \mathcal T$ is in $\mathcal{SC}(H)$ and  $T(\tau)$ is similar to a contraction on $H$. We show that $\mathcal T \in \mathcal{SC}(H)$ and that the similarity constant $\mathcal C(\mathcal T)$ satisfies \eqref{boundSimEq}.
	To this aim we apply Holbrook-type criterion and construct a Hilbert space $K,$ operators $A \in \mathcal L(H, K)$ and $B \in \mathcal L(K, H)$
	and a contraction semigroup $\mathcal S=(S(t))_{t \ge 0}$ on $K$ such that
	$T(t)=AS(t)B$ for all $t \ge 0.$ 
	The construction is divided into two parts, treating separately the cases where $t$ is close to zero and where $t$ is bounded away from zero, and these are then merged via the auxiliary semigroup $\mathcal W$ defined in \eqref{WDef}.
		\medskip

	First, starting from $t$ near zero, we will employ an argument
	used for the proof of Theorem \ref{nearzero}.
	Let $\|\cdot\|_1$ be an equivalent Hilbertian norm on $H$ satisfying $\|T(t)\|_1\leq e^{\lambda t}$ for all $t\geq0$ and
	\begin{align*}
		\| \cdot \| \leq \|\cdot\|_1 \leq \mathcal C(e_{-\lambda} \mathcal T) \|\cdot\|.
	\end{align*}
	Recalling the definition of spaces $L^2_{\lambda, \nu}$
	in \eqref{QuasiAuxRemarkEq1} and of operators $Q_{\lambda,\nu}$ in \eqref{QuasiAuxRemarkEq1bis}, fix $\nu=2$ and consider
	$$
	K_1 := L_{\lambda,2}^2(H_1),
	$$
	where $H_1=(H, \|\cdot\|_1)$.
	Define $A_{1} \in \linearOp(K_1, H),$  $B_{1} \in \linearOp(H, K_1)$, and the
	$C_0$-semigroup $\mathcal S_{1} = (S_{1}(t))_{t\geq0}$ on $K_1$  by
	\begin{equation}\label{a1b1s1}
		\begin{aligned}
		A_{1} f &:= Q_{\lambda,2}f =\int_{0}^{2}f(x)\, dx, \qquad f \in K_1,
		\\ B_{1} h &:= F_{\lambda,2}h= h \chi_{[0,1]}, \qquad h \in H, 
		\\ S_{1}(t) &= R_{\lambda,2}(t) \auxMult_{T(t)}, \qquad t\geq0, 
		\end{aligned}
	\end{equation}
	where $\mathcal R_{\lambda, 2} = (R_{\lambda, 2}(t))_{t\geq0}$ is the right shift semigroup on $K_1$. 
	Since $\|R_{\lambda, 2}(t)\| \leq e^{-\lambda t}$ for all $t\geq0$,
	we conclude that $\mathcal S_{1}$ is a contraction $C_0$-semigroup on $K_1$ and,
	in view of \eqref{Interp1},
	\begin{equation}\label{InterpEq1}
		A_{1} S_{1}(t) B_1 = \begin{cases}
			T(t), \quad & 0\leq t \leq 1,
			\\ (2-t) T(t), \quad & 1\leq t \leq 2,
			\\ 0, \quad & t \geq 2.
		\end{cases}
	\end{equation}
	In addition, by \eqref{ABboundsEq},
	\begin{equation}\label{AB1BoundsEq}
		\|B_1\|
		= \left(\frac{1-e^{-2\lambda}}{2\lambda}\right)^{1/2}, \qquad \|A_{1} \| \leq \mathcal C (e_{-\lambda} \mathcal T) \left(\frac{e^{4\lambda}-1}{2\lambda}\right)^{1/2}.
	\end{equation}
		
	Second, using Theorem \ref{awayzero}, we infer that there exist a Hilbert space $K_0$,
	operators $A_0 \in \linearOp(K_0, H),$ $B_0 \in \linearOp(H, K_0)$, and a contraction $C_0$-semigroup $(S_0(t))_{t\geq0}$ on $K_0$ such that
	$$
	T(t) = A_0 S_0(t) B_0, \qquad t \geq 1.
	$$
	By applying Theorem \ref{awayzeroSCboundTh} with $\tau_1 = \tau$ and $\tau_2 = \min\{1, \tau\}$, we can assume that
	\begin{equation}\label{badABbound}
		\|B_0\| \leq 1, \qquad \|A_0\| \leq
		C(T(\tau)) M^2 \max \{1,\sqrt{\tau}\},
	\end{equation}
	where $M=\sup_{t\in [0,\tau]} \|T(t)\|$.
	Consider now the semigroup $\mathcal W=(W(t))_{t \ge 0}$ and the Hilbert space
	 \[
	K_2 := L^2_{\oplus, \rm{eq}}(K_0)
	\]
	 given by
	\eqref{WDef} and \eqref{space}, respectively. 
	For any $T \in \linearOp (H)$ define a bounded operator $\auxMult^{\oplus}_T$ 
	on $L^2_{\oplus}(K_0)$ with any equivalent norm by
	\[
	\auxMult^\oplus_T (f,g)=(\auxMult_T f, \auxMult_T g), \qquad (f,g) \in L^2_{\oplus}(K_0).
	\]
	
	Then, observing that the semigroups $\mathcal M^{\oplus}_{\mathcal S_0}$
	and $\mathcal W$ commute, let the $C_0$-semigroup  $\mathcal S_2 = (S_2(t))_{t\geq0}$ on $K_2$ be given by 
		$$
	S_2(t):= \auxMult_{S_0(t)}^{\oplus} W(t/2), \qquad t\geq0.
	$$
	Lemma \ref{WcontractLemma} implies that
	$\mathcal S_2$ is a contraction semigroup on $K_2$.
	Finally, define the operators $A_2 \in \linearOp(K_2, H)$ and $B_2 \in \linearOp(H, K_2)$  by
	\begin{align*}
		A_2 (f,g) &:= 2A_0 \int_0^1 f(x) \, dx, \qquad (f,g) \in L^2_{\oplus, \rm{eq}}(K_0), 
		\\ B_2 h &:= \begin{pmatrix} 0 \\ B_0 h \chi_{[0,1/2]}  \end{pmatrix}, \qquad h \in H,
	\end{align*}
	so that
	$$A_2 = \begin{pmatrix} 2  A_0 Q   & 0 \end{pmatrix} = \begin{pmatrix} 2  Q  & 0 \end{pmatrix}  \auxMult^{\oplus}_{A_0}
	$$
	(where $Q$ is defined by \eqref{QOpDef}).
	Now, taking into account Lemma \ref{IntQLemma} and the fact that $\mathcal W$ commutes with the multiplication operators $\auxMult_{A_0}^{\oplus}$ and $\auxMult_{S_0(t)}^{\oplus}$ for all $t \ge 0$, we infer that
	\begin{equation}\label{AnotherAuxEq}
		\begin{aligned}
			A_2 S_2(t) B_2 h &=  \begin{pmatrix} 2  Q  & 0 \end{pmatrix} \auxMult_{A_0}^{\oplus}  \auxMult_{S_0(t)}^{\oplus} W(t/2)  \begin{pmatrix} 0 \\ B_0 h \chi_{[0,1/2]}  \end{pmatrix}  
			\\ &= \begin{pmatrix} 2 Q  & 0 \end{pmatrix} W(t/2)  \begin{pmatrix} 0 \\ A_0 S_0(t) B_0 h \chi_{[0,1/2]}  \end{pmatrix} 
			\\ & =  2  Q V(t/2) \,  A_0 S_0(t) B_0 h \chi_{[0,1/2]}
			\\ &= \begin{cases}
				0, \quad & 0\leq t \leq 1,
				\\ (t-1) A_0 S_0(t) B_0 h, \quad & 1\leq t \leq 2,
				\\ A_0 S_0(t) B_0 h, \quad & t \geq 2, 
			\end{cases}
			\\ &= \begin{cases}
				0, \quad & 0\leq t \leq 1,
				\\ (t-1) T(t)h, \quad & 1\leq t \leq 2,
				\\ T(t)h, \quad & t \geq 2,
			\end{cases} \qquad h \in H, \, t \geq 0.
		\end{aligned}
	\end{equation}
	\medskip

	Finally, let
	$$K := K_1 \oplus K_2, \qquad \qquad S(t) := S_{1}(t) \oplus S_2(t), \quad t\geq0,
	$$
	so that $\mathcal S = (S(t))_{t\geq 0}$ is a contraction $C_0$-semigroup on $K$
	by our construction.
	If
	\[
	A: = \begin{pmatrix}A_{1} & A_2\end{pmatrix} \qquad \text{and} \qquad B: = \begin{pmatrix}B_1 \\ B_2\end{pmatrix},
	\]
	then $A \in \linearOp(K,H)$ and $B \in \linearOp(H, K),$ and  from \eqref{InterpEq1} and \eqref{AnotherAuxEq} it follows that
	\begin{align*}
		A S(t) B &=  \begin{pmatrix}A_{1} & A_2\end{pmatrix} \begin{pmatrix}S_{1}(t) & 0 \\ 0 & S_2(t)\end{pmatrix} \begin{pmatrix}B_1 \\ B_2\end{pmatrix}
		\\ &= A_{1} S_{1}(t) B_1 + A_2 S_2(t) B_2 = T(t), \qquad t \geq0,
	\end{align*}
	as required in \eqref{holbeq}.
	Note that the above equality holds if $A_{1}$ and $B_1$ are replaced by $\gamma A_{1}$ and $\frac{1}{\gamma} B_1$ for any $\gamma >0,$ respectively. Thus, Proposition \ref{compressionProp} yields $\mathcal T \in \mathcal{SC}(H)$ with
	\begin{align}\label{finalIneqSim}
		\mathcal C (\mathcal T) \leq \|A\|\|B\| &= \inf_{\gamma>0} \left\{ \left( \gamma^2 \|A_{1}\|^2 + \|A_2\|^2 \right)^{1/2}  \left( \gamma^{-2} \|B_1\|^2 + \|B_2\|^2 \right)^{1/2} \right\}
		\\ &= \|A_{1} \| \|B_1\|  + \|A_2\| \|B_2\|,\notag
	\end{align}
	where the infimum is attained at
	$$\gamma = \sqrt{\frac{\|A_2\| \|B_1\| }{\|A_1\| \|B_2\|}}.
	$$
	Here and in the sequel of this proof, for convenience,
	we omit subscripts in operator norms of $A_i$'s and $B_i$'s, as this will not cause confusion. 
	
	It remains to establish suitable upper bounds for the norms of $A_{1}, \, A_2, \, B_1$ and $B_2$.
	By \eqref{AB1BoundsEq},
	\begin{equation}\label{a1b1}
		\|A_{1} \| \|B_1\| \leq \mathcal C (e_{-\lambda} \mathcal T) \left(\frac{e^{4\lambda}-1}{2\lambda}\right)^{1/2} \left(\frac{1-e^{-2\lambda}}{2\lambda}\right)^{1/2} \leq  \sqrt{2} \mathcal C (e_{-\lambda}\mathcal T) \frac{e^{2\lambda}-1}{2\lambda}.
	\end{equation}
	Moreover, since in view of \eqref{badABbound} we have $\|B_0\|\leq 1$,
	\begin{align*}
		\|B_2 h \|_{K_2} & = \left\| \left(0,  B_0 h \chi_{[0,1/2]}  \right)\right\|_{K_2} 
		\leq   \left\| \left(0,h \chi_{[0,1/2]}  \right) \right\|_{L_{\oplus, \rm{eq}}^2(H)}  =
		\sqrt{2}  \left\|h \chi_{[0,1/2]}  \right\|_{L^2([0,1], H)}   = \|h\| 
	\end{align*}
	for all $h\in H$. Hence
	\begin{equation}\label{b2}
		\|B_2 \| \leq 1.
	\end{equation}
	To obtain norm bounds for $A_2,$ observe that if $\embedding \in \linearOp(L^2_\oplus(K_0), L^2_{\oplus, \rm{eq}}(K_0))$ is the natural embedding, then, recalling \eqref{SimDefEq}, we have 
	\[
	\|\embedding (f,g)\|_{L^2_{\oplus,\rm{eq}}(K_0)}=
	\|\Lambda_I(f,g)\|_{L^2_{\oplus}(K_0)},
	 \qquad (f,g) \in L_\oplus^2(K_0).
	\] 	
	Since $\|Q\|=1$, we then obtain
	\begin{align*}
		\left\|\begin{pmatrix} Q & 0 \end{pmatrix} \right\|_{\linearOp(K_2, K_0)} &= \left\|\begin{pmatrix} Q & 0 \end{pmatrix} \Lambda_I^{-1}\right\|_{\linearOp(L_{\oplus}^2(K_0) , K_0)} 
				 = \sqrt{2}\|Q\| = \sqrt{2}.
	\end{align*}
	Hence, employing the bound for $\|A_0\|$ in \eqref{badABbound},
	we conclude that
	\begin{align}\label{a2}
		\|A_2 \| \leq &  2 \left\|\begin{pmatrix} Q & 0 \end{pmatrix} \right\|_{\linearOp(K_2, K_0)} \|A_0\|  
		\leq 2 \sqrt{2} \, C(T(\tau)) M^2 \max \{1,\sqrt{\tau}\},
	\end{align}
	Combining  \eqref{a1b1}, \eqref{b2} and \eqref{a2} to estimate the right hand side of \eqref{finalIneqSim}, we obtain \eqref{boundSimEq}, and thus finish the proof.
\end{proof}
\begin{remark}\label{rough}

We have not aimed for the sharpest estimate in \eqref{boundSimEq}, and the bound may be quite rough; for instance, its right-hand side cannot be smaller than $3\sqrt{2}$. Importantly, it allows us to control 
$\mathcal{C}(\mathcal{T})$ in terms of $\mathcal{C}(T(\tau))$, $\mathcal{C}(e_{-\lambda}\mathcal{T})$, and $\sup_{t \in [0,\tau]} \|T(t)\|$. This will play a crucial role in the study of similarity
of infinite tensor products of semigroups  to semigroups of contractions in \cite{oliva2025tensorbis}.
\end{remark}

The construction of the equivalent Hilbertian norm in Theorem \ref{SimConstTh}, which makes the semigroup
$\mathcal T$ contractive, is somewhat implicit and and involves several intermediate steps.
Therefore, it is of interest to provide an explicit expression for this norm, making it more suitable for further use.
The following formula serves this purpose and reveals finer details of the renorming employed in Theorem \ref{SimConstTh}.

Let $H$ be a Hilbert space and let $\mathcal T = (T(t))_{t\geq0}$ be a $C_0$-semigroup on $H$ such that $\mathcal T \in \mathcal{SQC}(H)$ and  $T(\tau)$ is similar to a contraction operator for some $\tau>0$. Let $\lambda > 0$ and let $\|\cdot\|_1$ be an equivalent Hilbertian norm on $H$ such that
$$\|T(t)\|_1 \leq e^{\lambda t}, \qquad t >0;
$$
and let $\|\cdot\|_2$ be an equivalent Hilbertian norm on $H$ satisfying $\|T(\tau)\|_2 \leq 1$.

Set, for $h \in H$,
\begin{align}
	\|h\|_{\mathscr H}^2 := \inf \Bigg\{&
	\int_0^2 e^{-2\lambda s} 
	\left\| \sum_{j=1}^n T(t_j) h_j \, \chi_{[t_j,\,t_j+1]}(s) \right\|_1^2 \, ds  \notag \\
	&\quad + \sum_{j=1}^n a(t_j) \int_0^\tau \|T(s)h_j\|_2^2 \, ds \label{EqNormEq} \\
	&\quad + 2 \sum_{\substack{i,j=1 \\ i>j}}^n 
	b(t_i,t_j) \int_0^\tau 
	\Re \langle T(t_i-t_j+s)h_i, \, T(s)h_j \rangle_2 \, ds
	\Bigg\}, \nonumber
\end{align}
where the infimum is taken over all choices of 
$n \in \mathbb N$, 
$0 \leq t_1 < \cdots < t_n$, 
and $h_1,\dots,h_n \in H$ such that $h = \sum_{j=1}^n T(t_j)h_j$.
The coefficients $a(t)$ and $b(s,t)$ are given by
\[
a(t) := 1 + \mathrm{meas}\bigl([t,t+1]\cap[0,2]\bigr) =
\begin{cases}
	2, & t \leq 1, \\
	3-t, & 1 \leq t \leq 2, \\
	1, & t \geq 2,
\end{cases}
\]
and, for $s \geq t \geq 0$,
\[
b(s,t) := 1 + \mathrm{meas}\bigl([s,s+1]\cap[t,t+1]\cap[0,2]\bigr) =
\begin{cases}
	2+t-s, & s \leq 2,\ t \leq 1, \\
	3-s, & s \leq 2,\ t \geq 1, \\
	1, & s \geq 2 \text{ or } s-t \geq 1.
\end{cases}
\]

Then $\|\cdot\|_{\mathscr H}$ is an equivalent Hilbertian norm on $H$ such that $\|T(t)\|_{\mathscr H}\le 1$ for all $t\geq0$. \medskip

 The norm $\|\cdot\|_{\mathscr H}$ is an equivalent Hilbertian norm on $H$ such that
 $\|T(t)\|_{\mathscr H} \leq 1$ for all $t \geq 0$. 
 Its construction follows by tracing the renormings used in the proofs of
 Theorems~\ref{awayzero}, \ref{nearzero}, and~\ref{SimConstTh}, with the $C_0$-semigroup
 $\mathcal W$ there replaced by the variant $\widetilde{\mathcal W}$ defined in \eqref{tildeW}. In particular, the second and the third lines in~\eqref{EqNormEq} arise from the
 explicit form of unitary dilations for contraction semigroups (see, e.g.,
 \cite{davies1980one}). 
 The coefficients $a(t)$ and $b(s,t)$ account for overlaps of shifted indicator functions within $[0,2]$, reflecting the unitary dilation structure.
 The choice of $\widetilde{\mathcal W}$ provides a more transparent expression for $\|\cdot\|_{\mathscr H}$, while in the proof of Theorem~\ref{SimConstTh} the semigroup $\mathcal W$ was preferred, since it leads to simpler estimates for the similarity constants and for certain technical lemmas. 
 Consequently, the norm $\|\cdot\|_{\mathscr H}$ defined here differs from the one implicitly constructed in Theorem~\ref{SimConstTh}, and it may not satisfy the inequality in~\eqref{boundSimEq}.

  First, given a Hilbert space $\auxHilbert,$ we define, for each $t>0$,
  the operator $\widetilde V(t): L^2([0,2], \auxHilbert) \to L^2([0,1], \auxHilbert)$ by
 \begin{align*}
 	(\widetilde V(t) g)(x) =
 	\begin{cases}
 		0, & t < x, \\
 		g(2 + x - t), & t-1 < x \leq t, \\
 		g([2+x-t] \bmod 2) + g([3+x-t] \bmod 2), & x \leq t-1,
 	\end{cases}
 	\qquad x \in [0,1],
 \end{align*}
 for all $g \in L^2([0,2],\auxHilbert).$ For $t=0$ we set $\widetilde V(0) := 0$.
 
 Now define
 \begin{align}\label{tildeW}
 	\widetilde W(t) :=
 	\begin{pmatrix}
 		R_p(t) & \widetilde V(t) \\
 		0 & R(t)
 	\end{pmatrix}, 
 	\qquad t \geq 0,
 \end{align}
 where $R_p$ and $R$ are the semigroups introduced earlier. Then
 $\widetilde{\mathcal W} = (\widetilde W(t))_{t \geq 0}$ is a $C_0$-semigroup
 on
 \[
 \widetilde L^2_\oplus(\auxHilbert) := L^2([0,1], \auxHilbert) \oplus L^2([0,2], \auxHilbert).
 \]
 
 A direct computation shows that, for $h \in \auxHilbert$ and $t \geq 0$,
 \begin{align*}
 	\widetilde W(t)(0, h \chi_{[0,1]}) &=
 	\begin{cases}
 		(0,\, h \chi_{[t,1+t]}), & 0 \leq t \leq 1, \\
 		(h \chi_{[0,t-1]},\, h \chi_{[t,2]}), & 1 \leq t \leq 2, \\
 		(h \chi_{[0,1]}, 0), & t \geq 2.
 	\end{cases}
 \end{align*}
 
 Finally, introduce the operator
 \[
 (Mg)(x) := g(x) + g(x+1), \qquad x \in [0,1], \ g \in L^2([0,2], \auxHilbert),
 \]
 which coincides with the $L^2([0,2], \auxHilbert) \to L^2([0,1], \auxHilbert)$ component 
 of $\widetilde W(2)$. If we equip $\widetilde L^2_\oplus(\auxHilbert)$ with the equivalent norm
 \[
 \|(f,g)\|_{\oplus,{\rm eq}}^2 := \|f + Mg\|^2 + \|g\|^2, 
 \qquad (f,g) \in \widetilde L^2_\oplus(\auxHilbert),
 \]
 then, by an argument similar to Lemma~\ref{WcontractLemma}, the semigroup
 $\widetilde{\mathcal W}$ is contractive on 
 $\widetilde L^2_{\oplus,{\rm eq}}(\auxHilbert)$.

 The remaining arguments concerning $\|\cdot\|_{\mathscr H}$ are somewhat technical and omitted here, since the explicit formula is not essential  for the main development of the paper. A detailed account of \eqref{EqNormEq} will be given elsewhere.
 
\section{Joint similarity constants in terms of individual ones}\label{jointSect}

In this section we clarify the relevance of similarity constants $C(T(t)), t \ge 0,$ and $\mathcal C(\mathcal T)$
for a $C_0$-semigroup $\mathcal T=(T(t))_{t \ge 0}$ and emphasize their interplay
with our preceding results. Example \ref{CounterEx} below will further illustrate the importance of
similarity constants
and make the connections to similarity properties of concrete $\mathcal T$ explicit.

The next proposition links the individual similarity of
$T(t), t\ge 0,$ to contractions to their joint similarity to contractions for all
$t \ge 0$ via bounds on the individual similarity constants.
It also provides a characterization of joint similarity to contractions in terms of the resolvent of the generator, which may be of interest for applications.
The equivalence of statements (i) and (ii) is inspired by \cite[Proposition 2.5]{holbrook1977distortion}, which concerns the lower semicontinuity of similarity constants with respect to the strong operator topology.
It is worth noting that statement (ii) is local in nature,
depending only on the similarity constants near zero.

\begin{proposition}\label{CrSimSetProp}
	Let $H$ be a Hilbert space, and let $\mathcal T = (T(t))_{t\geq0}$ be a $C_0$-semigroup on $H,$ with generator $\generator$. Then the following statements are equivalent.
	\begin{itemize}
		\item [(i)]  $\mathcal T \in \mathcal{SC}(H)$.
		\item [(ii)] $\liminf_{t\to0} C(T(t)) < \infty$.
		\item [(iii)] $\liminf_{\lambda \to \infty} C(\lambda(\lambda - \generator)^{-1}) < \infty$.
	\end{itemize}
	If any of the above equivalent statement holds, then $(0,\infty) \subset \rho(\generator)$ and
	\begin{align*}
 		\mathcal C(\mathcal T) &= \lim_{t\to 0} C(T(t)) = \sup_{t>0} C(T(t))
		\\ &= \lim_{\lambda \to \infty} C \left(\lambda (\lambda - \generator)^{-1} \right) = \sup_{\lambda > 0 }C \left(\lambda (\lambda - \generator)^{-1}\right).
	\end{align*}
\end{proposition}
\begin{proof} $ $ \newline
	
	(i) $\implies$ (ii): This is trivial. \medskip
	
	(ii) $\implies$ (i): As $\mathcal C(\mathcal T) \geq C(T(t))$ for every $t>0$, one has
	\begin{align}\label{TrivialIneq1}
		 \mathcal C(\mathcal T) \geq \sup_{t>0} C(T(t)) \geq \liminf_{t\to 0} C(T(t)).
	\end{align}
	Let $\mathcal C:=\liminf_{t\to0} C(T(t))$, assume that $\mathcal C$ is finite,
	and choose a decreasing sequence $(t_k)_{k\in \NN}$ in $(0,\infty)$ such that $\lim_{k\to \infty} t_k = 0$ and
	$$\lim_{k \to \infty} C(T(t_k)) = \mathcal C.
		$$
	For each $k \in \NN$, let $\|\cdot\|_k$ be a Hilbertian norm on $H$ satisfying $\|T(t_k)\|_{k} \leq 1$ and
	\begin{align}\label{UniformBoundEq}
		\|h\| \leq \|h\|_k \leq C(T(t_k)) \|h\|, \qquad h\in H.
	\end{align}
	We define an inner product $\langle \cdot, \cdot\rangle_{\mathscr H}$ on $H$ by
	$$ \langle h, f\rangle_{\mathscr H} :=
		{\operatorname{LIM}} \left[\langle h, f\rangle_k\right]\qquad h, f \in H,
	$$
	where $\langle \cdot, \cdot\rangle_k $ is the inner product on $H$ induced by $\|\cdot\|_k$ and 
$\operatorname{LIM}$ is a fixed Banach limit on $\ell^\infty(\mathbb N).$ It follows from \eqref{UniformBoundEq} that 
$\|\cdot\|_{\mathscr H}$ is an equivalent Hilbertian norm on $H,$
	and
	\begin{equation}\label{equival}
	\|h\| \le \|h\|_{\mathscr H}\le \mathcal C \|h\|, \qquad h \in H.
	\end{equation}
	Let $\mathscr H$ denote the Hilbert space $(H,\|\cdot\|_{\mathscr H})$. Fix now $t>0$ and $h \in H$. For each $k\in \NN$, let $m_k \in \NN$ be such that $| t - m_k t_k|$ realizes the distance between $t$ and the set $t_k \NN$, that is,
	$$|t - m_k t_k| = \min_{m \in \NN} |t - m t_k|.
	$$
	Note that $\lim_{k\to \infty} |t-m_k t_k| = 0$ since $\lim_{k\to \infty} t_k = 0$. Consequently,
	$$\|T(t)h - T(m_k t_k)h\|_k \leq C(T(t_k)) \|T(t) h - T(m_k t_k)h\| \to 0, \quad \text{as}\quad k \to \infty.
	$$
	Also, $\|T(m_k t_k)\|_k \leq \|T(t_k)\|_k^{m_k} \leq 1$, $k \in \NN$. Therefore,
	\begin{align*}
		\|T(t)h\|_{\mathscr H}^2 &=
				{\operatorname{LIM}} \, \left[ \|T(t)h\|_k^2 \right]\leq
				{\operatorname{LIM}} \,  \left [\left(\|T(m_k t_k)h\|_k + \|T(t)h - T(m_k t_k)h\|_k\right)^2\right]
		\\ &=
				{\operatorname{LIM}} \, \left[ \|T(m_k t_k) h\|_k^2 \right] \leq
				{\operatorname{LIM}} \, [\|h\|_k^2] = \|h\|_{\mathscr H}^2.
	\end{align*}
		Thus $\|T(t)\|_{\mathscr H} \leq 1$ for all $t\geq0$, and $\mathcal T$ is a contraction $C_0$-semigroup
		on $\mathscr H$, so that
		$\mathcal C(\mathcal T) \leq \mathcal C$. It then follows from \eqref{TrivialIneq1} that
	$$\mathcal C(\mathcal T) = \sup_{t>0} C(T(t)) = \liminf_{t\to0} C(T(t)),
	$$
		hence $\lim_{t\to0} C(T(t))$ exists and equals to $\mathcal C(\mathcal T)$. \medskip
	
	(i) $\implies$ (iii): This implication, as well as $(0,\infty) \subset \rho(\generator)$ and the inequalities
	\begin{align}\label{TrivialIneq2}
		\mathcal C(\mathcal T) \geq \sup_{\lambda>0} \, C \left(\lambda (\lambda - \generator)^{-1}\right)  \geq \liminf_{\lambda \to \infty} \, C \left(\lambda (\lambda - \generator)^{-1}\right),
	\end{align}
	follow immediately from the fact that  $\mathcal T$ is a contraction semigroup if and only if
		$\|(\lambda-\generator)^{-1}\| \leq \frac{1}{\lambda}$ for all $\lambda>0$. \medskip
	
	(iii) $\implies$ (i): Assume $\mathcal C:= \liminf_{\lambda \to \infty} C \left(\lambda (\lambda-\generator)^{-1}\right) < \infty$, and fix a sequence $(\lambda_k)_{k\in \NN} \subset (0,\infty) \cap \rho(\generator)$ with $\lim_{k \to \infty} \lambda_k = \infty$ such that  $\lambda_k (\lambda_k -\generator)^{-1}$ is a contraction with respect to an equivalent Hilbertian norm $\|\cdot\|_k$ on $H$ and  $C_k := C \left( \lambda_k (\lambda_k - \generator)^{-1}\right)$ satisfy $\lim_{k\to \infty} C_k = \mathcal C$ and $\|h\|\leq \|h\|_K \leq C_k \|h\|$ for all $h\in H$. Arguing as above, set
	$$\langle h, f\rangle_{\mathscr H} :=
		{\operatorname{LIM}} \, [\langle h, f\rangle_k], \qquad h \in H,
	$$
	where again $\operatorname{LIM}$ is a Banach limit on $\ell^\infty(\mathbb N)$. In view of the choice of $(\lambda_k)_{k\in \NN}$ we conclude that  the corresponding norm $\|\cdot\|_{\mathscr H}$ is equivalent to the original norm on $H$
	and satisfies \eqref{equival}. In particular, $\mathscr H:=(H, \|\cdot\|_{\mathscr H})$ is a Hilbert space.
		
	Let $t>0$ be fixed. For each $k \in \NN$, choose  $t_k \in \frac{1}{\lambda_k}\NN$ to realize the distance between $\frac{1}{\lambda_k}\NN$ and $t$, i.e.,
	$$|t - t_k| = \min_{n\in \NN} \left|t - \frac{n}{\lambda_k}\right|,
	$$
	and let $n_k \in \NN$ be such that $t_k = \frac{n_k}{\lambda_k}$. Note that $\lim_{k\to\infty} t_k = t$ as $\lim_{k\to \infty} \lambda_k = \infty$. By the Post-Widder inversion formula \eqref{post}, for all $h \in H,$
		\begin{align*}
		\left\| T(t)h - \left(I - {t_k \generator} /{n_k}  \right)^{-n_k} h\right\|_k &\leq \left\| T(t) h - T(t_k) h \right\|_k + \left\| T(t_k) h - \left( I -{t_k}\generator/{n_k}  \right)^{-n_k} h \right\|_k
		\\ &\leq C_k \left( \left\| T(t) h - T(t_k) h \right\| + \left\| T(t_k) h - \left( I -{t_k}\generator /{n_k}  \right)^{-n_k} h \right\| \right)
		\\ &\rightarrow 0, \qquad \text{as}\,\, k \to \infty.
	\end{align*}
	where we used that the convergence in \eqref{post} is uniform on compacts from $[0,\infty).$
	 Therefore
\begin{align*}
	\|T(t)h\|_\mathscr H^2 &=
	{\operatorname{LIM}} \, \left[ \|T(t)h\|_k^2 \right]
	\\ &\leq
	{\operatorname{LIM}} \, \left[ \left( \left\| T(t)h - \left(I - t_k \generator/{n_k} \right)^{-n_k} h \right\|_k
	+ \left\| \left(I - t_k \generator/{n_k} \right)^{-n_k} h \right\|_k \right)^2 \right]
	\\ &=
	{\operatorname{LIM}} \, \left[ \left\| \left( \lambda_k \left(\lambda_k - \generator\right)^{-1} \right)^{n_k} h \right\|_k^2 \right]
	\\ &\leq
	{\operatorname{LIM}} \,\left[  \| h \|_k^2 \right] = \| h \|_\mathscr H^2, \qquad h \in H,
\end{align*}
		that is, $\|T(t)\|_\mathscr H \leq 1$. As $t>0$ was arbitrary, we conclude
	that $\mathcal T$ is a contraction $C_0$-semigroup on $\mathscr H$, so $\mathcal C(\mathcal T) \leq \mathcal C$. Hence $(0,\infty) \subset \rho(\generator)$, and  from \eqref{TrivialIneq2} it follows that $\lim_{\lambda \to \infty} C(\lambda(\lambda-\generator)^{-1})$ exists and
	$$\mathcal C(\mathcal T) = \sup_{\lambda>0} C\left(\lambda(\lambda-\generator)^{-1}\right) = \lim_{\lambda \to \infty} C \left(\lambda(\lambda-\generator)^{-1}\right).
	$$
\end{proof}

Now as a consequence of Theorem \ref{SimConstTh} and Proposition \ref{CrSimSetProp}, we present the following ``trichotomy'' corollary that clarifies the interplay between $\mathcal{C}(\mathcal{T})$ and ${C(T(t)),\, t \geq 0}$, highlighting their significance in characterizing semigroups in $\mathcal{SC}(H)$ and clarifying fine structural properties of
such semigroups.
  Note that it is immediate that $\mathcal{T} \in \mathcal{SC}(H)$ implies that $T(t)$ is similar to a contraction for every $t > 0$ and that $\sup_{t > 0} C(T(t)) < \infty$. Conversely, Proposition \ref{CrSimSetProp} ensures that this condition is also sufficient.

\begin{corollary}\label{trichotomyCor}
	Let $H$ be a Hilbert space, and let $\mathcal T = (T(t))_{t\geq0}$ be a $C_0$-semigroup on $H$. Then only one of the following alternatives holds.
	\begin{itemize}
		\item [(i)] $\mathcal T \in \mathcal{SC}(H)$. Furthermore,
		$$
		\mathcal C(\mathcal T) = \sup_{t>0} C(T(t)) = \lim_{t\to 0} C(T(t)).
		$$
		\item [(ii)] For each $\tau>0$, the operators $(T(t))_{t\geq \tau}$ are jointly similar to contractions, i.e., $\sup_{t>\tau} C(T(t))<\infty$. At the same time, $\mathcal T$ does not belong to $\mathcal{SQC}(H)$, and thus
		$$
		\mathcal C(\mathcal T) = \lim_{t\to 0} C(T(t)) = \infty.
		$$  
		\item [(iii)] For each $t>0$, $T(t)$ is not similar to a contraction. Thus,
		$$
		\mathcal C(\mathcal T) = C(T(t)) = \infty, \qquad t>0.
		$$
	\end{itemize}
\end{corollary}

\begin{proof}
	It is clear that the properties (i)--(iii) are mutually exclusive. Letting $\mathcal T = (T(t))_{t\geq0}$ be a $C_0$-semigroup on $H,$ we show that $\mathcal T$ satisfies at least one of them.
	
	By Corollary \ref{SimSetCor}, either there is no $t>0$ such that $T(t)$ is similar to a contraction, or the operators in $(T(t))_{t\geq\tau}$ are jointly similar to contractions for every $\tau > 0$. In the first case, (iii) holds. In the second case, $C(T(t))$ is finite for every $t \ge 0$, and there are two further alternatives. If $\liminf_{t\to 0} C(T(t)) < \infty$, Proposition \ref{CrSimSetProp} implies that (i) holds. Otherwise, $\lim_{t\to 0} C(T(t)) = \infty$, so $\mathcal T$ does not belong to $\mathcal{SC}(H)$. It then follows from Theorem \ref{quasiContrTh_int} that $\mathcal T$ does not belong to $\mathcal{SQC}(H)$ either, so (ii) holds.
\end{proof}

\section{Failure of similarity to contraction semigroups}\label{ExamplesSect}

The theory developed so far showed that
the study of similarity to contraction $C_0$-semigroups
leads to new phenomena absent in the discrete setting.
It is thus desirable to illustrate and clarify them 
with instructive examples. The present section serves this purpose
by improving and strengthening several well-known examples from the literature.
Consequently, we show that assumptions (i) and (ii) in Theorem \ref{SimConstTh} are independent of each other,
and moreover, they may fail simultaneously.
The latter allows us to produce semigroups satisfying very strong requirements yet
far from $\mathcal{SC}(H).$
The present section serves this purpose
by improving and strengthening several well-known examples from the literature.

\subsection{Packel type semigroups and similarity to contractions}\label{PackelSect}

The first example of a bounded $C_0$-semigroup on a Hilbert space which is not similar to a semigroup of contractions was found by Packel in \cite{packel1969semigroup}.
Here we propose a more general version of this example allowing us
to simplify and generalize other examples in the literature,
and provide an alternative construction for some of the examples presented 
in \cite{oliva2025tensor}. 

Let $J$ be either $\ZZ$, $\ZZ_+ = \{n \in \ZZ \, : \, n \geq 0\}$ or $\ZZ_- = \{n \in \ZZ \, : \, n \leq 0\}$, and let $\mathbf a = (a_n)_{n\in J} \subset (0,\infty)$ be an increasing sequence such that
\begin{align}\label{aEq}
\begin{cases}
	\mbox{if } J = \ZZ \mbox{ or } J = \ZZ_+, \qquad &a_{n+1} \geq 2a_n, \qquad n\in J,
	\\\mbox{if } J = \ZZ_-, \qquad &a_n \geq 2a_{n-1}, \qquad n\in J.
\end{cases}
\end{align}
Fix $t>0$. If there exists $n \in J$ such that $a_{n} < t \leq a_{n+1}$, then define $n_0(t) = n$. 
If no such $n$ exists, then for $J=\mathbb Z_-$ and $t>a_0$ set $n_0(t)=0$, while for $J=\mathbb Z_+$ and $t\le a_0$ set $n_0(t)=-1$ with $a_{-1}=0$.
{In the latter case, we allow the temporary index $-1\notin J$ solely to define $n_0(t)$ and the interval family below, all subsequent sums use indices $n\ge n_0(t)$.}

For each $n \in J$ such that $n > n_0(t)$, set $I_n(t) := (a_n-t, a_n]$ and let $I_{n_0}(t) = [0, 2 a_{n_0} -t]$ for $n_0= n_0(t)$. Thus, the operator $V_{\mathbf a}(t) \in \linearOp(L^2(\RR_+))$ given by
\begin{align*}
(V_{\mathbf a}(t)f)(x): = \begin{cases}
	f  \left(2a_n - x - t \right), &\quad x \in I_n(t), \, n\geq n_0(t), 
	\\ 0, &\quad \mbox{otherwise},
\end{cases}
\qquad \qquad x \in (0,\infty), \,f \in L^2(\RR_+),
\end{align*}
is well-defined,
and we let $V_{\mathbf a}(0)$
stand for the zero operator. 
Note that for fixed $t$, the intervals $\{I_n(t)\}_{n \geq n_0(t)}$ are disjoint, which makes the operator behave ``isometrically'' on disjoint supports.
Using this, it is direct to see that
$\|V_{\mathbf a}(t)\| \leq 1.$ 
{Moreover, 
\[
\|V_{\mathbf a}(t)f\|_2^2=\sum_{n\ge n_0(t)}\!\int_{a_n-t}^{a_n}\!|f(u)|^2\,du
=\int_0^\infty |f(u)|^2 g_t(u)\,du
\]
 with $0\le g_t\le1$ and $g_t(u)\to0$ a.e. as $t\to 0$, 
so $V_{\mathbf a}(t)\to0$ as $t\to 0$ strongly by the dominated convergence theorem.}

Let now $\mathcal R = (R(t))_{t\geq0}$ and $\mathcal L = (L(t))_{t\geq0}$ be the right and the left shift 
semigroups on $L^2(\RR_+),$ respectively.
{Since}
\begin{equation}\label{semiprop}
{V_{\mathbf a} (s+t) = L(s) V_{\mathbf a}(t) + V_{\mathbf a}(s) R(t), \qquad s,t\geq0,}
\end{equation}
(cf. the proof of Lemma \ref{WcontractLemma})  setting 
\begin{equation}\label{tmatr}
T_{\mathbf a}(t): = \begin{pmatrix}
	L(t) & V_{\mathbf a}(t)\\[2pt]
	0 & R(t)
\end{pmatrix},
\qquad t\geq0,
\end{equation}
{yields a $C_0$-semigroup.
Following closely the proof of \cite[property (I), pp. 242-243]{packel1969semigroup},
we infer that $\mathcal T_{\mathbf a}$} 
is bounded  on $L^2(\RR_+) \oplus L^2(\RR_+)$. 
Figure \ref{packelFig} illustrates the nature of  $\mathcal T_{\mathbf a}.$
\begin{center}
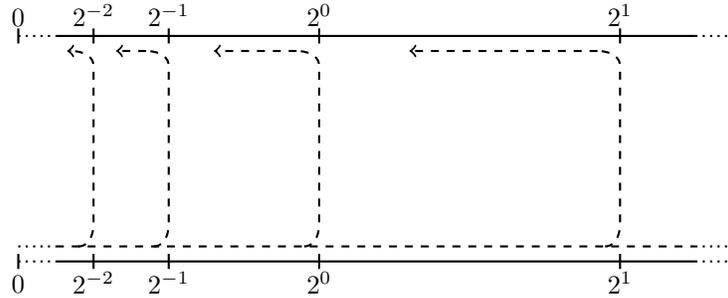
\begin{figure}
	\begin{tikzpicture}
		\draw[thick] (0.5,1.5) -- (9,1.5); 
		\draw[thick] (0.5,-1.5) -- (9,-1.5); 
		
		\draw[dotted,thick] (0,-1.5) -- (0.5, -1.5);
		\draw[dotted,thick] (9,-1.5) -- (9.5, -1.5);
		\draw[thick] (0,-1.4) -- (0,-1.6);
		\node[below] at (0,-1.55) {\small $0$};
		
		\draw[dotted,thick] (0, 1.5) -- (0.5, 1.5);
		\draw[dotted,thick] (9, 1.5) -- (9.5, 1.5);
		\draw[thick] (0,1.4) -- (0,1.6);
		\node[above] at (0,1.5) {\small $0$};

		\foreach \x in {-2,-1,0,1} {
			\pgfmathsetmacro\pos{4*2^\x} 
			
			\draw[thick] (\pos,1.6) -- (\pos,1.4);
			\node[above] at (\pos,1.5) {\small $2^{\x}$};
			
			\draw[thick] (\pos,-1.4) -- (\pos,-1.6);
			\node[below] at (\pos,-1.5) {\small $2^{\x}$};
		}
		
		\draw[dashed,thick] (0.5, -1.3) -- (9, -1.3);
		\draw[dotted,thick] (0,-1.3) -- (0.5, -1.3);
		\draw[dotted,thick] (9,-1.3) -- (9.5, -1.3);

		\foreach \x in {-2,-1,0,1} { 
			\pgfmathsetmacro\posA{2*2^\x} 
			\pgfmathsetmacro\posB{2*2^(\x+1)} 
			
			\draw[dashed,thick]  (\posB-0.2,-1.3) to[out=0, in=-90] (\posB,-1.3+0.3) -- (\posB,1.1) to[out=90,in=0] (\posB-0.3,1.3);
			
			\draw[dashed,thick,->] (\posB - 0.3, 1.3) to[out=180, in=0] (\posA*1.3,1.3);
		}
	\end{tikzpicture}
	\caption{Visualization of $\mathcal T_{\mathbf a}$ for $\mathbf a=(2^n)_{n\in\mathbb Z}$, illustrating how function supports are shifted under $\mathcal T_{\mathbf a}$.}
	\label{packelFig}
\end{figure}
\end{center}

Packel proved in \cite{packel1969semigroup} that if  ${\mathbf a} = (4^{n})_{n\in \ZZ}$, then 
$\mathcal T_{\mathbf a} \notin\mathcal{SC}(L^2(\RR_+) \oplus L^2(\RR_+))$. Later on, Chernoff \cite{chernoff1976two} used an infinite direct sum of scaled versions of Packel's semigroup (which applies in fact to any  bounded $C_0$-semigroup not similar to a contraction one) to prove, in particular, the existence of a semigroup that is not similar to a quasi-contraction semigroup. The examples by Chernoff and Packel became
basic in the literature. We show below that Packel's semigroup does not belong to $\mathcal{SQC}(L^2(\RR_+) \oplus L^2(\RR_+))$, and thus, Chernoff's additional argument is not needed for this purpose.
Moreover, we construct a family of analogous semigroups which, being not similar to a contraction semigroup may or may not  belong to $\mathcal{SQC}(L^2(\RR_+) \oplus L^2(\RR_+))$,
depending on the choice of $\mathbf a.$

\begin{theorem}\label{PackelSemTh}
Let $J$ be either $\ZZ, \, \ZZ_-$ or $\ZZ_+$, and let $\mathbf a = (a_n)_{n\in J} \subset (0,\infty)$ satisfy \eqref{aEq}.
\begin{enumerate}
	\item [(i)] For $t>0$, $T_{\mathbf a}(t)$ is similar to a contraction if and only if $J = \ZZ_-$.
	\item [(ii)] $\mathcal T_{\mathbf a} \in \mathcal{SQC}(L^2(\RR_+) \oplus L^2(\RR_+))$ if and only if $J = \ZZ_+$. 
\end{enumerate}
\end{theorem}

To prove (ii), we first show that if $J = \ZZ_-,$ then the compression of $\mathcal T_{\mathbf a}$ to a suitable subspace is a \emph{nilpotent} $C_0$-semigroup that is not similar to a quasi-contraction $C_0$-semigroup. This result is interesting in its own right as a strong counterexample to an analogue of Rota's theorem for $C_0$-semigroups. Moreover, in contrast to the
examples presented in \cite{chernoff1976two, lemerdy2000bounded, oliva2025tensor}, the semigroup is 
given by an explicit and simple formula  amenable for other computational purposes.
A different example of nilpotent semigroup not in $\mathcal{SC}(H),$
and, in addition, immediately compact,  
was constructed in \cite{oliva2025tensor}.

\begin{proposition}\label{nilpotentExplProp}
Let $b>0$, and let $\mathbf a = (a_n)_{n\in \ZZ_-} \subset (0,b]$ satisfy \eqref{aEq}. Then, $H = L^2(0,b) \oplus L^2(0,b)$ is the orthogonal difference\footnote{identified with a subspace of $L^2(\RR_+) \oplus L^2(\RR_+)$.} 
of two $\mathcal T_{\mathbf a}$-invariant subspaces, and the compression $\mathcal N_{\mathbf a} = (N_{\mathbf a}(t))_{t\geq0}$ of $\mathcal T_{\mathbf a}$ to $H$ is a nilpotent $C_0$-semigroup that does not belong to $\mathcal{SQC}(H)$.
\end{proposition}
\begin{proof}
Making the natural identification $L^2(\Delta) = \{f \in L^2(\RR_+) \, : \, \operatorname{supp}(f) \subseteq \Delta\}$ with $\Delta \subseteq (0,\infty)$, note that  $\{0\} \oplus L^2(b,\infty) \subset L^2(0,b) \oplus L^2(\RR_+)$ are $\mathcal T_{\mathbf a}$-invariant subspaces and
$$ \left(L^2(0,b) \oplus L^2(\RR_+)\right) \ominus \left(\{0\} \oplus L^2(b,\infty) \right) = L^2(0,b) \oplus L^2(0,b) = H.
$$
Hence, by Sarason's characterization of dilations (see Section \ref{PreliminariesSect}), the compression $\mathcal N_{\mathbf a}$ of $\mathcal T_{\mathbf a}$ to $H$ is a well-defined $C_0$-semigroup on $H$. Clearly, $\mathcal N_{\mathbf a}$ is a nilpotent $C_0$-semigroup since $N_{\mathbf a}(2b) = 0.$
In view of Theorem \ref{quasiContrTh_int}, it suffices to prove that 
$\mathcal N_{\mathbf a} \not \in \mathcal{SC}(H)$.

We proceed by contradiction. Assume to the contrary that $\mathcal N_{\mathbf a} \in \mathcal{SC}(H)$ 
and let $\|\cdot\|_{\mathscr H}$ be a Hilbertian norm on $H$ such that $\mathcal N_{\mathbf a}$ is a contraction $C_0$-semigroup on $\mathscr H:=(H, \|\cdot\|_{\mathscr H})$ and
$$\|\cdot\| \leq \|\cdot\|_{\mathscr H} \leq \mathcal C (\mathcal N_{\mathbf a}) \|\cdot\|.
$$
By Sz.-Nagy's dilation theorem, there exists a unitary dilation $\mathcal U_{\mathbf a} = (U_{\mathbf a}(t))_{t\in \RR}$ of $\mathcal N_{\mathbf a}$ on a Hilbert space $\mathscr K$ containing $\mathscr H$ as a subspace. Set
$$\mathscr F := \overline{\operatorname{span}} \{U_{\mathbf a}(t)h \, : \, h\in H, \, t\geq 0\}.
$$
Clearly, $\mathscr F$ is a $\mathcal U_{\mathbf a}$-invariant subspace containing $\mathscr H$ as a subspace, the restriction $\mathcal S_{\mathbf a} =(S_{\mathbf a}(t))_{t\geq0}$ of $\mathcal U_{\mathbf a}$ to $\mathscr F$ is an isometric $C_0$-semigroup, and  $\mathcal T_{\mathbf a}$ is the compression of 
$\mathcal S_{\mathbf a}$ to $\mathscr H$. Let 
\begin{align*}
	S_{\mathbf a}(t) = \begin{pmatrix}
		S_{\mathbf a}^1(t) & S_{\mathbf a}^2(t)
		\\0	&	T_{\mathbf a}(t)
	\end{pmatrix}, \qquad t \geq 0,
\end{align*}
be the matrix representation of $\mathcal S_{\mathbf a}$ with respect to the decomposition $\mathscr F = (\mathscr F \ominus \mathscr H) \oplus \mathscr H.$
Since $\mathcal S_{\mathbf a}$ is an isometric semigroup, we have
\begin{equation}\label{isomPack}
	\begin{aligned}
		\|h\|_{{\mathscr H}}^2 &= \| S_{\mathbf a} (t) h\|_{{\mathscr H}}^2 
		=  \|T_{\mathbf a}(t)h\|_{{\mathscr H}}^2 + \|S_{\mathbf a}^2(t) h\|_{{\mathscr H}}^2, \qquad t \geq 0, \, h \in \mathscr H.
	\end{aligned}
\end{equation}

Given $h\in \mathscr H$, $N \in \NN$ and $0=t_0<t_1<\ldots<t_N$, an iteration of \eqref{isomPack} yields
\begin{equation}\label{AuxIneqPack}
	\begin{aligned}
		\|h\|_{{\mathscr H}}^2 & = \|T_{\mathbf a}(t_1) h\|_{{\mathscr H}}^2 + \|S^2_{\mathbf a} (t_1-t_0) h\|_{{\mathscr H}}^2
		\\ &= \|T_{\mathbf a}(t_2) h\|_{{\mathscr H}}^2 + \|S^2_{\mathbf a} (t_2-t_1) T_{\mathbf a} (t_1)h\|_{{\mathscr H}}^2 + \|S^2_{\mathbf a} (t_1-t_0) h\|_{{\mathscr H}}^2
		\\ &= \ldots = \|T_{\mathbf a}(t_N)h\|_{{\mathscr H}}^2 + \sum_{j=0}^{N-1} \|S_{\mathbf a}^2(t_{j+1} - t_j) T_{\mathbf a}(t_j) h\|_{{\mathscr H}}^2.
	\end{aligned}
\end{equation}
Fix $n\in \ZZ_-$ for the rest in the proof. For all $n\leq k \leq 0$, let $h_n = \left(0,\chi_{(0,a_n)}\right) \in \mathscr H$, and set $t_n(k) := 2 a_k - a_n$ so that $t_n(k) \in (0,2b)$ and $V_{\mathbf a}(t_n(k)) \chi_{(0,a_n)} = \chi_{(0,a_n)}$. Set also $t_n(n-1) = 0$. Note that the sequence $(t_n(k))_{k = n-1,\ldots,0}$ is increasing in $k$, and
\begin{align*}
	T_{\mathbf a}(t_n(k))h_n &= \left(\chi_{(0, a_n)}, R(t_n(k)) \chi_{(0, a_n)} \right), \qquad n \leq k \leq 0.
\end{align*}
{As $a_n = t_n(n)$ and
$$t_n(k+1) - t_n(k)  = 2(a_{k+1} - a_k) \geq 2 a_k > a_n$$ 
for all $n\leq k \leq 0$, 
the functions $R(t_n(k))\chi_{(0,a_n)}=\chi_{(t_n(k),\,t_n(k)+a_n)}$ have pairwise disjoint supports for $n\le k\le0.$ Hence}
\eqref{AuxIneqPack} applied to the time-sequence 
$$0,\, a_n,\, t_n(n)+a_n,\, t_n(n+1),\, t_n(n+1)+a_n,\,\ldots,\, t_n(0), \, t_n(0) + a_n,
$$
yields
\begin{align*}
	\|h_n\|_{{\mathscr H}}^2 & = \|T_{\mathbf a}(t_n(0)+a_n)h_n\|_{{\mathscr H}}^2 + \|S_{\mathbf a}^2 (a_n) T_{\mathbf a}(t_n(0)) h_n\|_{{\mathscr H}}^2 
	\\ &\qquad +  \sum_{k=n}^{0} \bigg( \left\|S_{\mathbf a}^2(t_n(k)-t_n(k-1) -a_n) T_{\mathbf a}(t_n(k-1)+a_n) h_n \right\|_{{\mathscr H}}^2 
	\\ & \qquad \qquad + \left\|S_{\mathbf a}^2(a_n) T_{\mathbf a}(t_n(k-1)) h_n \right\|_{{\mathscr H}}^2 \bigg)
	\\ &\geq \sum_{k=n}^{0} \left\|S_{\mathbf a}^2(a_n) T_{\mathbf a}(t_n(k)) h_n \right\|_{{\mathscr H}}^2 = \sum_{k=n}^{0} \left\|S_{\mathbf a}^2(a_n) \left(\chi_{(0, a_n)},R(t_n(k)) \chi_{(0, a_n)} \right) \right\|_{{\mathscr H}}^2.
\end{align*}
Hence, {by the Cauchy-Schwarz inequality}, we obtain
\begin{align}\label{CSchwPackEq}
	\sum_{k=n}^{0} \left\|S_{\mathbf a}^2(a_n) \left(\chi_{(0, a_n)}, R(t_n(k)) \chi_{(0, a_n)} \right) \right\|_{{\mathscr H}} \leq (|n|+1)^{1/2} \|h_n\|_{{\mathscr H}}. 
\end{align}
Moreover, since $L(a_n)\chi_{(0,a_n)} = 0$,  we have $T_{\mathbf a}(a_n) \left(\chi_{(0,a_n)},0 \right) = 0,$ and  \eqref{isomPack} implies 
$$\left\| \left( \chi_{(0,a_n)}, 0\right)  \right\|_{{\mathscr H}} =  \left\| S_{\mathbf a}^2(a_n)  \left(\chi_{(0,a_n)}, 0\right) \right\|_{{\mathscr H}}.
$$
Hence,
\begin{equation}\label{ineqq}
	\begin{aligned}
	& \left\| S_{\mathbf a}^2 (a_n) \left(0,\sum_{k=n}^{0} R(t_n(k)) \chi_{(0, a_n)} \right)\right\|_{{\mathscr H}}
	\\ \geq &  \left\|\sum_{k=n}^{0} S_{\mathbf a}^2 (a_n) \left(\chi_{(0,a_n)}, 0 \right)\right\|_{{\mathscr H}} - \sum_{k=n}^{0} \left\| S_{\mathbf a}^2 (a_n) \left(\chi_{(0,a_n)}, R(t_n(k)) \chi_{(0, a_n)} \right)\right\|_{{\mathscr H}}
	\\ \geq &  (|n|+1) \left\| \left(\chi_{(0,a_n)}, 0 \right)\right\|_{{\mathscr H}} - (|n|+1)^{1/2} \|h_n\|_{{\mathscr H}}
	\\ \geq & \left((|n|+1) - \mathcal C(\mathcal N_{\mathbf a}) (|n|+1)^{1/2}\right) (a_n)^{1/2}, 
	\end{aligned}
\end{equation}
where we used that $$\left\|\left(\chi_{(0,a_n)}, 0\right)\right\| = \left\|\left(0, \chi_{(0,a_n)}\right)\right\| = \|\chi_{(0,a_n)}\| = (a_n)^{1/2}.$$ 

On the other hand, note that the functions  
$$R(t_n(k))\chi_{(0,a_n)} = \chi_{(t_n(k), a_n + t_n(k))}, \qquad n\leq k \leq 0,
$$
have mutually disjoint supports, so
$$\left\|\sum_{k=n}^{0} R(t_n(k)) \chi_{(0, a_n)}\right\| \leq (|n|+1)^{1/2} (a_n)^{1/2}. 
$$
Taking into account that $\|S_{\mathbf a}^2(t)\| \leq \|S_{\mathbf a}(t)\|= 1$ for all $t\geq0,$ we infer that
\begin{align}\label{estimate_s}
	\left\| S_{\mathbf a}^2 (a_n) \left(0,\sum_{k=n}^{0} R(t_n(k)) \chi_{(0, a_n)} \right)\right\|_{{\mathscr H}} &\leq  \left\| \left(0,\sum_{k=n}^{0} R(t_n(k)) \chi_{(0, a_n)} \right)\right\|_{{\mathscr H}}
	\\ & \leq   \mathcal C(\mathcal N_{\mathbf a}) \, (|n|+1)^{1/2} (a_n)^{1/2}, 
\end{align}
Letting $n\to -\infty$, it follows that \eqref{estimate_s} contradicts \eqref{ineqq}, which finishes the proof.
\end{proof}

In addition to Lemma \ref{nilpotentExplProp}, we will need a well-known asymptotic property of 
semigroup orbits, going back to \cite{foguel1964counterexample} and already used in similar contexts.
Given a $C_0$-semigroup $\mathcal T = (T(t))_{t\geq0}$ on a Hilbert space $H$, set
\begin{align*}
W(\mathcal T) := \{h \in H \, : \, \lim_{t\to\infty} T(t)h = 0 \;\text{weakly}\},
\end{align*}
and, given a bounded operator $T \in \linearOp(H)$, define analogously
\begin{align*}
	W(T) := \{h \in H \, : \, \lim_{n\to\infty} T^n h = 0 \;\text{weakly}\}.
\end{align*}
Recall the following result proved in \cite[Lemma 1]{packel1969semigroup}, see also 
\cite[Theorem 2.5]{foguel1963power} and \cite[Theorems 4 and 5]{holbrook1971iterates}.
\begin{lemma}\label{weakLimitLemma}
Let $\mathcal T = (T(t))_{t\geq0}$ be a $C_0$-semigroup on a Hilbert space $H$. If $\mathcal T \in \mathcal{SC}(H)$, then 
\[
W(\mathcal T) \cap W(\mathcal T^\ast)^\perp = \{0\}.
\]
\end{lemma}

The proof of Lemma \ref{weakLimitLemma} in \cite{packel1969semigroup} relies on the fact that 
$\mathcal T\in \mathcal{SC}(H)$ is equivalent to similarity of $\mathcal T$ to a semigroup
admitting a unitary dilation, and the statement then essentially reduces to the unitary case. 
We will need a different approach, relying on the discrete setting, to obtain the following stronger version.
Its proof is comparatively direct and puts together several well-known facts from semigroup theory.
\begin{lemma}\label{strongPackLemma}
Let $\mathcal T = (T(t))_{t\geq0}$ be a $C_0$-semigroup on a Hilbert space $H$. 
If there exists $\tau>0$ such that $T(\tau)$ is similar to a contraction, then 
\[
W(\mathcal T) \cap W(\mathcal T^\ast)^\perp = \{0\}.
\]
\end{lemma}

\begin{proof}
Fix $\tau > 0$ such that $T(\tau)$ is similar to a contraction, and let $R\in \linearOp(H)$ be positive and invertible such that 
$S:= R T(\tau) R^{-1}$ is a contraction. It follows from \cite[Theorem 3.1]{foguel1963power} that $W(S) = W(S^\ast)$. 
Since $h \in W(T(\tau))$ if and only if 
\[
\lim_{n\to \infty} S^n R h = \lim_{n\to \infty} RT(\tau)^n h = 0 \quad \text{weakly},
\]
we infer that 
\begin{equation}\label{STWeak}
	 W (T(\tau)) = R^{-1} (W(S)).
\end{equation}
Similarly,
\begin{equation}\label{STWeakAst}
	W(T(\tau)^\ast) = R(W(S^\ast)) = R(W(S)),
\end{equation}
since $S^\ast = R^{-1} T(\tau)^\ast R$. 
In addition, following the arguments given in \cite[Lemma 3.3.1]{benchimol1977feedback} 
(or \cite[p.~788]{foguel1964counterexample}), we have
\[
h \in R^{-1}(W(S)^\perp) \;\iff\; Rh \in W(S)^\perp \;\iff\; h \in (R(W(S)))^\perp,
\]
and using \eqref{STWeakAst} we deduce that $R^{-1}(W(S)^\perp) = (W(T(\tau)^\ast))^\perp$. 
Combining this with \eqref{STWeak}, we conclude that 
\begin{equation}\label{WTdiscrete}
	W(T(\tau)) \cap (W(T(\tau)^\ast))^\perp = \{0\}.
\end{equation}

Next we prove that 
\begin{equation}\label{SemDiscreteWeak}
	W(T(\tau)) = W(\mathcal T).
\end{equation}
Clearly $W(\mathcal T) \subseteq W(T(\tau))$. For the opposite inclusion, fix $h \in W(T(\tau))$. 
Since $\mathcal T$ is bounded, for any relatively compact subset $M \subset H$, we have
\[
\lim_{n\to \infty} \langle T(n \tau) h , f \rangle = 0 \quad \text{uniformly for } f \in M.
\]
For every $f \in H$, the set $M_f:= \{T(t)^\ast f \, : \, t \in [0,\tau)\}$ is relatively compact, so if $(t_n)_{n\in \NN}\subset [0,\infty)$ with $t_n \to \infty$, we may write $t_n = m_n \tau + r_n$ with $m_n \to \infty$, $m_n  \in \NN$ and $r_n \in [0,\tau)$. Then
\[
\langle T(t_n)h, f \rangle = \langle T(m_n \tau ) h, T(r_n)^\ast f \rangle \to 0,
\]
since $T(r_n)^\ast f \in M_f$. As $(t_n)_{n \in \mathbb N}$ was arbitrary, this shows $\lim_{t\to\infty} T(t)h = 0$ weakly, i.e. $h \in W(\mathcal T)$. 
Thus \eqref{SemDiscreteWeak} holds. An analogous argument yields $W(\mathcal T^\ast) = W(T(\tau)^\ast)$. 
Together with \eqref{WTdiscrete}, this gives the claim.
\end{proof}
\begin{remark}
Any invertible similarity $R$ suffices in the first part, its positivity is convenient but not essential.
The property \label{SemDiscreteSWeak} could be well-known, but we were not able to find a
direct reference. For its global counterpart see e.g. \cite[Theorem 3.4]{eisner2010stability}.
\end{remark}

We are now ready to prove Theorem \ref{PackelSemTh}.

\begin{proof}[Proof of Theorem \ref{PackelSemTh}]

To prove (i), assume first $J= \ZZ$ or $J = \ZZ_+$. For $k \in \ZZ_+$ such that $a_k\geq 1$, set $t_k := 2 a_k - 1$ so that 
$$T_{\mathbf a}(t_k) (0,\chi_{(0,1)}) = (\chi_{(0,1)}, \chi_{(t_k, t_k+1)}).
$$
Since $\lim_{k\to \infty} t_k = \infty$, it follows that 
\[
\lim_{k\to \infty} T_{\mathbf a}(t_k) (0,\chi_{(0,1)}) = (\chi_{(0,1)},0) \quad \text{weakly},
\]
 and then $(\chi_{(0,1)},0) \in W(\mathcal T_{\mathbf a}^\ast)^\perp$. 
On the other hand, clearly  $(\chi_{(0,1)},0) \in W(\mathcal T_{\mathbf a})$ since $T_{\mathbf a}(1)(\chi_{(0,1)},0) = 0$.
Thus our claim is a consequence of Lemma \ref{strongPackLemma}.

If $J = \ZZ_-$, then $V_{\mathbf a}(t) = 0$ for all $t \geq 2 a_0$. Hence $T_{\mathbf a}(t)$ is a contraction for every $t\geq 2 a_0$, and our claim then follows from Corollary \ref{SimSetCor}. \medskip

Now we turn to the proof of (ii).
Assume first $J = \ZZ_-$. If  $\mathcal T_{\mathbf a} \in \mathcal{SQC}(L^2(\RR_+) \oplus L^2(\RR_+))$, then
Remark \ref{similar_dil} applies to $e_{-\lambda }\mathcal T_{\mathbf a}$ for some $\lambda \geq 0.$ Therefore, from  $e_{-\lambda } \mathcal T_{\mathbf a} \in \mathcal{SC}(L^2(\RR_+) \oplus L^2(\RR_+))$ (see also (iv) $\implies$ (i) in Theorem \ref{nearzero}) it follows that every compression of $\mathcal T_{\mathbf a}$ to semi-invariant
subspace  is similar to a quasi-contraction semigroup, which contradicts Proposition \ref{nilpotentExplProp}. \medskip

Consider now $J = \ZZ$, set $b = a_0>0$ and identify $H = L^2(0,b) \oplus L^2(0,b)$ with a subspace of $L^2(\RR_+) \oplus L^2(\RR_+)$. Set also $\mathbf c = (a_n)_{n\in \ZZ_-}$.  By Proposition \ref{nilpotentExplProp}, the compression $\mathcal N_{\mathbf c} = (N_{\mathbf c}(t))_{t\geq0}$ of $\mathcal T_{\mathbf c}$ 
to $H$ does not belong to $\mathcal{SQC}(H).$
It is straightforward to verify that
\begin{align*}
	N_{\mathbf c} (t) = \proj_H T_{\mathbf a}(t)\restriction_H, \qquad t \in [0,2b],
\end{align*}
where $\proj_H$ is the orthogonal projection from $L^2(\RR_+) \oplus L^2(\RR_+)$ onto $H$. Hence, from Theorem \ref{nearzero}, (iv) $\implies$ (i), it follows that $\mathcal T_{\mathbf a} \notin \mathcal{SQC}(L^2(\RR_+) \oplus L^2(\RR_+))$.
\medskip

Suppose now $J = \ZZ_+.$ We identify $L^2(\RR_+) \oplus L^2(\RR_+)$ with a closed subspace of $L^2(\RR) \oplus L^2(\RR).$ 
For each $t>0$ and each $n\in \ZZ_+$, set $I^n(t) := (a_n-t,a_n]$ and define $W_{\mathbf a}(t)\in \linearOp(L^2(\RR))$ by
\begin{align*}
	(W_{\mathbf a}(t)f)(x): = \sum_{n\in \ZZ^+} \chi_{I^n(t)}(x) f(2a_n-x-t),\qquad x \in \RR,\, f \in L^2(\RR),
\end{align*}
where the sum above is finite for each $x\in \RR$. The boundedness of $W_{\mathbf a}$ follows from the fact  that the intervals $\{I_n(t)\}_{n>n_0(t)}$ are mutually disjoint. 
It is also straightforward that $\lim_{t\to 0} W_{\mathbf a}(t) = 0$ in the strong operator topology.

Setting $W_{\mathbf a}(0)$ to be the zero operator and arguing as in the proof of \cite[property (I), pp. 242-243]{packel1969semigroup},
it is easy to  see that
\begin{equation}\label{pack1}
	W_{\mathbf a}(s+t) = L_\RR(s) W_{\mathbf a}(t) + W_{\mathbf a}(s) R_\RR(t), \qquad s,t\geq0,
\end{equation}
where $(L_\RR(t))_{t\in \RR}$
and $(R_\RR(t))_{t\in \RR}$ are the left shift and right shift semigroups on $L^2(\RR)$ (so $L_\RR(t) = R_\RR(-t)$), respectively.

Thus the family $\mathcal S_{\mathbf a} = (S_{\mathbf a}(t))_{t\geq 0} \subset \mathcal L(L^2(\RR) \oplus L^2(\RR))$ 
given by
\begin{equation*}
	S_{\mathbf a}(t) = \begin{pmatrix}
		L_\RR(t) & W_{\mathbf a}(t)
		\\0 & R_\RR(t)
	\end{pmatrix},
	\qquad t\geq0,
\end{equation*}
is a $C_0$-semigroup on $L^2(\RR) \oplus L^2(\RR)$. Moreover, since $(L_\RR(t))_{t\in \RR}$ and $(R_\RR(t))_{t\in \RR}$ are groups, it follows that $S_{\mathbf a}(t)$ is invertible for all $t>0$, 
so $\mathcal S_{\mathbf a}$ extends to 
a $C_0$-group on $L^2(\RR) \oplus L^2(\RR).$
On the other hand, recalling \eqref{tmatr},  observe that $V_{\mathbf a}(t)$ is the compression of $W_{\mathbf a}(t)$ to $L^2(\RR_+) \oplus L^2(\RR_+)$ for every $t\geq0$. So, $\mathcal T_{\mathbf a}$ is a compression of $\mathcal S_{\mathbf a}$ and then Corollary \ref{C0groupdilation_int} implies that $\mathcal T_{\mathbf a} \in \mathcal{SQC}(L^2(\RR_+) \oplus L^2(\RR_+))$, as required.
\end{proof}

Thus,  for Packel type semigroups $\mathcal T=(T(t))_{t \ge 0}$ considered above,
$T(t)$ may or may not be similar to a contraction for every $t>0,$ and $\mathcal T$ may or may not be similar to a quasi-contraction semigroup. In particular, this shows that conditions (i) and (ii) in Theorem \ref{SimConstTh} are independent of each other even for bounded semigroups, and may fail simultaneously, so that the formulation of Theorem \ref{SimConstTh}
is optimal, in a sense.

\subsection{(Non-)similarity to contractions and the Bhat-Skeide interpolating semigroup}\label{SkeideSect}

While the similarity theory for a single operator $T$ is well-developed, its continuous-parameter counterpart is far less so. It is thus
tempting to transfer the results obtained in the discrete setting
to the setting of $C_0$-semigroups.
For a long time, this task was quite challenging and in most of cases whenever such a transfer was possible, one had to repeat the proof closely following the arguments in the discrete case.
Recently, Bhat and Skeide  proposed in \cite[Lemma 2.4]{bhat2015pure} a neat way to circumvent the difficulty and given $T \in \mathcal L(H)$ to interpolate the powers of the tensor product $I \otimes T$ 
on $L^2(\mathbb T)\otimes H$ by an appropriate semigroup $\mathcal T$
on the same space.
Their technique was further developed, in particular, in \cite{dahya2024interpolation} and \cite{shalit2023semigroups}.
Thus, as it was proved in \cite{oliva2025tensor}, many examples of $T$ related to similarity to contractions appeared to have semigroup analogues, and the Bhat-Skeide construction was  used  in \cite{oliva2025tensor} to produce $C_0$-semigroups not similar to semigroups of contractions and having various additional properties arising in applications.

It is also natural to try to use the Bhat-Skeide technique in the present studies as well.
However, this technique has strong limitations as we show below.
{For $t \ge 0$, set}
\[
\varphi_t(z) = e^{-2\pi i t} z, \qquad  
{p_t(e^{2\pi i r}) = \chi_{[0,\,1-\{t\})}(r)}, \quad 
z = e^{2\pi i r}, \ r \in [0,1),
\]
so that $p_t$ is simply the indicator of the arc $[0,\,1-\{t\})$ on the unit circle $\TT$,
and where $\{t\} = t - \lfloor t \rfloor$ denotes the fractional part of $t$.

Define the bounded operators 
\begin{align*}
U(t) f = f \circ \varphi_t, \qquad P(t) f = p_t f, \qquad f \in L^2(\TT), \qquad t\geq 0,
\end{align*}
and note that $U(1) = P(1) = I$. Then the $C_0$-semigroup $\mathcal T_T$, given on $L^2(\TT) \otimes H$ by
\begin{align*}
T(t) &:= \left(U(t) \otimes I\right) \left(P(t) \otimes T^{\lfloor t \rfloor}+ \left(I - P(t) \right) \otimes T^{\lfloor t \rfloor +1}\right), \qquad t \geq 0,
\end{align*}
will be called the Bhat-Skeide interpolating semigroup associated with $T$.
A direct (but somewhat non-trivial)  verification shows that $T(n)=I\otimes T^n$ for all $n \in \mathbb Z_+,$
which explains our terminology.
Moreover,  $\mathcal T$ is contractive if and only if $T$ is a contraction, and $\mathcal T$ is bounded
if and only if $T$ is power bounded.
For justification of the above properties and further details, 
see \cite[Section 2.1]{dahya2024interpolation}, and also the original result \cite[Lemma 2.4]{bhat2015pure} and the discussion in \cite[Example 8.2]{shalit2023semigroups}. 
Moreover, it is crucial that $\mathcal T$ belongs to $\mathcal{SC}\left(L^2(\TT) \otimes H\right)$ if and only if $T$ is similar to a contraction, see Proposition \ref{BhatProp} below.
Thus, the Bhat-Skeide interpolation provides a convenient method for constructing semigroup counterparts of operators that are not similar to contractions, avoiding 
the lengthy and technical computations found, for instance, in
 \cite{packel1969semigroup} and Section \ref{PackelSect}. This leads to the following observation, that is somewhat surprising.

\begin{proposition}\label{BhatProp}
Let $H$ be a Hilbert space, $T \in \linearOp(H)$, and let $\mathcal T_T = (T(t))_{t\geq0}$ the Bhat-Skeide interpolating semigroup associated with $T$. Then 
\begin{itemize}
	\item [(i)] $\mathcal T_T \in \mathcal{SC}(H)$ if and only if $T$ is similar to a contraction.
	\item [(ii)] $\mathcal T_T$ is a quasi-contraction $C_0$-semigroup if and only if $T$ is a contraction.
	\item [(iii)] $\mathcal T_T \in \mathcal{SQC} \left(L^2(\TT) \otimes H\right)$.
\end{itemize}
\end{proposition}
\begin{proof}

The proof of (i) is direct since the ``if'' implication is obvious and the opposite implication
holds since $T$ can be identified with the restriction of $T(1) = I\otimes T$ to $\{\chi_{\TT}\} \otimes H$
(via the canonical identification $h \to 
\chi_{\TT} \otimes h$).

 To see (ii), note that for $0<t<1$ one has
 \[
 T(t) = (U(t)\otimes I)\,\big(P(t)\otimes I + (I-P(t))\otimes T\big),
 \]
 so that
 \[
 \|T(t)\| = \max\{1,\|T\|\}, \qquad 0<t<1.
 \]

To show that (iii) holds we construct an explicit equivalent norm $\|\cdot\|_{\rm{eq}}$ making $\mathcal T_T$ quasi-contractive. We identify $L^2(\TT)\otimes H$ with $L^2(\TT,H)$, which leads to the following representations
\begin{align}\label{BSRep}
	(U(t)\otimes I) f = f\circ \varphi_t, \qquad
	(P(t) \otimes I) f = p_t f, \qquad
	((I\otimes T) f)(e^{2\pi i r}) = T (f (e^{2\pi i r})), 
\end{align}
for all $t\geq 0$,  $r\in [0,1)$ and for every simple function $f \in L^2(\TT,H)$. As the set of simple functions is dense in $L^2(\TT, H)$, one infers that \eqref{BSRep} holds for every $f\in L^2(\TT)$ and that $\mathcal T_T$ can be identified with
\begin{align*}
	(T(t) f) (e^{2\pi ir}) &= T^{\lfloor t \rfloor} ((p_t f)  (\varphi_t(e^{2\pi i r}))) + T^{\lfloor t \rfloor+1}(((1-p_t) f) (\varphi_t(e^{2\pi i r})))
	\\ &= \begin{cases}
		T^{\lfloor t \rfloor + 1} (f(e^{2\pi i(r-t)})), \qquad & r \in [0, \{t\}),
		\\T^{\lfloor t \rfloor} (f(e^{2\pi i(r-t)})), \qquad & r \in [\{t\}, 1),
	\end{cases}
\end{align*}
where $f\in L^2(\TT, H)$ and $t\geq 0$. Put
$$\|f\|_{{\rm eq}}^2 := \int_0^1 \left(\|f(e^{2\pi i r})\|^2 + r \|T f(e^{2\pi i r})\|^2\right) \, dr, \qquad f \in  L^2(\TT, H).
$$
Clearly $\|\cdot\|_{{\rm eq}}$ defines an equivalent Hilbertian norm on $L^2(\TT, H)$. 
Moreover, if $t \in (0,1)$, then
\begin{align*}
	\|T(t) f\|_{{\rm eq}}^2 &= \int_t^{1} \left(\|f(e^{2\pi i (r-t)})\|^2 + r \|T(f(e^{2\pi i(r-t)}))\|^2\right)\, dr
	\\ & \qquad + \int_0^t \left(\| T (f(e^{2\pi i (r-t)}))\|^2 + r\|T^2 f (e^{2\pi i (r-t)})\|^2\right) \, dr
	\\&= \int_0^{1-t} \left(\|f(e^{2\pi i r})\|^2 + (r+t) \|T (f(e^{2\pi ir}))\|^2\right)\, dr
	\\ & \qquad + \int_{1-t}^1 \left(\| T (f(e^{2\pi i r}))\|^2 + (r-1+t)\|T^2 (f (e^{2\pi i r}))\|^2\right) \, dr
	\\ & \leq \int_0^{1-t} \left(\|f(e^{2\pi i r})\|^2 + r \|T (f(e^{2\pi i r}))\|^2\right) \, dr 
	\\ & \qquad + t\|T\|^2 \int_0^{1-t} \|f(e^{2\pi i r})\|^2 \, dr
	\\ & \qquad  + \left( \frac{1}{1-t}  +  t\|T\|^2 \right) \int_{1-t}^1 r \|T (f(e^{2\pi ir}))\|^2 \, dr 
	\\ & \leq \left(\frac{1}{1-t} + t\|T\|^2\right) \|f\|_{{\rm eq}}^2, \qquad f \in L^2(\TT, H), 
\end{align*}
so 
\begin{equation}\label{boundd}
	\|T(t)\|_{{\rm eq}} \leq \left(\frac{1}{1-t} + t\|T\|^2\right)^{1/2},  \qquad t\in (0,1). 
\end{equation}
Since  the right hand side of \eqref{boundd} is differentiable at zero,  the characterization of quasi-contractivity by \eqref{quasi_loc} implies the claim.

\end{proof}

Thus, in view of Theorem \ref{quasiContrTh_int}, the Bhat-Skeide interpolation is not suitable for constructing 
$C_0$-semigroups $\mathcal T = (T(t))_{t\ge0}$ that are not in $\mathcal{SC}(H)$, even if each $T(t)$ is similar to a contraction for $t>0$. In this case, one has to resort to either
classical Packel's ideas or use Le Merdy's type examples and their elaborations discussed in the next section.

\subsection{(Non-)similarity to quasi-contraction semigroups}\label{standardexamples}

A distinct feature of similarity theory for $C_0$-semigroups
is that standard similarity criterions such as, e.g., Rota's theorem do not in general hold
for continuous parameter semigroups, making thus the theory more involved
and motivating our studies. To justify this point, and to illustrate
the similarity criteria from the preceding section,  
we consider examples of $C_0$-semigroups $\mathcal T = (T(t))_{t\geq0}$ on a Hilbert space $H$ 
such that $\mathcal T \notin \mathcal{SC}(H)$ and, at the same time,
$T(t)$ is similar to a contraction for every $t>0$.
In this case, it follows from Corollary \ref{trichotomyCor} that $\mathcal C(\mathcal T)=\lim_{t \to 0}\mathcal C(T(t))=\infty,$
while $C(T(t)) <\infty$ for all $t \ge 0.$
Such examples are in sharp contrast with the situation of discrete parameter,
and showcase semigroups that are, in a sense, as far from $\mathcal{SC}(H)$ as possible.
In particular, such semigroups are not similar even to quasi-contraction semigroups.

\begin{example}\label{CounterEx}
\begin{enumerate} 
	\item [(i)] 
	Let $\mathcal S = (S(t))_{t\geq0}$ be a bounded $C_0$-semigroup on $H$ that does not belong to $\mathcal{SC}(H).$ 
	Define  a $C_0$-semigroup  on the infinite direct sum $\bigoplus_{n=1}^\infty H$ by setting $T(t) = \bigoplus_{n=1}^\infty S(nt)$ for all $t \ge 0.$ Chernoff \cite{chernoff1976two} proves that, for each $\lambda>0$, $e_{-\lambda} \mathcal T$ is an exponentially stable semigroup that does not belong to $\mathcal{SC}\left(\bigoplus_{n=1}^\infty H \right)$, i.e., $\mathcal T \notin \mathcal{SQC}\left(\bigoplus_{n=1}^\infty H \right)$.
	In particular, a continuous analogue of Rota's theorem for $\mathcal{SQC}(H)$ cannot hold as well.
	
	\item [(ii)] Let $\mathcal H$ be a separable Hilbert space and let $(h_j)_{j\in \NN}$ be a conditional basis of $\mathcal H$. Define
	\begin{equation}\label{leMerdySem}
		T_{LM}(z) h_j: = e^{-2^j z} h_j, \qquad j\in \NN, \, \Real z >0,
	\end{equation}
	and extend  $T_{LM}$ to $\mathcal H$ by linearity and density,
	 denoting the extension by the same symbol. Le Merdy proved in \cite{lemerdy2000bounded} that $\mathcal T_{LM} = (T_{LM}(t))_{t\geq0}$ is a (sectorially) bounded holomorphic $C_0$-semigroup of angle $\pi/2$ 
	 which is, in addition,
	 immediately compact.\footnote{Recall that a $C_0$-semigroup $\mathcal T = (T(t))_{t\geq0}$ is said to be immediately compact if $T(t)$ is a compact operator for all $t>0$.}
	 Moreover, it is easy to see that $\mathcal T_{LM}$ is exponentially stable.
	  At the same time,  $\mathcal T_{LM}$ does not belong to $\mathcal{SC}(\mathcal H)$.
	 Thus,
	continuous counterparts to either Nagy's or Rota's discrete similarity criteria fail dramatically.

	On the other hand, let $\generator$ be the generator of $\mathcal T_{LM}.$,
	Then $\generator^{-1}$, given by $\generator^{-1} h_j = -2^{-j} h_j$ for all $j\in \NN$, is a bounded operator on $H.$ It is  sectorial of angle $0$ since $\generator$ is so. 
	Hence, $\generator^{-1}$ generates a bounded holomorphic semigroup on $H$ given by 
	$$
	e^{\generator^{-1}z} h_j = e^{-2^{-j}z} h_j, \qquad j \in \NN, \, \Real z >0.
	$$
	Note that $(e^{\generator^{-1} t})_{t\geq0}$ is a quasi-contraction $C_0$-semigroup (since its generator is bounded) that does not belong to $\mathcal{SC}(\mathcal H)$. Indeed, if it belonged to $\mathcal{SC}(\mathcal H)$, then $\generator^{-1}$ would be similar to a dissipative operator, which would imply that $\generator$ is similar to a dissipative operator, arriving at contradiction since $\mathcal T_{LM} \notin \mathcal{SC}(\mathcal H)$; see also \cite[Remark 5.5.4]{haase2006functional}. 
	
	This counterexample shows that the assumption (ii) in Theorem \ref{quasiContrTh_int} cannot be omitted even in the case when $\mathcal T$ is bounded holomorphic of maximal angle $\pi/2$ and has bounded generator (so that $\mathcal T$ extends to a group on $H$). 
	
	\item [(iii)] Consider  the Riemann-Liouville semigroup $\mathcal T_{RL} = (T_{RL}(t))_{t\geq0}$ on $L^2[0,1]$ given by
	\begin{equation}\label{RiemannEq}
		(T_{RL}(t)f)(x) = \frac{1}{\Gamma(t)} \int_0^x (x-y)^{t-1} f(y)\, dy, \qquad x \in [0,1], \, f \in L^2[0,1], \, t >0. 
	\end{equation}
	It is well known that $\mathcal T_{RL}$ is an immediately compact, quasi-nilpotent and bounded holomorphic of angle $\frac{\pi}{2}$ $C_0$-semigroup on $L^2[0,1]$; see, for instance, \cite[Theorem 3.1 \& Theorem 3.4]{carvalho2022riemann}.
	
	Using arguments similar to Proposition \ref{OT24quasi}, we proved in \cite{oliva2025tensor} that the tensor product $\mathcal T_{RL} \otimes \mathcal T_{LM} = (T_{RL}(t) \otimes T_{LM}(t))_{t\geq0}$, where $\mathcal T_{LM} = (T_{LM}(t))_{t\geq0}$ is as in \eqref{leMerdySem}, is an immediately compact, quasi-nilpotent and bounded holomorphic of angle $\frac{\pi}{2}$ $C_0$-semigroup on $L^2[0,1] \otimes \auxHilbert$ that is not similar to a quasi-contraction one. 
	This example is a substantial strengthening	of Le Merdy's example.
	
	\item [(iv)] 
	Another improvement of Le Merdy's example can be given via a slightly different tensor product construction. Let $\mathcal R = (R(t))_{t\geq0}$ be the right shift semigroup on $L^2[0,1].$ 
	It was proved in \cite{oliva2025tensor} that $\mathcal R \mathcal T_{RL} \otimes \mathcal T_{LM} = (R(t) T_{RL}(t) \otimes T_{LM}(t))_{t\geq0}$, where $\mathcal T_{RL} = (T_{RL}(t))_{t\geq 0}$ is the Riemann-Liouville semigroup \eqref{RiemannEq} and $\mathcal T_{LM} = (T_{LM}(t))_{t\geq0},$ defined by \eqref{leMerdySem}, is a nilpotent and immediately compact $C_0$-semigroup on $L^2[0,1] \otimes \auxHilbert$ that is not similar to a quasi-contraction $C_0$-semigroup. 
\end{enumerate}
\end{example}

It follows from Theorem \ref{quasiContrTh_int} that none of the examples mentioned above, except for $\mathcal T_{\generator^{-1}}$, is similar to a quasi-contraction semigroup. This property can also be proved directly,  without invoking Theorem \ref{quasiContrTh_int}.

\section{Similarity criteria  and their applications to control theory}\label{DefectOperSect}

In this section, we provide criteria for a $C_0$-semigroup $\mathcal T=(T(t))_{t \ge 0}$ to be similar to a contraction, quasi-contraction, or isometric $C_0$-semigroup.
The criteria are formulated  in terms of the
size of ``weighted'' orbits $C\mathcal T(\cdot)$ of $\mathcal T,$
where a linear operator $C$, interpreted as an operator weight, may not be even closable
and plays a role analogous to that of an observation operator in control theory.
In Section \ref{ObservSubsect} we make this precise and
discuss implications of our results for control theory.

\subsection{Similarity to contractions in terms of operator means and its interplay with control theory}

The first result describes $\mathcal T\in \mathcal{SC}(H)$ in terms of local two-sided $L^2$-bounds
for $CT(\cdot),$
and it apparently addresses the largest class of $C$ making the bounds meaningful.
\begin{theorem}\label{PetitcunotProp}
	Let $\mathcal T = (T(t))_{t\geq0}$ be a $C_0$-semigroup on a Hilbert space $(H, \|\cdot\|)$ with generator $\generator$. Then $\mathcal T \in \mathcal{SC}(H)$ if and only if
	there exist a Hilbert space $(K, \|\cdot\|_K)$ and a densely defined linear operator $C: \dom(C) \subseteq H \to K$ such that
	 the following conditions hold:
				\begin{itemize}
			\item[(i)] $T(t)(\dom(C)) \subseteq \dom(C)$ for all $t\geq0$.
			\item [(ii)] For every $h \in \dom(C)$, the mapping $F_h: [0,\infty) \to K$ defined by
			\begin{equation}\label{fh}
			F_h(t)=CT(t)h, \qquad t \geq 0,
			\end{equation}
			is Bochner measurable.
			\item [(iii)] There exist $\alpha,\beta>0$ such that
			\begin{equation}\label{PetitcunotAssumpt}
				\alpha \|h\|^2 \leq \|T(t) h\|^2 + \int_0^t \|C T(s) h\|^2_{K} \, ds \leq \beta \|h\|^2, \qquad  h\in \dom(C), \, t\geq 0.
			\end{equation}
		\end{itemize}
	If \textnormal{(i)--(iii)} hold, then $C$ can be chosen to be a bounded operator from $\dom(\generator)$ to $K,$
	where $\dom(\generator)$ is endowed with the graph norm.
\end{theorem}

\begin{proof}
	Assume first that  (i)--(iii) hold. Taking into account (i) and (ii),
	for every $t\geq0,$
	define a linear operator $G(t): \dom(C) \subseteq H \to L^2((0,\infty),K)$ by
	$$
	G(t)h =
	CT(t-\cdot)h  \chi_{[0,t]}(\cdot), \qquad h \in \dom(C),
		$$
		{which is Bochner square-integrable by (ii)–(iii).}
		By (ii) and (iii), $G(t)$ continuously extends to a bounded operator from $H$ to $L^2((0,\infty),K),$ denoted also by $G(t).$
	
	Let  $\mathcal R = (R(t))_{t\geq0}$ be the right shift semigroup on $L^2((0,\infty),K)$ and let
	$\mathcal S = (S(t))_{t\geq0} \subset \mathcal L\left(H \oplus L^2((0,\infty),K)\right) $ be given by
	\begin{align}\label{SGrabowskiSem}
		S(t) = \begin{pmatrix}
		T(t) 	&	0
		\\ G(t)	& R(t)
	\end{pmatrix}, \qquad t \geq0.
	\end{align}
	Note that $\mathcal S = (S(t))_{t\geq0}$ satisfies the semigroup law since 
	\begin{align*}
		R(s) G(t) h + G(s) T(t) h &= C T(t+s - \cdot) h \, \chi_{[s,s+t]} + C T(s-\cdot)T(t) h\, \chi_{[0,s]}
		\\ &= C T(s+t-\cdot) h \, \chi_{[0,s+t]} = G(s+t)h, \qquad h \in \dom(C), \, s,t\geq0.
	\end{align*}
	From the dominated convergence theorem it follows that $\lim_{t\to0^+} G(t)h = 0$ for every $h\in \dom(C),$ and then for all $h \in H$ since  $\dom(C)$ is dense in $H.$ Thus, $\mathcal S$ is a $C_0$-semigroup on $H \oplus L^2((0,\infty),K)$.

Let $\|\cdot\|_{\oplus}$ stand for the norm on $H \oplus L^2((0,\infty),K).$
	Taking into account (ii), for all $h \in \dom(C)$, $f\in L^2((0,\infty),K)$, and $t\geq0$, we have
	\begin{equation*}
		\begin{aligned}
			\|S(t)(h,f)\|_\oplus^2 &= \|T(t)h\|^2 + \|G(t)h\|^2 + \|R(t)f\|^2
			\\&= \|T(t)h\|^2 + \int_0^t \|C T(s)h\|^2_{K} \, ds + \|f\|^2.
		\end{aligned}
	\end{equation*}
	Hence by (iii),
	\begin{align}\label{NagyAsEq}
		\min\{1,\alpha\} \|(h,f)\|_\oplus^2 \leq \|S(t)(h,f)\|_\oplus^2 \leq \max\{1,\beta\} \|(h,f)\|_\oplus^2,
	\end{align}
	so that $\mathcal S$ satisfies \eqref{nagy_isom_in}.
	Thus, 
	$\mathcal S$ is similar to a $C_0$-semigroup of isometries. As $\mathcal S$ is a dilation of $\mathcal T$, taking into account Remark \ref{similar_dil}, we conclude that $\mathcal T \in \mathcal{SC}(H)$. \medskip 
	
	Now, assume that $\mathcal T \in \mathcal{SC}(H).$
		Let $\mathcal S = (S(t))_{t\geq0}$ be a $C_0$-semigroup of contractions on $H$ such that $S(t) := R T(t) R^{-1}$, $t\geq0,$ for an invertible $R \in \linearOp(H).$ Without loss of generality, we can assume that $\|R\|\leq 1.$
	Since the generator $\widetilde \generator$ of $\mathcal S$ is dissipative,
 we have
	\begin{align*}
		\Real \langle R^\ast R \generator h, h \rangle &= \Real \langle R \generator h , R h\rangle = \Real \langle \widetilde \generator Rh, Rh \rangle \leq 0, \qquad h \in \dom(\generator),
	\end{align*}
	so that $R^\ast R \generator$ is dissipative as well. Hence if we define
	$$\langle h,f \rangle_K := - \langle R^\ast R \generator h, f \rangle  - \langle h, R^\ast R \generator f\rangle, \qquad h,f \in \dom(\generator),
	$$
	then $\langle \cdot, \cdot \rangle_K$ is a positive semidefinite, sesquilinear form on $\dom(\generator).$
	Setting $N: = \{ h \in \dom(\generator) \, : \, \langle h, h \rangle_K = 0\},$
	define the Hilbert space $K$ as the completion of quotient space $\dom(\generator)/N$ in the norm induced by $\langle \cdot, \cdot \rangle_K,$
	and  $C: \dom(\generator) \to K$ as the composition of the quotient mapping from $\dom(\generator)$ to $\dom(\generator)/N$ with the inclusion mapping from $\dom(\generator)/N$ to $K$.
	Clearly, (i) holds. Moreover, since
	$$\|C h\|_K^2 = -2 \Real \langle R^\ast R \generator h, h\rangle \leq 2 \|R\|^2 \|\generator h\| \|h\| \leq 2  \left(\|\generator h\|^2 + \|h\|^2\right)^{1/2}, \qquad h \in \dom(\generator),
	$$
	we conclude that  $C \in \mathcal L (\dom(\generator), K).$
	In particular, the mapping $t \mapsto CT(t)h $ is continuous on $\mathbb R_+$ for every $h \in H$, hence Bochner measurable, and (ii) follows.
	
	To prove (iii), observe that
	\begin{align*}
		\int_0^t \|C T(s) h\|^2_{K} \,ds &= - \int_0^t \left(\langle R^\ast R \generator T(s) h, T(s) h\rangle + \langle T(s) h, R^\ast R \generator T(s) h\rangle \right)\, ds
		\\ &= - \int_0^t \left(\langle \widetilde \generator S(s) Rh, S(s) R h\rangle + \langle S(s) R h, \widetilde \generator S(s) R h \rangle\right)\, ds
		\\ &= - \int_0^t \frac{d}{d s}\left(\|S(s) Rh\|^2\right)\, ds
		\\&= \|Rh\|^2 - \|S(t)Rh\|^2 = \|Rh\|^2 - \|RT(t)h\|^2, \qquad h \in \dom(\generator), \, t \geq 0.
	\end{align*}
	Therefore, for all $t\geq0$ and $h\in \dom(C)$,
	we have
	\begin{align*}
		\|T(t) h\|^2 + \int_0^t \|C T(s) h \|^2_{K} \,ds =& \|T(t) h\|^2 +  \|Rh\|^2 - \|RT(t)h\|^2\\
		\le& \|T(t)h\|^2 + \|Rh\|^2 \leq \left(\sup_{s\geq0} \|T(s)\|^2  + 1 \right) \|h\|^2,
	\end{align*}
	and
	\begin{align*}
		\|T(t) h\|^2 + \int_0^t \|C T(s) h \|^2_{K} \,ds \geq \|Rh\|^2 \geq \|R^{-1}\|^{-2} \|h\|^2,
	\end{align*}
	where we used  that $\|R\|\leq 1$.
	This yields \eqref{PetitcunotAssumpt}, and finishes the proof.
\end{proof}

Note that the proof given above follows ideas similar to the ones in \cite{grabowski1996admissible} (see also \cite{gao1999infinite}), where the authors use more technically involved arguments, although with a slightly different set-up and with a different goal.

\begin{remark}
	Theorem \ref{PetitcunotProp} is a continuous counterpart of the discrete similarity criterion given in the essentially inaccessible preprint \cite{eckstein1978operators}.
	Its formulation and proof
	can be found in \cite[Chapter 3]{petitcunot2008problemes}. The criterion states that $T \in \mathcal L(H)$
	is similar to a contraction if and only if there are $\alpha,\beta >0$  and $C\in \mathcal L(H)$ such that
	\[
	\begin{aligned}
				\alpha \|h\|^2 \leq  \|T^n h\|^2 + \sum_{j=0}^n \|C T^j h\|^2 \, ds &\leq \beta \|h\|^2, \qquad h \in H, \, n \in \mathbb Z_+.
			\end{aligned}
	\]
	This result appears to be equivalent to the one
	 proved \cite[Theorem II]{bello2019operator}, claiming that $T \in \mathcal L(H)$ is similar to a contraction if and only if there exist 	$C \in \linearOp(H)$	and $\alpha,\beta >0$  such that
		\begin{equation}\label{DiscrBelloEq}
			\alpha \|h\|^2 \leq \limsup_{n\to \infty} \|T^n h\|^2 + \sum_{j=0}^\infty \|C T^{j} h\|^2 \leq \beta \|h\|^2, \qquad h\in H.
		\end{equation}
		We omit justification of this fact, but give a similar argument showing that
		\eqref{PetitcunotAssumpt} can be recast as
		\begin{equation}\label{PetitcunotAssumpt2}
			\begin{aligned}
				\alpha \|h\|^2 \leq  \limsup_{t\to \infty} \|T(t)h\|^2 + \int_0^\infty \|C T(s) h\|^2_{K} \, ds &\leq \beta \|h\|^2, \qquad  h \in \dom(C).
			\end{aligned}
		\end{equation}
		Indeed, it is easy to see that \eqref{PetitcunotAssumpt} implies \eqref{PetitcunotAssumpt2}.
		For the converse implication, note first that since \eqref{PetitcunotAssumpt2} yields $\sup_{t\ge 0} \|T(t)\|<\infty$ by the uniform boundedness principle, the second inequality in \eqref{PetitcunotAssumpt}
		follows. To obtain the first inequality in \eqref{PetitcunotAssumpt}
		it suffices to note that
		\begin{align*}
			\alpha \|h\|^2 &\leq \limsup_{s\to \infty} \|T(s) h\|^2 + \int_0^\infty \|C T(s) h \|^2_{K} \, ds \\
		&\leq  \left(\sup_{s>0} \|T(s)\| \right)^2 \|T(t) h\|^2 + \int_0^{t} \|C T(s) h\|^2_{K} \, ds + \int_0^\infty \|C T(s) T(t)h\|^2 \, ds
		\\
		&\leq  \left(\beta + \sup_{s>0} \|T(s)\|^2\right) \|T(t) h\|^2 + \int_0^{t} \|C T(s) h\|^2_{K} \, ds, \qquad h \in H, \, t \geq 0.
		\end{align*}
\end{remark}

\begin{remark}\label{nonclosableObsRem}
	Even if  $\mathcal T=(T(t))_{t\geq0}$ belongs to $\mathcal{SC}(H),$
	 the operator $C$ constructed in the proof of Theorem \ref{PetitcunotProp} may not be closable. Indeed, let $\mathcal T$ be the right shift semigroup on $L^2[0,1],$ with generator $\generator$.
	Recall that
	$$\dom(\generator) = \{f \in L^2[0,1] \, : \, f \mbox{ absolutely continuous in } [0,1], \,\, f' \in L^2[0,1] \mbox{ and } f(0) = 0\},
	$$
	and $\generator f = - f'$ for $f \in \dom(\generator)$. Using the same notation as in the proof of Proposition \ref{PetitcunotAssumpt}, observe that
	\begin{align*}
		\|f\|_K^2 &= -\langle \generator f, f\rangle - \langle f, \generator f\rangle 
		= \int_0^1 \frac{d }{d s} \left(|f(s)|^2\right)\, ds 
		= |f(1)|^2, \qquad f \in \dom(\generator).
	\end{align*}
	Thus,
	$$N = \{ f \in \dom(\generator) \, : \, \|f\|_K = 0\} = \{f \in \dom(\generator) \, : \, f(1) = 0\},
	$$
		$K= \dom(\generator)/N$ can be identified with $\CC,$  and $C: \dom(\generator) \to \CC$ is then given by $Cf = f(1)$, $f\in \dom(\generator)$.
	Now if $(f_n)_{n\in \NN} \subset \dom(C),$ $\lim_{n\to \infty} f_n = 0$ in $L^2[0,1]$,  and $f_n(1) = 1,$ then $\lim_{n\to \infty} C f_n = 1 \neq 0 = C (\lim_{n \to \infty} f_n)$, hence $C$ is not closable.
\end{remark}

On the other hand, the simple proposition below
shows that assumption (i) in Theorem \ref{PetitcunotProp} along with closability of $C$
imply (ii). 
\begin{proposition}\label{meas-suff}
	Let $\mathcal T = (T(t))_{t\geq0}$ be a $C_0$-semigroup on a Hilbert space $H$. Let $K$ be a Hilbert space, and let 
		$C:\dom(C) \subseteq H \to K$ 
	be a closable linear operator on $K$ such that:
	\begin{itemize}
		\item [(i)]$T(t)(\dom(C)) \subseteq \dom(C)$, $t\geq0$.
		\item [(ii)] For each $h \in \dom(C)$, the mapping $F_h: [0,\infty) \to K$ given by \eqref{fh}
		is almost-separably valued (i.e., its essential range is contained in a separable subspace of $K$).
	\end{itemize}
	Then $F_h$
		is Bochner measurable for each $h \in \dom(C)$.
\end{proposition}
\begin{proof}
	Let $\overline C$ be  the closure of $C,$ and let $\|\cdot\|_{\dom(\overline{C})}$ be the Hilbert space graph norm on $\dom(\overline{C})$, so that $\|h\|_{\dom(\overline C)}^2 = \|h\|^2 + \|\overline Ch\|_K^2$, $h \in \dom(\overline C)$. Fix $h \in \dom(C)$ and let $\Omega$ be a null-set of $\RR_+$ such that the
	set $\{CT(t)h \,:\, t \in \RR_+\setminus \Omega\}$ is separable in $K$. As the sets $\{T(t)h \, : \, t \in \RR_+ \setminus \Omega\}$ and $\{CT(t)h \, : \, t\in \RR_+\setminus \Omega\}$ are separable in $H$ and $K$, respectively, it follows that $\{T(t)h \, : \, t\in \RR_+\setminus \Omega\}$ is a separable subset of
	$\left(\dom(\overline C), \|\cdot\|_{\dom(\overline C)}\right)$. 
	Hence,  the closure of ${\operatorname{span}}\{T(t) h \, : t\in \RR_+ \setminus \Omega\}$ in
	$\left(\dom(\overline C), \|\cdot\|_{\dom(\overline C)}\right)$ is a separable Hilbert space, and we denote it by $Y_h$. Consequently, the inclusion mapping $\embedding_h: (Y_h, \|\cdot\|_{\dom(\overline C)}) \to (Y_h, \|\cdot\|_H)$ is a Borel isomorphism; see, for instance, \cite[Corollary 4.5.5]{srivastava1998course}. Since $\overline C: (\dom(\overline C), \|\cdot\|_{\dom(\overline C)}) \to K$ is continuous, it follows that
	$\overline C \embedding_h^{-1}:Y_h \to H$ is a Borel map, so the map
	$t \mapsto \overline C  \embedding_h^{-1} T(t)h = CT(t)h$ from $\RR_+ \setminus \Omega$ to $K$ is Borel as well.  Hence the latter map is Bochner measurable as $\{CT(t)h \, : \, t \in \RR_+\setminus \Omega\}$ is separable by hypothesis. Since $\Omega$ is a null-set of $\RR_+$, we conclude that
	$F_h$ is Bochner measurable.
\end{proof}

\begin{remark}
	The argument given in the proof of Proposition \ref{meas-suff} also applies, with some modifications, if $C$ is closable in a stronger norm $\|\cdot\|_X$ on $\dom(C)$ such that, for each $h\in \dom(C)$, $\{T(t)h \, : \, t \geq0\}$ is contained in a standard Borel\footnote{A measurable space $X$  is said to be standard Borel if it is isomorphic to a Borel subset of a Polish space.} subset of $\left(\dom(C), \|\cdot\|_X\right)$. This more general case covers, in particular, the example given in Remark \ref{nonclosableObsRem}.
\end{remark}

We proceed with a counterpart  of Theorem \ref{PetitcunotProp} for quasi-contraction $C_0$-semigroups.

\begin{theorem}\label{PetitcunotQuasiProp}
	Let $\mathcal T = (T(t))_{t\geq0}$ be a $C_0$-semigroup on a Hilbert space $(H, \|\cdot\|),$ with generator
	$\generator$. Then $\mathcal T \in \mathcal{SQC}(H)$ if and only if there exist a Hilbert space $(K, \|\cdot\|_K)$,  a densely defined linear operator $C:\dom(C)\subseteq H \to K$ satisfying
	the conditions {\rm (i)} and {\rm (ii)} of Theorem \ref{PetitcunotProp}, and
	 $\alpha, \beta>0$ and $\tau>0$ such that
		 	\begin{equation}\label{PetitQuasi}
		 		\alpha \|h\|^2 \leq  \|T(\tau)h\|^2 + \int_0^\tau \|C T(s) h\|^2_{K}\, ds \leq \beta \|h\|^2, \qquad h\in \dom(C).
		 	\end{equation}
		If (i), (ii), and \eqref{PetitQuasi} hold, and  $\dom(\generator)$ is equipped with the graph norm, then $C$ can be chosen to be a bounded operator from $\dom(\generator)$ to $K$.
\end{theorem}
\begin{proof}
To prove the ``if'' part, using \eqref{PetitQuasi} and arguing as in the proof of Theorem \ref{PetitcunotProp}, we consider the $C_0$-semigroup $\mathcal S = (S(t))_{t\geq0}$ on $H \oplus L^2((0,\infty),K)$ given by \eqref{SGrabowskiSem}.
	Then	\eqref{PetitQuasi} implies that $S(\tau)$ is bounded from below, and thus,
	by Proposition \ref{boundedBelowContractProp}(ii), $\mathcal S$ is similar to a quasi-contraction $C_0$-semigroup. Since $\mathcal S$ is a dilation of $\mathcal T$, Remark \ref{similar_dil} implies that $\mathcal T \in \mathcal{SQC}(H)$. \medskip

The ``only if'' implication is direct.	Assuming that
$e_{-\lambda} \mathcal T$ belongs to   $\mathcal{SC}(H)$
for some $\lambda \in \mathbb R,$
and applying Theorem \ref{PetitcunotProp} to $e_{-\lambda} \mathcal T$ and any fixed $\tau>0,$
we deduce (i), (ii), and \eqref{PetitQuasi}.
\end{proof}
\begin{remark}\label{remark_c}
Theorems \ref{PetitcunotProp} and \ref{PetitcunotQuasiProp} hold if their assumptions (i)
 are replaced by $\dom(\generator) \subset \dom(C)$,
and (ii) and (iii) hold for all $h \in \dom(\generator)$ rather than $\dom(C)$.
For the proof it suffices to replace $C$ by its restriction to $\dom(\generator).$
\end{remark}

\subsection{Similarity to isometries in terms of operator means}

There are several classical results characterizing when a linear operator 
$\generator$ with $\sigma(\generator)\subseteq \RR$ is similar to a 
self-adjoint operator, i.e. when $i\generator$ generates a uniformly 
bounded $C_0$-group.
 Some of these results are formulated in terms of two-sided bounds 
 for means of $L^2$-norms of semigroup orbits, see for instance \cite[Section 2, p. 43]{malamud1985criterion}, \cite[Remark 3, p. 15]{naboko1984conditions}, \cite[Theorem 3.1]{van1983operators}. However, unlike the approach in the preceding section, 
they do not employ operator weights.

As a prototype, Naboko \cite{naboko1984conditions} showed that 
similarity of $\generator$ to a self-adjoint operator follows if there 
exist constants $\alpha,\beta>0$ such that
\begin{align}\label{NabokoEq}
	\alpha \|h\|^2 \leq \limsup_{\varepsilon \to 0^+} \, \varepsilon \int_{-\infty}^\infty \|(\xi + i \varepsilon - \generator)^{-1} h\|^2 \, d\xi \leq \beta \|h\|^2, \qquad h \in H,
\end{align}
see also Theorem \ref{NabokoTh}.
Apart from their theoretical importance,
these results appeared to be useful in the study of concrete differential operators, see, for example, \cite{karabash2009similarity, malamud2003similarity} and the references therein.
Our next result extends these criteria to the semigroup setting, providing
two-sided bounds for orbit means with operator weights, in the spirit of 
Theorems \ref{PetitcunotProp} and \ref{PetitcunotQuasiProp},
 cf. \cite[Proposition 1.15]{kubrusly1997introduction}.

\begin{theorem}\label{KubruslyNabokoProp}
	Let $(H, \|\cdot\|)$ be a Hilbert space and let $\generator$ be a linear operator in $H$. Then the following are equivalent.
	\begin{itemize}
		\item [(i)] $\generator$ generates a  $C_0$-semigroup $\mathcal T$ similar to an isometric $C_0$-semigroup.
		\item [(ii)] $\generator$ generates a $C_0$-semigroup $\mathcal T = (T(t))_{t\geq0}$, and there exist a Hilbert space $(K, \|\cdot\|_K)$ and  a densely defined, linear operator $C$ on $K$ 
		satisfying assumptions {\rm (i)} and {\rm (ii)} of Theorem \ref{PetitcunotProp} such that
			\begin{align}\label{kubruslyEq1}
				\alpha \|h\|^2 \leq  \limsup_{t \to \infty} \frac{1}{t} \int_0^t  \|C T(s) h\|^2_{K} \, ds &\leq \beta \|h\|^2,
			\end{align}
			or
			\begin{align}\label{kubruslyEq2}
				\alpha \|h\|^2 \leq  \liminf_{t \to \infty} \frac{1}{t} \int_0^t   \|C T(s) h\|^2_{K} \, ds &\leq \beta \|h\|^2,
			\end{align}
			for some $\alpha,\beta >0$ and  all $h\in \dom(C)$.
		\item [(iii)] $\sigma(\generator) \subseteq \{z \in \CC \, :  \, \Real z \leq 0\}$ and there exist a Hilbert space $K$ and a densely defined, linear operator $C: \dom(C) \subseteq H \to K,$
		such that the following holds.
		\begin{enumerate}
			\item [(a)] $\dom(\generator) \subseteq \dom(C)$.
			\item [(b)] For each $\varepsilon>0$ and $h\in H$, the mapping from $\RR$ to $H$ given by
			\begin{equation}\label{esigmah}
				\xi \mapsto  C (\varepsilon + i\xi - \generator)^{-1} h,
			\end{equation}
			is (Bochner) measurable.
			\item [(c)] There exist $\alpha,\beta>0$ such that
			\begin{equation}\label{resAssumpt}
				\begin{aligned}
				\limsup_{\varepsilon \to 0^+} \varepsilon \int_{-\infty}^\infty \|C (\varepsilon + i \xi - \generator)^{-1} h \|^2_K \, d \xi &\leq \beta \|h\|^2,
				\\ \liminf_{\varepsilon \to 0^+} \varepsilon \int_{-\infty}^\infty \|C (\varepsilon + i \xi - \generator)^{-1} h \|^2_K \, d \xi &\geq \alpha \|h\|^2,
				\end{aligned}
			\end{equation}
			for all $h \in H$.
		\end{enumerate} 
	\end{itemize}
\end{theorem}

\begin{remark}
	Note that, a priori, the operator $E$ in (iii) need not be densely defined, although the density of its domain will be established in the course of the proof.
\end{remark}

\begin{proof}
	(i) $\iff$ (ii): If $\generator$ generates a $C_0$-semigroup $\mathcal T$ that is similar to an isometric one, then 
	 it follows that
	\eqref{kubruslyEq1} and \eqref{kubruslyEq2} hold with $C=I$.
	
	To prove the opposite implication, assume that $\generator$ generates a $C_0$-semigroup $\mathcal T = (T(t))_{t\geq0}$ satisfying \eqref{kubruslyEq1}.
	Then, one has
	\begin{align}
		\|T(\tau)h\|^2
		& \geq \frac{1}{\beta} \limsup_{t\to \infty} \frac{1}{t} \int_\tau^{t+\tau} \|C T(s) h \|^2_K \, ds
		\notag \\ &= \frac{1}{\beta} \limsup_{t\to \infty} \left( \frac{1}{t} \int_0^{t+\tau} \|C T(s) h\|^2_K\,ds  -  \frac{1}{t} \int_0^\tau \|C T(s) h\|^2_K \,ds\right)
	\label{2bounds}	\\
		\notag
		 &= \frac{1}{\beta} \limsup_{t\to \infty} \frac{1}{t} \int_0^{t+\tau} \|C T(s) h\|^2_K\,ds
		\geq \frac{\alpha}{\beta} \|h\|^2, \qquad h\in \dom(C), \, \tau > 0.
	\end{align}
	A similar argument yields
	$\|T(\tau) h\|^2 \leq \frac{\beta}{\alpha}\|h\|^2$ for all $h \in \dom(C)$ and $\tau>0$. Since $\dom(C)$ is dense in $H$, one obtains $\sqrt{\frac{\alpha}{\beta}}\|h\| \leq  \|T(\tau) h\| \leq \sqrt{\frac{\beta}{\alpha}} \|h\|$ for all $h\in H$ and $\tau>0$. Hence, \eqref{nagy_isom_in} implies that $\mathcal T$ is similar to a semigroup of isometries, and our claim follows.
	If \eqref{kubruslyEq2} holds then the argument is completely analogous.\medskip

	(i) $\implies$ (iii): This is Theorem \ref{NabokoTh} with $C= I$. \medskip
	
	(iii) $\implies$ (i): Fix $h \in H$. By \eqref{resAssumpt} and the resolvent identity, one gets
	\begin{align*}
		& \limsup_{\varepsilon \to 0^+} \, \varepsilon \int_{-\infty}^\infty \|(\varepsilon + i\xi -\generator)^{-1}h \|^2 \, d \xi \\ 
		\leq & \frac{1}{\alpha} \limsup_{\varepsilon \to 0^+} \, \varepsilon \int_{-\infty}^\infty \left(\liminf_{\delta \to 0^+} \, \delta \int_{-\infty}^\infty\|C (\delta + i\eta  - \generator)^{-1}(\varepsilon+ i\xi -\generator)^{-1}h \|^2_K \, d\eta \right) d \xi \\
		\leq & \frac{1}{\alpha} \limsup_{\varepsilon \to 0^+} \left(\liminf_{\delta \to 0^+} \, \varepsilon \delta \int_{-\infty}^\infty   \int_{-\infty}^\infty\|C (\delta + i\eta - \generator)^{-1}(\varepsilon + i\xi -\generator)^{-1}h \|^2_K \, d\eta d \xi\right) \\
		\le & \frac{2}{\alpha} \limsup_{\varepsilon \to 0^+} \bigg(\liminf_{\delta \to 0^+} \, \varepsilon \delta \int_{-\infty}^\infty   \int_{-\infty}^\infty \frac{1}{(\varepsilon - \delta)^2 + (\xi - \eta)^2} \bigg(\|C (\delta + i\eta - \generator)^{-1}h \|^2_K  \\ 
		& \qquad \qquad + \|C(\varepsilon + i\xi -\generator)^{-1}h \|^2_K
		\bigg)\, d\eta d \xi\bigg).
	\end{align*}
	Furthermore, note that by \eqref{resAssumpt} and Fubini's theorem,
	\begin{align*}
		&\liminf_{\delta\to 0^+} \, \varepsilon \delta \int_{-\infty}^\infty   \int_{-\infty}^\infty \frac{1}{(\varepsilon - \delta)^2 + (\xi - \eta)^2} \|C (\delta + i\eta - \generator)^{-1}h \|^2_K \, d\eta d\xi \\
		=&\liminf_{\delta\to 0^+} \frac{\pi\varepsilon \delta}{\varepsilon-\delta} \int_{-\infty}^\infty \|C (\delta + i\eta - \generator)^{-1}h \|^2_K \, d \eta
		\leq \pi \beta \|h\|^2, \qquad \xi > 0, \, h \in H,
	\end{align*}
	and in addition,
	\begin{align*}
		&\liminf_{\delta \to 0^+} \, \varepsilon \delta \int_{-\infty}^\infty   \int_{-\infty}^\infty \frac{1}{(\varepsilon - \delta)^2 + (\xi - \eta)^2} \|C (\varepsilon + i\xi - \generator)^{-1}h \|^2_K \, d\eta d\xi \\
		=& \liminf_{\delta \to 0^+} \, \frac{\pi \varepsilon \delta}{\varepsilon - \delta} \int_{-\infty}^\infty \|C (\varepsilon + i\xi - \generator)^{-1} h \|^2_K \, d \xi = 0,
	\end{align*}
{for all $h\in H$ and sufficiently small $\varepsilon>0$ (so that the integrals above are finite).}
	Thus, taking the two above displays into account,
	we obtain that
	\begin{align*}
		\limsup_{\varepsilon \to 0^+} \, \varepsilon \int_{-\infty}^\infty \|(\varepsilon + i\xi -\generator)^{-1}h \|^2 \, d \xi \leq 2\pi \frac{\beta}{\alpha} \|h\|^2, \qquad h \in H.
	\end{align*}
	A similar reasoning using $\limsup_{\delta \to 0^+}$ instead of $\liminf_{\delta \to0^+}$
	and relying on the elementary inequality $\|a-b\|^2\ge (1/2)\|a\|^2 -\|b\|^2, a, b \in H,$ shows that
	\begin{align*}
		\limsup_{\varepsilon \to 0^+} \, \varepsilon \int_{-\infty}^\infty \|(\varepsilon + i\xi -\generator)^{-1}h \|^2 \, d \xi \geq \frac{\pi}{2} \frac{\alpha}{\beta} \|h\|^2, \qquad h \in H.
	\end{align*}
	It then follows from Proposition \ref{NabokoTh} that $\generator$ generates a $C_0$-semigroup $\mathcal T$ that is similar to an isometric one.
\end{proof}

\begin{remark}
	If the operator $C$ is closed, then the equivalence (ii) $\iff$ (iii)  in Proposition \ref{KubruslyNabokoProp} follows by standard arguments involving Plancherel's theorem and two-sided inequalities between Cesàro means and Abel means; see, for instance, \cite[Lemma 1.1]{van1983operators}.
\end{remark}

\subsection{Observability and controllability via similarity}\label{ObservSubsect}

We continue with elaborating several significant applications of the abstract similarity
criteria from the preceding section to characterization of observability and controlabillity
of infinite-dimensional control systems.
To make the applications accessible to non-experts, all of the necessary background
is developed in Appendix \ref{ObservControlSect}.

We start with a criterion for observability/controllability in infinite
time, unifying and strengthening several results in the literature (see the discussion following the theorem).
\begin{theorem}\label{infiniteTimeObsTh}
	Let $(H, \|\cdot\|)$ be a Hilbert space and let $\mathcal T = (T(t))_{t\geq0}$ be a $C_0$-semigroup on $H,$ with generator $\generator$. Then the following are equivalent.
	\begin{itemize}
		\item [(i)]  $\mathcal T \in \mathcal{SC}(H)$ and $\lim_{t \to \infty}\|T(t)h\|=0$ for all $h \in H.$
		\item [(ii)] There exist a Hilbert space $(K, \|\cdot\|_K)$ 
		and an infinite time admissible observation operator
		$C \in \mathcal L(\dom(\generator),K)$
		  for $\mathcal T$ 
		such that the pair $(\generator, C)$ is infinite time exactly observable.
	\end{itemize}
	Moreover, the following are also equivalent.
	\begin{itemize}
		\item [(iii)] $\mathcal T \in \mathcal{SC}(H)$ and $\lim_{t \to \infty}\|T^*(t)h\|=0$ for all $h \in H.$
		\item [(iv)] There exist a Hilbert space $(K, \|\cdot\|_K)$
		and an infinite time admissible control operator 
		 $B \in \mathcal L(K, H_{-1})$
		for $\mathcal T$
		such that the pair $(\generator, B)$ is infinite time exactly controllable.
		\end{itemize}
\end{theorem}
\begin{proof}
	(i) $\implies$ (ii): This follows from Theorem \ref{PetitcunotProp} in view of $\lim_{t\to \infty} \|T(t)h\|^2 = 0$ for any $h \in H$. \medskip

	(ii) $\implies$ (i): Let
	$$\|h\|_{\mathscr H}^2 := \int_0^\infty \|C T(t) h\|^2_K\, dt, \qquad h \in \dom(\generator).
	$$
	Clearly, $\|\cdot\|_\mathscr H$ is a Hilbertian norm on $\dom(\generator)$, and since the pair $(\generator,C)$ is infinite time exactly observable,
	it is equivalent to the original norm on $H$. Thus,
	$\mathscr H := (H, \|\cdot\|_{\mathscr H})$ is a Hilbert space
	and
	$$\|T(t)h\|_\mathscr H \leq \|h\|_\mathscr H, \qquad t\geq 0, \, h \in \dom(\generator),
	$$
	so that
	$\mathcal T$
	is similar to a contraction semigroup on $H$. Moreover, the dominated convergence theorem implies that
	$\lim_{t\to \infty} \|T(t) h \|_{\mathscr H} = 0$ for all $h \in \dom(\generator)$, hence the same holds
	for all $h \in \mathscr H$ since $\dom(\generator)$ is dense in $\mathscr H$.
	Since $\|\cdot \|_{\mathscr H}$ is equivalent to the original norm, the convergence holds on $H$ as well.
		\medskip
	
	(iii) $\iff$ (iv): This is a straightforward consequence of the equivalence (i) $\iff$ (ii) and the duality of controllability and observability, see Theorem \ref{DualObsTh}.
\end{proof}

While we omit the discussion of discrete variants of the above results, it is instructive to recall that Helton proved in \cite{helton1974discrete}  that a bounded operator $A$ is similar to a strongly stable contraction, with its adjoint operator $A^\ast$ also strongly stable, if and only if there exist \textit{bounded} admissible operators $C$ (observation) and $B$ (control) such that $(A,C)$ is exactly observable and $(A,B)$ is exactly controllable (see \cite{helton1974discrete} for precise definitions
and more details).

In the setting of Theorem \ref{infiniteTimeObsTh}, Grabowski and Callier showed in their seminal work \cite{grabowski1996admissible} that infinite time exact observability of $(\generator, C)$ (with unbounded $C$) is equivalent to similarity of $\mathcal T$ to a contraction semigroup, assuming that $\mathcal T$
is exponentially stable (see Theorem $3.1$ and Corollary $3.1$ in \cite{grabowski1996admissible}). Their argument is rather intricate, as they aim at showing the dissipativity for the generator of $\mathcal{S}$ (given by \eqref{SGrabowskiSem}) rather than establishing  upper and lower norm bounds for orbits of $\mathcal S$, which is considerably more direct.
Furthermore, \cite[Proposition 5.1]{haak2012exact} shows that infinite time exact observability implies similarity to a contraction semigroup, though without addressing strong stability. Among other related results in the literature, we mention  \cite[Proposition 2]{levan1990left} where one proves that exponential stability of $\mathcal T$, combined with the exact controllability of $(\generator, I)$, implies similarity of $\mathcal S$ to a semigroup of contractions.

The next result is our main application underlining the limitations of the notion of \emph{finite time} observability/controllability for well-posed control systems.
\begin{theorem}\label{finiteTimeObsTh}
	Let $(H, \|\cdot\|)$ be a Hilbert space and let $\mathcal T = (T(t))_{t\geq0}$ be a $C_0$-semigroup on $H,$ with generator $\generator$. Then the following are equivalent.
	\begin{itemize}
		\item [(i)] $\mathcal T \in \mathcal{SQC}(H)$.
		\item [(ii)] There exist a Hilbert space $(K, \|\cdot\|_
		K)$ and 
		a finite time admissible observation operator 
		$C \in \mathcal L(\dom(\generator), K)$ 
		for $\mathcal T$
		such that the pair $(\generator, C)$ is finite time exactly observable.
		\item [(iii)] There exist a Hilbert space $K$ and a finite time admissible control operator $B \in \linearOp(K, H_{-1})$ 
		for $\mathcal T$
		such that the pair $(\generator,B)$ is finite time exactly controllable.
	\end{itemize}
	\end{theorem}
\begin{proof}
	(i) $\implies$ (ii): From (i) it follows that there exists
	$\lambda\geq0$ such that $e_{-\lambda} \mathcal T$ belongs to $\mathcal{SC}(H)$, thus we assume without loss of generality that $e_{-\lambda} \mathcal T$ is a contraction $C_0$-semigroup. As $e_{-(\lambda+1)}\mathcal T$ is also a contraction $C_0$-semigroup, Theorem \ref{PetitcunotProp} implies that there exist $\alpha>0$, $\beta>0$, a Hilbert space $K$ and $C\in \mathcal L(\dom(\generator),K)$ such that
	$$\alpha \|h\|^2 \leq \|e^{-(\lambda +1)t} T(t)h\|^2 + \int_0^t \|Ce^{-(\lambda+1)s}T(s)h\|^2_K\, ds \leq \beta \|h\|^2, \qquad h \in \dom(\generator), \, t \geq 0.
	$$
	Fix $\tau>0$ with $e^{-2\tau} < \alpha$. Then
	$$
	\int_0^\tau \|CT(s)h\|^2_K \, ds \leq \beta e^{2(\lambda+1)\tau} \|h\|^2, \qquad h \in \dom(\generator),
	$$
	so $C$ is an admissible observation operator for $\mathcal T$ in finite time  $\tau$. Moreover,
	since $\|e^{-\lambda t}T(t)\|\leq 1$ for all $t \ge 0,$ we have
	\begin{align*}
		(\alpha - e^{-2\tau}) \|h\|^2 & \leq \alpha \|h\|^2 - e^{-2\tau} \|e^{-\lambda \tau}T(\tau) h\|^2
		\\ &\leq \int_0^\tau e^{-2(\lambda+1)s}\|CT(s)h\|^2_K\, ds
		\\ & \leq \int_0^\tau \|CT(s)h\|^2_K \, ds, \qquad h \in \dom(\generator).
	\end{align*}
	Hence the pair $(\generator,C)$ is finite time exactly observable, and the claim follows.
	\medskip
	
	(ii) $\implies$ (i): This is an immediate corollary of Theorem \ref{PetitcunotQuasiProp} and the definition of finite time exact observability. \medskip
	
	(ii) $\iff$ (iii): As in the proof of Theorem \ref{infiniteTimeObsTh},  this is a direct consequence of the equivalence (i)$\iff$(ii) combined with   
	Theorem \ref{DualObsTh}.
	\end{proof}

Versions of Theorem \ref{finiteTimeObsTh} for \emph{bounded} $C$ where obtained in \cite[Theorem 1.3]{dolecki1977general} and \cite[Theorem 1]{zwart2013left}. However, most of PDE applications require $C$ unbounded, and $C$ could be even non-closable even if $(C, \generator)$ is observable/controllable (as in Remark \ref{nonclosableObsRem}). If $C$ is bounded, then, as proved in \cite[Theorem 1]{zwart2013left}, the finite time observability of $(C, \generator)$ is equivalent to the boundedness from below (equivalently, left-invertibility)
 of $\mathcal T$ and thus, by \cite[Theorem 7.3]{batty2017lower}, to the extendability of $\mathcal T$ to a $C_0$-group into a larger Hilbert space $K\supseteq H.$ Thus, in view of Corollary \ref{C0groupdilation_int}, the difference between the cases of bounded and unbounded $C$'s in Theorem \ref{finiteTimeObsTh} can be interpreted as the difference between the properties of $\mathcal T$ to be dilatable and to be extendable to a $C_0$-group on larger Hilbert space containing $H$.

\begin{remark}
	The results above provide a different proof of Rota's type theorem for quasi-contraction $C_0$-semigroups,   \cite[Proposition 5.2]{oliva2025tensor}, stating that every exponentially stable, quasi-contraction $C_0$-semigroup belongs to $\mathcal{SC}(H)$. Indeed, let $\mathcal T$ be an exponentially stable $C_0$-semigroup in  $\mathcal{SQC}(H)$. By Theorem \ref{finiteTimeObsTh}, there exists a finite time admissible observation operator $C$ for $\mathcal T$ such that the pair $(\generator, C)$ is finite time exactly observable.
	As $\mathcal T$ is exponentially stable, $C$ is an \emph{infinite} time admissible observation operator
 for $\mathcal T$ making  $(\generator,C)$ exactly observable in \textit{infinite} time; see, for instance, \cite[Proposition 6.5.2]{tucsnak2009observation}. It then follows from Theorem \ref{infiniteTimeObsTh} that $\mathcal T \in \mathcal{SC}(H)$, as required.
\end{remark}

\section{Final remarks}

In this section, we provide some additional remarks on the results obtained in this paper,
which might be of importance for further developments.

\subsection{Local criterion for similarity to quasi-contraction $C_0$-semigroups}

In the same spirit as \cite{holbrook1973operators}, see also Proposition \ref{HolSemigroup}, we provide a simple criterion for the similarity to quasi-contraction $C_0$-semigroups. This is motivated by potential applications of our main result (Theorem \ref{quasiContrTh_int}) and the fact that similarity to a quasi-contraction $C_0$-semigroup depends only on the behaviour of the semigroup $\mathcal T = (T(t))_{t\geq0}$ for values of $t$ close to the origin, see Theorem \ref{nearzero}.

\begin{proposition}\label{zeroagain}
	Let $\mathcal T = (T(t))_{t\geq0}$ be a $C_0$-semigroup on a Hilbert space $H$. Assume there exist a Hilbert space $K$, a surjective operator $A \in \linearOp(K, H)$, and a family of bounded operators $\mathcal S = (S(t))_{t\geq0}$ on $K$ such that
	\begin{equation}
	\|S(t)\| = 1 + O(t), \qquad \mbox{ as } t \to 0,
	\end{equation}
	and
	\begin{equation}\label{commut}
	\|T(t) A - A S(t)\| = O(t), \qquad \mbox{ as } t \to 0.
	\end{equation}
	Then $\mathcal T \in \mathcal{SQC}(H)$.
\end{proposition}
\begin{proof}
Since	the induced operator $\mathscr A \in \linearOp(K \ominus \ker A ,  H)$ is invertible,
 the semigroup $\mathcal V =(V(t))_{t\geq0}$ on $K \ominus \ker A$ defined as $V(t) = \mathscr A^{-1} T(t) \mathscr A$ for $t\geq0$, is similar to $\mathcal T$. Then we have $A = \mathscr A \proj_{K \ominus \ker A}$, $\mathscr A = A\restriction_{K \ominus \ker A}$, 
	and
	\begin{align*}
			\|V(t) - \proj_{K \ominus \ker A} S(t)\restriction_{K \ominus \ker A}\| &= \|\mathscr A^{-1} T(t) \mathscr A - \proj_{K \ominus \ker A} S(t)\restriction_{K \ominus \ker A}\| 
			\\&\leq \|\mathscr A^{-1} \| \|T(t) A - A S(t) \|  = O(t), \quad \mbox{ as } t\to 0.
	\end{align*}
	Hence, we conclude
	\begin{align*}
			\|V(t)\| &\leq \|\proj_{K \ominus \ker A} S(t) \restriction_{K \ominus \ker A}\| + \|V(t) - \proj_{K \ominus \ker A} S(t)\restriction_{K \ominus \ker A}\| \leq 1 + O(t), \quad \mbox{ as } t\to 0,
		\end{align*}
		which by \eqref{quasi_loc} is equivalent to $\mathcal T \in \mathcal{SQC}(H).$
\end{proof}
Clearly, the assumptions of Proposition \ref{zeroagain} are also necessary. However, they lie quite close to the conclusion of the proposition, and it seems plausible that they could be replaced by more revealing conditions.

\subsection{Semigroups discontinuous at $0$}

Let $H$ be a Hilbert space.  A family of operators $\mathcal T = (T(t))_{t\geq0} \subset \linearOp(H)$ is said to be a \textit{degenerate semigroup} if the following holds:
\begin{itemize}
	\item [(i)] $T(s+t) = T(s)T(t)$ for $s,t \geq 0$, and $T(0)= I$.
	\item [(ii)] The mapping $t \mapsto T(t)$ is strongly continuous from $(0,\infty)$ to $\linearOp(H)$.
	\item [(iii)] $\limsup_{t\to0^+} \|T(t)\| < \infty$.
\end{itemize}

Most of the arguments in this paper extend naturally to the broader context of degenerate semigroups with only minor modifications, since these arguments do not rely on the existence of a generator or on strong continuity at the origin. 
For ease of reference, we state below the corresponding extension of Theorems \ref{quasiContrTh_int} and \ref{SimConstTh}, which will play a crucial role in the study of infinite tensor products of semigroups in  \cite{oliva2025tensorbis}.

\begin{theorem}\label{nonstrongTh}
	Let $H$ be a Hilbert space and let $ \mathcal T = (T(t))_{t\geq0}$ be a degenerate semigroup on $H$. 
	Assume that there exist $\lambda >0$ and $\tau>0$ such that
	\begin{enumerate}
		\item [(i)] $e_{-\lambda} \mathcal T \in \mathcal{SC}(H);$
		\item [(ii)] $T(\tau)$ is similar to a contraction.
	\end{enumerate}
	Then $\mathcal T$ belongs to $\mathcal{SC}(H).$
	Moreover, the similarity constant $\mathcal C(\mathcal T)$ of $\mathcal T$ satisfies
	\begin{equation}\label{boundSimEq2}
		\mathcal C(\mathcal T) \leq  \sqrt{2} \, \mathcal C(e_{-\lambda} \mathcal T) \frac{e^{2\lambda}-1}{2\lambda} + 2 \sqrt{2} \, C(T(\tau)) M^2 \max \{1,\sqrt{\tau}\},
	\end{equation}
	where $M = \sup_{t\in [0,\tau]} \|T(t)\|$.
\end{theorem}

A different and more direct approach to studying the similarity problem for a degenerate semigroup $\mathcal T$ on $H$ is the following. Set
\begin{equation}\label{HandK}
	N := \{h \in H \, :\, T(t)h = 0 \ \text{for all } t>0\}, 
	\qquad 
	K := \{h \in H \, : \, \lim_{t\to0^+} T(t) h = h\}.
\end{equation}
By adapting the arguments from the proof of \cite[Corollary~2.2]{arendt2001approximation} (see also \cite[Proposition~2.1]{arendt2001approximation}), we conclude that $N$ and $K$ are closed $\mathcal T$-invariant subspaces such that $H$ is the (in general, non-orthogonal) direct sum $H=N \dot{+} K,$ and that the restriction 
\[
\mathcal T \restriction_K = (T(t)\restriction_K)_{t\geq0} \subset \linearOp(K)
\] 
is a $C_0$-semigroup on $K$.

Let $P$ be the projection onto $K$ along $N$, define the equivalent Hilbertian norm $\|\cdot \|_{\mathscr H}$ on $H$ by
\[
\|h\|_{\mathscr H}^2=\|Ph\|^2 + \|(I-P)h\|^2, \qquad h \in H,
\]
and set $\mathscr H=(H, \|\cdot\|_{\mathscr H})$. Then $K=N^{\perp}$ in $\mathscr H$, and if $O: H \to \mathscr H$ denotes the identity isomorphism, the projection $OPO^{-1}$ onto $K$ is orthogonal.
Passing to $\mathscr H$ and employing the
description \eqref{eqNormSim} of similarity constants in terms of equivalent norms,
it is direct to show that 
$\mathcal T \in \mathcal{SC}(H)$ (resp.\ $\mathcal{SQC}(H)$) 
if and only if 
$\mathcal T\restriction_K \in \mathcal{SC}(K)$ (resp.\ $\mathcal{SQC}(K)$). 
Moreover,
\begin{equation}\label{SimConstRestriction}
	\mathcal C (\mathcal T\restriction_K) \;\leq\; \mathcal C (\mathcal T) 
	\;\leq\; \|O\|\,\|O^{-1}\| \,\mathcal C (\mathcal T\restriction_K).
\end{equation}
Using this observation, most of our results extend to the setting of degenerate semigroups, with the upper bounds for the similarity constants (such as \eqref{boundSimEq2}) increasing by the factor $\|O\|\|O^{-1}\|$. 
On the other hand, since $\lim_{t \to 0}T(t) = P$ strongly, 
an application  of \eqref{const_simm}
yields
\[
\mathcal C (\mathcal T) \geq\sup_{t>0} \|T(t)\| \geq \|P\|.
\]
{As $C(\mathcal T\restriction_K)$ is independent of $P,$ no universal constant $\kappa \geq 1$, independent of
$P,$}
 can satisfy 
\[
\mathcal C (\mathcal T) \leq \kappa \, \mathcal C(\mathcal T\restriction_K)
\]
for all degenerate semigroups $\mathcal T$. This loss of uniformity in the bounds prevents a formal reduction of Theorem~\ref{nonstrongTh} to Theorem~\ref{SimConstTh}, and becomes particularly relevant in the study of similarity properties of \emph{families} of operator semigroups. Such situations arise naturally, for instance, in the context of infinite tensor products of operator semigroups, see \cite{oliva2025tensorbis}.

\subsection{Similarity by complete boundedness}\label{cbdSubsection}

 It is well known that a bounded operator 
 $T$ on a Hilbert space 
 $H$ is similar to a contraction if and only if 
 the homomorphism $p \to p (T)$ defined on the algebra $\mathcal P$ of polynomials $p$ extends to a completely
 bounded map from the disc algebra $A(\mathbb D)$ into $\mathcal L(H)$; see \cite[Theorem 9.11]{paulsen2002completely} or \cite[Theorem 4.13]{pisier2001similarity}. The corresponding result in the context of $C_0$-semigroups follows as a direct corollary and is essentially contained in the proof of \cite[Theorem 5.1]{fackler2015regularity}; see also \cite[Theorem 2.2]{lemerdy1996dilation}. Since this criterion plays a fundamental role in the discrete theory, we formulate explicitly its semigroup analogue.

To this end, we need a few definitions. Recall that $\CC^+ = \{z \in \CC \, : \, \Real z>0\}$ and let $\mathcal A$ be the (non-closed) subalgebra of $L^\infty (\CC^+)$ generated by the functions $\{e_{-t} \, : \, t \geq 0\}$, that is,
$$\mathcal A = \left\{\sum_{k=1}^N a_k e_{-t_k} \, :\, N\in \NN, \, a_1,\ldots, a_N \in \CC, \, t_1, \ldots, t_N \geq 0\right\}.
$$
 Clearly, $ \|f\|_{\mathcal A} = \sup_{z\in \CC^+} |f(z)|$ for all $f\in \mathcal A$.
Note that $\mathcal A$ can be identified with a subalgebra of multiplication operators on $L^2(\CC^+),$ and thus $\mathcal A \subset \mathcal L(L^2(\CC^+)).$

Now, given a $C_0$-semigroup $\mathcal T = (T(t))_{t\geq0}$ on a Hilbert space $H$, we define the homomorphism $\Theta_\mathcal T: A \to \linearOp(H)$ by setting
$$
\Theta_\mathcal T \left(e_{-t} \right) = T(t), \qquad t \geq 0,
$$
and extending this map to $A$ by linearity.

The following result characterizes  $\mathcal T \in \mathcal{SC}(H)$ in terms of the complete boundedness of $\Theta_\mathcal T$, thus invoking a well-known concept in operator space theory, see \cite[Chapters 8 \& 9]{paulsen2002completely} and \cite[Chapter 3]{pisier2001similarity} for precise definitions and background.

\begin{theorem}\label{completelyBoundedTh}
	{Let $\mathcal T = (T(t))_{t\geq0}$ be a $C_0$-semigroup on a Hilbert space $H$. Then $\mathcal T \in \mathcal{SC}(H)$ if and only if  
	$\Theta_\mathcal T: \mathcal A \to \linearOp(H)$ extends to a completely bounded map from
	the closure $\bar{\mathcal A}$ of $A$ into $\mathcal L(H).$}
	\end{theorem}

The ``if'' part of the claim follows directly from Paulsen's similarity criteria (see, for instance, \cite[Theorem 9.1]{paulsen2002completely}), whereas the arguments for the ``only if'' implication can be found in the proof of \cite[Theorem 5.1]{fackler2015regularity}. Theorem \ref{completelyBoundedTh} has a version formulated in more traditional terms. Let $\widehat M(\mathbb R_+)$ denote  the Banach algebra of  Laplace transforms $\widehat \mu$ of bounded Borel measures $\mu$ on $\mathbb R_+,$ 
equipped with the  (induced) variation norm.
Then one can prove a similar statement 
by replacing  $\Theta_\mathcal T$ with the homomorphism
$\widehat M(\mathbb R_+)\ni \mu \to \int_{0}^\infty T(t)\,d\mu(t),$
where the integral is understood in the strong operator topology.

At present, unlike in the discrete setting, we have not found a way to exploit  Theorem \ref{completelyBoundedTh} or its variants.

\appendix

\section{Similarity to isometric semigroups}\label{NabokoResult}

Here we prove a generalization of a result stated by Naboko in \cite[Section 3]{naboko1984conditions} for generators of bounded $C_0$-groups (i.e. operators similar to skew-adjoint ones). Such a generalization is needed in the proof of Theorem \ref{KubruslyNabokoProp}, and  since the proof in \cite{naboko1984conditions} is omitted,  we provide a complete argument  in the case of isometric $C_0$-semigroups
and with more general assumptions on resolvent two-sided bounds.

\begin{theorem}\label{NabokoTh}
	Let $H$ be a Hilbert space, and let $\generator$ be a linear operator in $H$. Then the following are equivalent.
	\begin{itemize}
		\item [(i)] $\generator$ generates a $C_0$-semigroup $\mathcal T$  similar to an isometric one.
		\item [(ii)] $\sigma(\generator) \subseteq \{z \in \CC \, : \, \Real z \leq 0\}$ and there exist $\alpha>0$ and $\beta>0$ such that either
		\begin{align*}
			\alpha \|h\|^2 \leq \limsup_{\varepsilon \to 0^+} \, \varepsilon \int_{-\infty}^\infty \|(\varepsilon + i \xi - \generator)^{-1} h \|^2 \, d\xi \leq \beta \|h\|^2,
		\end{align*}
		or
		\begin{align*}
			\alpha \|h\|^2 \leq \liminf_{\varepsilon \to 0^+} \, \varepsilon \int_{-\infty}^\infty \|(\varepsilon + i \xi - \generator)^{-1} h \|^2 \, d\xi \leq \beta \|h\|^2
		\end{align*}
		holds for all $h \in H$.
	\end{itemize}
\end{theorem}

\begin{proof}
	(i) $\implies$ (ii): Even though this implication is well-known, we include a proof of it for the sake of completeness. Assume $\generator$ generates a $C_0$-semigroup $\mathcal T = (T(t))_{t\geq0}$ that is similar to an isometric semigroup. Then there exist $\alpha>0$ and $\beta>0$ such that
	$$\alpha \|h\| \leq \|T(t)h\| \leq \beta \|h\|, \quad h \in H, \, t\geq 0.
	$$
	As $\omega_0(\mathcal T) = 0$, one has
	$\sigma(\generator) \subseteq \{z \in \CC \, : \, \Real z \leq 0\}$ 
		and, given $\varepsilon>0$ and $h\in H$, the mapping from $\RR$ to $H$ defined by $\xi \mapsto (\varepsilon + i\xi - \generator)^{-1}h$ is the Fourier transform of $t \to e^{-\varepsilon t}T(t)h, t \ge 0,$ extended by zero to $\mathbb R.$
		Hence, by the Plancherel theorem,
	\begin{align}\label{FourierEq}
		\varepsilon \int_{-\infty}^\infty \|(\varepsilon + i\xi - \generator)^{-1} h\|^2 \, d\xi &= 2\pi \varepsilon \int_0^\infty e^{-2\varepsilon t} \|T(t)h\|^2 \, dt
		\leq \pi \beta^2 \|h\|^2
	\end{align}
	for all $h \in H$ and $\varepsilon > 0$.	Similarly, one proves
	\[
		\varepsilon \int_{-\infty}^\infty \|(\varepsilon + i\xi - \generator)^{-1} h\|^2 \, d\xi  \geq \pi \alpha^2 \|h\|^2, \qquad h \in H, \, \varepsilon>0,
	\]
	and (ii) follows. \medskip
	
	(ii) $\implies$ (i): First, we prove that $E$ is the generator of a bounded $C_0$-semigroup on $H.$ Assume that (ii) holds with bounds involving $\limsup_{\varepsilon \to 0^+}$ (the proof is analogous for the inequalities involving $\liminf_{\varepsilon \to 0^+}$). By the resolvent identity, we have
	\begin{align*}
		(w - \generator)^{-1} (z - \generator)^{-n} = \frac{1}{(z-w)^n} (w-\generator)^{-1} - \sum_{j=1}^n \frac{1}{(z-w)^{n-j+1}} (z-\generator)^{-j},
	\end{align*}
	for all $n \in \NN$, and all $z, w \in \CC$ with $\Real z >0$, $\Real w>0$ and $z\neq w$.
	Therefore,
	\begin{equation}\label{resNEq}
	\begin{aligned}
		\|(z-\generator)^{-n}h\|^2 &\leq  \frac{1}{\alpha} \limsup_{\varepsilon \to 0^+} \, \varepsilon \int_{-\infty}^\infty \|(\varepsilon + i\xi - \generator)^{-1} (z-\generator)^{-n} h\|^2 \, d \xi \\
		& = \frac{1}{\alpha} \limsup_{\varepsilon \to 0^+} \, \varepsilon \int_{-\infty}^\infty \Bigg\| \frac{1}{(z-\varepsilon-i\xi)^n} (\varepsilon + i \xi - \generator)^{-1} h
		\\ & \qquad \quad - \sum_{j=1}^n \frac{1}{(z-\varepsilon-i\xi)^{n-j+1}} (z-\generator)^{-j}h \Bigg\|^2 \, d \xi
	\end{aligned}
	\end{equation}
	for all  $n \in \NN,$ $h \in H$ and $z$ with $\Real z>0.$
	Note that, for all $n,j \in \NN$ with $j\leq n$,
	\begin{equation}\label{auxEqApp1}
		\begin{aligned}
		&\limsup_{\varepsilon \to 0^+} \, \varepsilon \int_{-\infty}^\infty \left\| \frac{1}{(z-\varepsilon-i\xi)^{n-j+1}} (z-\generator)^{-j}h \right\|^2 \, d\xi
		\\ =& \|(z-\generator)^{-j}h\|^2 \limsup_{\varepsilon \to 0^+} \,  \int_{-\infty}^\infty \frac{\varepsilon}{((\Real z - \varepsilon)^2 + \xi^2)^{2(n-j+1)}} \, d\xi = 0, \quad \Real z > 0, \, h \in H,
			\end{aligned}
	\end{equation}
	where the last step follows from the dominated convergence theorem.
	On the other hand,
	\begin{equation}\label{auxEqApp2}
		\begin{aligned}
		& \limsup_{\varepsilon \to 0^+} \, \varepsilon \int_{-\infty}^\infty  \left\|\ \frac{1}{(z - \varepsilon - i\xi)^n} (\varepsilon + i \xi-\generator)^{-1}h \right\|^2 \, d\xi \\
		=&\limsup_{\varepsilon \to 0^+} \, \varepsilon \int_{-\infty}^\infty \frac{1}{((\Real z - \varepsilon)^2 + (\Imag z - \xi)^2)^n} \|(\varepsilon + i \xi-\generator)^{-1}h\|^2 \, d\xi \\
		\leq &  \limsup_{\varepsilon \to 0^+} \, \frac{\varepsilon}{(\Real z - \varepsilon)^{2n}} \int_{-\infty}^\infty \|(\varepsilon + i \xi-\generator)^{-1}h\|^2 \, d\xi \\
		\leq & \frac{\beta}{(\Real z)^{2n}} \|h\|^2, \qquad n \in \NN,\, \Real z >0, \, h \in H.
		\end{aligned}
	\end{equation}
	Expanding the squared norm in the last term of \eqref{resNEq}, applying H{\"o}lder's inequality, and taking into account \eqref{auxEqApp1} and \eqref{auxEqApp2}, we obtain
	\begin{equation}\label{gen}
		\|(z-\generator)^{-n}h\|^2 \leq \frac{\beta}{\alpha} \frac{1}{(\Real z)^{2n}} \|h\|^2, \qquad n \in \NN, \Real z > 0, \, h \in H.
	\end{equation}
	The  inequality \eqref{gen} with $n=1$ implies  that $-\generator$ is a sectorial operator, so it is densely defined; see, for instance, \cite[Proposition 2.1.1]{haase2006functional}.  In view of \eqref{gen} $\generator$ satisfies the assumptions of  Hille-Yosida's theorem and is thus generates a bounded $C_0$-semigroup $\mathcal T = (T(t))_{t\geq0}$ on $H$.
	
	Now, using the first equality of \eqref{FourierEq}, we have
	\begin{equation}\label{belowIneq}
		\begin{aligned}
		\|h\|^2 &\leq \beta \limsup_{\varepsilon \to 0^+} \, \varepsilon \int_{-\infty}^\infty \|(\varepsilon +i \xi - \generator)^{-1} h\|^2 \, d \xi
		\\ &= 2\pi \beta \limsup_{\varepsilon \to 0^+} \,  \varepsilon \int_0^\infty e^{-2 \varepsilon t} \|T(t)h\|^2\, dt
		\\ &= 4\pi \beta \limsup_{\varepsilon \to 0^+} \,  \varepsilon^2 \int_0^\infty e^{-2\varepsilon \tau} \int_0^\tau \|T(t)h\|^2 \, dt d\tau, \qquad h \in H.
		\end{aligned}
	\end{equation}
	Fix $h \in H$ and $\delta > 0$, and let $s_0>0$ be such that
	$$\frac{1}{s} \int_0^s \|T(t)h\|^2 \, dt \leq \delta + \limsup_{\tau\to \infty} \frac{1}{\tau} \int_0^\tau \|T(t)h\|^2 \, dt, \qquad s > s_0.
	$$
	Using \eqref{belowIneq} and $\limsup_{\varepsilon\to0^+} \varepsilon^2 \int_0^{s_0} e^{-2\varepsilon \tau} \int_0^\tau \|T(t)h\|^2 \, dt d\tau = 0$, one obtains
	\begin{equation}\label{belowIneq2}
		\begin{aligned}
			\|h\|^2 &\leq 4\pi \beta \limsup_{\varepsilon\to 0^+}\, \varepsilon^2 \int_{s_0}^\infty e^{-2\varepsilon \tau} \int_0^\tau \|T(t)h\|^2 \, dt d\tau 
			\\ &  \leq 4\pi \beta \left( \delta + \limsup_{\tau\to \infty} \frac{1}{\tau} \int_0^\tau \|T(t)h\|^2 \, dt \right) \left( \limsup_{\varepsilon \to 0^+}\, \varepsilon^2 \int_{s_0}^\infty \tau e^{-2\varepsilon \tau} \, d\tau \right)
			\\ &= \pi \beta \left( \delta + \limsup_{\tau\to \infty} \frac{1}{\tau} \int_0^\tau \|T(t)h\|^2 \, dt \right),
		\end{aligned}
	\end{equation}
	where we used that $\limsup_{\varepsilon \to 0^+}\, \varepsilon^2 \int_{s_0}^\infty \tau e^{-2\varepsilon \tau} \, d\tau = 1/4$. As $\delta >0$ and $h\in H$ were arbitrary, \eqref{belowIneq2} implies that
	\begin{equation}\label{bound_lower}
	\|h\|^2 \leq \pi \beta \limsup_{\tau\to \infty} \frac{1}{\tau} \int_0^\tau \|T(t)h\|^2 \, dt
	\end{equation}
	holds for all $h \in H$.
	Next, let $h \in H$ and $s>0$ be fixed. Observe that
	\begin{equation}\label{shiftinv}
	\limsup_{\tau\to \infty} \frac{1}{\tau} \int_0^\tau \|T(t)h\|^2 \, dt =	\limsup_{\tau\to\infty}\frac1\tau\int_s^{s+\tau}\|T(t)h\|^2\,dt
	\end{equation}
	since
		\begin{align*}
	\left| \frac{1}{\tau} \int_0^\tau \|T(t)h\|^2 \, dt -\frac1\tau\int_s^{s+\tau}\|T(t)h\|^2\,dt\right|
	=\left| \frac{1}{\tau} \int_0^s \|T(t)h\|^2 \, dt -\frac1\tau\int_\tau^{s+\tau}\|T(t)h\|^2\,dt\right|
	\end{align*}
	and, in view of the boundedness of $\mathcal T,$ the latter expression goes to zero as $\tau \to \infty.$ Hence,
	letting  $M:=\sup_{t>0} \|T(t)\|$ and using the shift-invariance property \eqref{shiftinv}, we obtain
	\begin{equation}\label{bound_low}
	\limsup_{\tau\to\infty}\frac{1}{\tau}\int_0^\tau \|T(t)h\|^2\,dt=\limsup_{\tau\to\infty}\frac1\tau\int_s^{s+\tau}\|T(t)h\|^2\,dt \le M^2\|T(s)h\|^2.
	\end{equation}
	Combining \eqref{bound_lower} with \eqref{bound_low}, we infer that
	$\|T(s)h\|\ge \frac{1}{M\sqrt{\pi\beta}}\|h\|.$
	This holds for all  $h$ and $s$ since their choice was arbitrary.
	Hence, $\mathcal T$ satisfies \eqref{nagy_isom_in}, and
	is thus  similar to an isometric semigroup   by Sz.-Nagy's theorem.
\end{proof}

\begin{remark}
	Assume that the inequalities of Proposition \ref{NabokoTh}(ii) are replaced by the stronger ones
	\begin{align*}
		\alpha \|h\|^2 &\leq \liminf_{\varepsilon \to 0^+} \varepsilon \int_{-\infty}^\infty \|(\varepsilon + i \xi - \generator)^{-1} h \|^2 \, d\xi
		\\ &\leq \limsup_{\varepsilon \to 0^+} \varepsilon \int_{-\infty}^\infty \|(\varepsilon + i \xi - \generator)^{-1} h \|^2 \, d\xi
		\leq \beta \|h\|^2, \qquad \varepsilon>0, \, h \in H.
	\end{align*}
	Then the first part of the proof for (ii)$\implies$ (i) in  Proposition \ref{NabokoTh} can be simplified by defining
	\begin{align*}
		\|h\|_{\rm{eq}}^2 := \operatorname{LIM}
			\left[\varepsilon_n \int_{-\infty}^\infty \|(\varepsilon_n + i \xi - \generator)^{-1} h \|^2 \, d\xi \right], \qquad h \in H,
	\end{align*}
	where $\{\varepsilon_n\}_{n\in \mathbb N} \subset (0,\infty)$ is any sequence with $\lim_{n\to \infty} \varepsilon_n = 0$. Then $\|\cdot\|_{{\rm eq}}$ is an equivalent Hilbertian norm on $H$, and following the same steps as in Proposition \ref{NabokoTh} for $n=1$, it can be proved that $\|(z - \generator)^{-1}\|_{{\rm eq}} \leq (\Real z)^{-1}$ for all $z \in \CC$ with $\Real z >0$. Thus, $\mathcal T$ is contractive in $\|\cdot\|_{\rm {eq}},$ and then $\mathcal T$  is similar to a semigroup of contractions.
	Arguing as in the final part of the proof of Theorem  \ref{NabokoTh} (starting from \eqref{belowIneq}), one deduces that $\mathcal T$ is similar to an isometric semigroup.
\end{remark}

\section{Basics on observability and controllability}\label{ObservControlSect}

Here we introduce several basic notions of abstract control theory needed in
Section \ref{ObservSubsect}. Most of them can be found in the books
\cite{staffans2005well} and \cite{tucsnak2009observation}, and we refer to these books
for more details and further comments.

We start with observability of systems governed by $C_0$-semigroups.
Let $(H,\|\cdot\|)$ and $(K, \|\cdot\|_K)$ be two Hilbert spaces and let $\mathcal T = (T(t))_{t\geq0}$ be a $C_0$-semigroup on $H$ with generator $\generator$. Given a bounded operator $C \in \linearOp(\dom(\generator), K)$, where $\dom(\generator)$ is equipped with the graph norm $\|h\|_{\dom(\generator)}^2 = \|h\|^2 + \|\generator h\|^2$, we define the operator $\Psi^C$ from $\dom(\generator)$ to $L^2_{{\rm loc}}((0,\infty),K)$ by
$$
(\Psi^C h)(t) = CT(t)h, \qquad h \in \dom(\generator), \, t>0.
$$
For $\tau>0$, let $\Psi_\tau^C$ be the operator  from $\dom(\generator)$ into  $L^2((0,\tau), K)$ induced by $\Psi^C$, so that $(\Psi_\tau^C h)(t) =  CT(t)h$ for all $h\in \dom(\generator)$ and $t\in (0,\tau)$.

With the above notation, a bounded operator $C \in \linearOp(\dom(\generator), K)$ is said to be a \textit{finite time admissible observation operator} for $\mathcal T$ if for some (hence all) $\tau>0$, the operator $\Psi_\tau^C$ can be continuously extended to a bounded operator from $H$ into $L^2((0,\tau),K)$; or equivalently, there exist $\tau>0$ and $\beta>0$ such that
$$ \int_0^\tau \|C T(t) h\|^2_K \, dt \leq \beta\|h\|^2, \qquad h \in \dom(\generator).
$$
If $C$ is a finite time admissible observation operator for $\mathcal T$, then the pair $(\generator,C)$ is said to be \textit{finite time exactly observable} if there exists $\tau>0$ such that $\Psi_\tau^C$ is bounded from below; or equivalently, there exist $\tau>0$ and $\alpha>0$ such that
$$
\alpha \|h\|^2 \leq \int_0^{\tau} \|C T(t) h\|^2_K \, dt, \qquad h \in \dom(\generator).
$$

The preceding definitions have natural extensions for $\tau=\infty$. Namely, we say that $C$ is an \textit{infinite time admissible observation operator for} $\mathcal T$ if $\Psi^C h \in L^2((0,\infty),K)$ for every $h \in \dom(\generator)$, and $\Psi^C$ can be continuously extended to a bounded operator from $H$ into $L^2((0,\infty),K)$; or equivalently, there exists $\beta>0$ such that
$$
 \int_0^\infty \|C T(t) h\|^2_K \, dt \leq \beta\|h\|^2, \qquad h \in \dom(\generator).
$$

If $C$ is an infinite time admissible observation operator for $\mathcal T$, then the pair $(\generator,C)$ is said to be \textit{infinite time exactly observable} if $\Psi^C$ is bounded from below; that is, there exists $\alpha>0$ such that
$$\alpha \|h\|^2 \leq \int_0^{\infty} \|C T(t) h\|^2_K \, dt, \qquad h \in \dom(\generator).
$$
It is well-known (see e.g. \cite[p. 145 \& Proposition 6.5.2]{tucsnak2009observation}) that if $\mathcal T$ is exponentially stable, then infinite time observability and finite time observability are equivalent. More precisely,
\begin{itemize}
	\item[(i)] $C$ is a finite time admissible observation operator for $\mathcal T$ if and only if $C$ is an infinite time admissible observation operator for $\mathcal T$;
	\item[(ii)] the pair $(\generator,C)$ is finite time exactly observable if and only if the pair $(\generator,C)$ is infinite time exactly observable.
\end{itemize}
Next we proceed with the notion of controllability, which is dual in a sense  to the observability discussed, see Theorem \ref{DualObsTh}. Let $H_{-1}$ be the extrapolated space of $H$, that is, given $\lambda \in \rho(\generator)$, $H_{-1}$ is the completion of $H$ endowed with the Hilbertian norm
$$\|h\|_{-1} := \|(\lambda - \generator)^{-1} h\|, \qquad h \in H.
$$
It is readily seen that the norms $\|\cdot\|_{-1}$ defined above for different $\lambda \in \rho(\generator)$ are equivalent. Thus, as a topological vector space $H_{-1}$ is independent on the choice of $\lambda \in \rho(\generator)$, and $H$ is continuously and densely embedded in $H_{-1}$. Recall also that, given $\|\cdot\|_{-1}$ as above, $\lambda - \generator \in \linearOp(H, H_{-1})$ is a unitary operator. Then $\mathcal T$ induces a $C_0$-semigroup $\mathcal T_{-1}$ on $H_{-1}$ which is unitarily equivalent through $\lambda - \generator$ to $\mathcal T$ on $H$, and such that $H$, as a subspace of $H_{-1}$, coincides with the domain of the generator of $\mathcal T_{-1}$.

Now, let $(K,\|\cdot\|_K)$ be a Hilbert space, and let $B \in \linearOp(K, H_{-1})$ be a bounded operator. For $\tau > 0$,  define the bounded operator $\phi_\tau^B$ from $L_{{\rm loc}}^2((0,\infty),K)$ to $H_{-1}$ as
$$\phi_\tau^B u = \int_0^\tau T_{-1}(\tau-t) Bu(t)\, dt, \qquad u \in L_{{\rm loc}
}^2((0,\infty),K).
$$
Then  $B$ is said to be a \textit{finite time admissible control operator} for $\mathcal T$ if for some (hence for all) $\tau>0$, one has $\operatorname{Ran} \phi_\tau^B \subseteq H$. It follows from the closed graph theorem that, if $B$ is a finite time admissible control operator for $\mathcal T$, then $\phi_\tau^B$ defines a bounded operator from $L^2((0,\infty),K)$ to $H$ for each $\tau>0$, that is, there exists $\beta_\tau>0$ such that
$$\left\| \int_0^\tau T_{-1}(\tau-t) B u(t) \, dt \right\| \leq \beta_\tau \int_0^\tau \|u(t)\|_K^2 \, dt, \qquad u \in L_{{\rm loc}}^2((0,\infty),K),
$$
where the norm in the left hand side is the norm in $H$.

Given a finite time admissible control operator $B$ for $\mathcal T$, one says that the pair $(\generator,B)$ is \textit{finite time exactly controllable} if there exists $\tau$ such that $\operatorname{Ran} \phi_\tau^B = H$, that is, if the operator $\phi_\tau^B$ from $L^2((0,\tau),K)$ to $H$ is surjective.

These control properties can also be defined for $\tau=\infty$ in a natural way, as in the case of observability. Namely, a bounded operator $B \in \linearOp(K, H_{-1})$ is an \textit{infinite time admissible control operator for} $\mathcal T$ if for every $u \in L^2((0,\infty),K)$, the following limit
$$\phi^B u := \lim_{s\to \infty} \int_0^s T_{-1}(t) B u(t) \, dt
$$
exists in $H_{-1}$, and moreover $\phi^B u \in H$ in a way that $\phi^B$ defines a bounded operator from $L^2((0,\infty),K)$ to $H$. That is, there exists $\beta>0$ such that
$$\left\| \int_0^\infty T_{-1}(t) B u(t) \, dt\right\| \leq \beta \int_0^\infty \|u(t)\|_K^2 \, dt, \qquad u \in L^2((0,\infty),K).
$$
Finally, if $B$ is an infinite time admissible control operator for $\mathcal T$, then the pair $(\generator,B)$ is said to be \textit{infinite time exactly controllable} if the operator $\phi^B$ from $L^2((0,\infty),K)$
to $H$ is surjective.

Controllability and observability are dual concepts in the following sense. Recall that the dual space $()\dom(\generator))^\ast$ can be canonically identified with $H_{-1}$ through the anti-linear isomorphism $L \in \linearOp(H_{-1}, ()\dom(\generator))^\ast)$ given by
$$L(h)f = \langle f, h \rangle_H, \qquad f \in \dom(\generator),
$$
for all $h$ in $H,$ and by a density argument, then continuously extended to every $h \in H_{-1}$.
\begin{theorem}\label{DualObsTh}
	Let $H$ be a Hilbert space and let $\mathcal T = (T(t))_{t\geq0}$ be a $C_0$-semigroup on $H$. Let $K$ be another Hilbert space, and let $C \in \linearOp(\dom(\generator), K)$ be a bounded operator. The following holds.
	\begin{itemize}
		\item [(i)] $C$ is a (in)finite time admissible operator for $\mathcal T$ if and only if $C^\ast$ is a (in)finite time admissible control operator for $\mathcal T^\ast$.
		\item [(ii)] Assume $C$ is a (in)finite time admissible operator for $\mathcal T$. Then the pair $(\generator,C)$ is (in)finite time exactly observable if and only if the pair $(\generator^\ast, C^\ast)$ is (in)finite time exactly controllable.
	\end{itemize}
\end{theorem}
For the proof of this result see e.g. \cite[Theorems 4.4.3 \& 11.2.1]{tucsnak2009observation}. 
\vspace{20 pt}

\noindent \textbf{Data availability} No data was used for the research described in the article. 
\vspace{20 pt}

\noindent \textbf{\large Declarations} \bigskip

\noindent \textbf{Competing interests} The authors declare that there is no conflict of interest.

\bibliography{mybib}
\bibliographystyle{amsplain-nodash}

\end{document}